\documentclass{amsart}
\usepackage{graphicx}
\usepackage[all]{xy}
\usepackage{latexsym}
\usepackage{amssymb}
\usepackage{amsmath}
\usepackage{amsthm}
\usepackage{amscd}
\usepackage{fancyhdr}
\usepackage{a4wide}
\usepackage{mathrsfs}

\def\id{{\mbox{1 \hskip -8pt 1}}}
 \newcommand{\lon}{\longrightarrow}
 
 \newcommand{\rar}{\rightarrow}

 \newcommand{\End}{{\mathsf E\mathsf n \mathsf d}}
\newcommand{\p}{{\partial}}
\newcommand{\Id}{{\mathrm I\mathrm d}}

 \newcommand{\Z}{{\mathbb Z}}
 \newcommand{\bS}{{\mathbb S}}

 \newcommand{\N}{{\mathbb N}}
 \newcommand{\K}{{\mathbb K}}
 \newcommand{\ot}{\otimes}

\newcommand{\sB}{{\mathsf  A\mathsf s\mathsf s\mathsf B}}
\newcommand{\sC}{{\mathsf C}}

\newcommand{\sP}{{\mathsf P}}
\newcommand{\sd}{{\mathsf d}}
\newcommand{\hsB}{{\mathsf  A\mathsf s\mathsf s\mathsf B}^+}

\newcommand{\Lieb}{{\mathsf L\mathsf i \mathsf e \mathsf B}}

\newcommand{\DefQ}{{\mathsf D\mathsf e\mathsf f \mathsf Q}}
\newcommand{\Def}{{\mathsf D\mathsf e\mathsf f}}
\newcommand{\Defq}{{\mathsf D\mathsf e\mathsf f \mathsf q}}

 \newcommand{\Beq}{\begin{equation}}
 \newcommand{\Eeq}{\end{equation}}
 \newcommand{\Beqr}{\begin{eqnarray}}
 \newcommand{\Eeqr}{\end{eqnarray}}
 \newcommand{\Beqrn}{\begin{eqnarray*}}
 \newcommand{\Eeqrn}{\end{eqnarray*}}
 \newcommand{\Ba}{\begin{array}}
 \newcommand{\Ea}{\end{array}}
 \newcommand{\Bi}{\begin{itemize}}
 \newcommand{\Ei}{\end{itemize}}
 \newcommand{\Bc}{\begin{center}}
 \newcommand{\Ec}{\end{center}}

 \newcommand{\f}{{\mathcal O}}
 
 \newcommand{\cB}{{\mathcal B}}

 \newcommand{\cE}{{\mathcal E}}
 \newcommand{\cF}{{\mathcal F}}

 \newcommand{\cK}{{\mathcal K}}
 \newcommand{\caL}{{\mathcal L}}
 
 \newcommand{\cP}{{\mathcal P}}

 \newcommand{\cV}{{\mathcal V}}

 \newcommand{\fg}{{\mathfrak g}}

 \newcommand{\fl}{{\mathfrak l}}

\newcommand{\fs}{{\mathfrak s}}

\newcommand{\ii}{{\mathfrak i}}
 \newcommand{\al}{\alpha}
 \newcommand{\be}{\beta}
 \newcommand{\ga}{\gamma}
 
 \newcommand{\Ga}{\Gamma}
 
 \newcommand{\var}{\varepsilon}

\newcommand{\hdelta}{{\delta^+}}

 \newcommand{\Hom}{{\mathrm H\mathrm o\mathrm m}}
 
 \newcommand{\sip}{\smallskip}
 \newcommand{\bip}{\bigskip}

\swapnumbers
\newtheorem{theorem}{Theorem}[subsection]

\newtheorem{lemma}[theorem]{Lemma}

\newtheorem{prop-def}[theorem]{Proposition-definition}

\newtheorem{f-theorem}{Formality Theorem}[section]
\newtheorem{main-theorem}{Main~Theorem}[section]
\newtheorem{section-theorem}{Theorem}[section]
\newtheorem{section-corollary}{Corollary}[section]

\theoremstyle{definition}

\newtheorem{fact-me}{Fact \cite{Me1}}[subsection]

\begin{document}



 \title{Formality theorem for
 quantizations of Lie bialgebras
 }
 \author{ S.A.\ Merkulov}
 \thanks{This work was
 partially supported by the G\"oran Gustafsson foundation.}
\address{
 Department of Mathematics, Stockholm University, Sweden and Mathematics Research Unit, Luxembourg University,  Grand Duchy of Luxembourg (present address) }
\email{sergei.merkulov@uni.lu}
 \date{}
 \begin{abstract}
 Using the theory of props we prove a formality theorem associated with universal
 quantizations of  Lie bialgebras.
 \bip



 \end{abstract}
 \maketitle

\section{Introduction}

\subsection{Two complexes related by a formality map}
 Let $V$ be a $\Z$-graded 
vector space over a
field $\K$ and $\f_V:= {\odot^{\bullet}} V= \oplus_{n\geq 0} \odot^n V$
the
graded commutative and
cocommutative bialgebra of polynomial functions on $V^*$. The
Gerstenhaber-Schack complex \cite{GS} 
\Beq\label{1: GS complex}
\left(\fg\fs(\f_V,\f_V)=\prod_{m,n\geq 1}\Hom(\f_V^{\ot m},
\f_V^{\ot n})[2-m-n], \
 d_{\fg\fs}\right),
\Eeq
has a $L_\infty$-algebra structure,
\[
\left\{\mu_n: \odot^n\fg\fs(\f_V,\f_V)\rar \fg\fs(\f_V,\f_V)[2-n]\right\}_{n\geq 1}\ \
\mbox{with}\ \ \mu_1=d_{\fg\fs},
\]
which controls deformations of the standard bialgebra structure on $\f_V$ \cite{MV} (an explicit formula for the differential $d_{\fg\fs}$ is given in \S 3.4.1 below).
This $L_\infty$-algebra
depends on the choice of a minimal resolution of the properad of bialgebras, but its isomorphism class
is defined canonically. Existence (and non-uniqueness) of such minimal resolutions was proven by Martin Markl in \cite{Ma}.

On the other hand, the completed graded commutative algebra,
\Beq\label{1: l_V}
\fl_V:=  {\prod_{m,n\geq 1}}\odot^m ( V[-1])\, {\ot}\,  \odot^{n}( V^*[-1]),
\Eeq
has a natural degree $-2$ Poisson structure, $\{\, ,\, \}: \fl_V\ot \fl_V\rar \fl_V[-2]$, which makes $\fl_V$ into a 3-algebra and which is
given on generators by
\[
\{sv, sw\}=0,\ \ \{s\al, s\be\}=0, \ \ \{s\al, sv\}=<\al,v>, \ \ \forall v,w\in V, \al,\be\in V^*.
\]
where $s: V\rar V[-1]$ and $s: V^*\rar V^*[-1]$ are natural isomorphisms. The Maurer-Cartan elements
of $(\fl_V,\ [\ ,\ \})$ are precisely strongly homotopy bialgebra structures in $V$ (see Corollary 5.1 in \cite{Me1}).

\bip

It was proven by Etingof and Kazhdan \cite{EK} that there exists a universal quantization of (possibly,
infinite-dimensional) Lie bialgebras. Any such a universal quantization morphism is highly non-trivial --- it depends on the choice of a Drinfeld associator \cite{Dr}. The main result of this paper is a proof of the following

\setcounter{f-theorem}{1}
\begin{f-theorem}\label{1.2}
{ Every universal quantization of Lie bialgebras lifts to a $L_\infty$ quasi-isomorphism},
\Beq\label{1:F}
F: \left(\fl_V[2], \{\, , \, \}\right) \lon  \left(\fg\fs(\f_V,\f_V), \mu_\bullet\right).
\Eeq
\end{f-theorem}

\sip
\noindent As the vector space $\fl_V[2]$ equals the cohomology of the
Gerstenhaber-Schack complex
$
\left(\fg\fs^\bullet(\f_V), d_{\fg\fs}\right)
$, Theorem~\ref{1.2} implies formality of the $L_\infty$ algebra
$ \left(\fg\fs(\f_V,\f_V), \mu_\bullet\right)$. Our proof of Theorem \ref{1.2} shows that a  lifting of a universal Etingof-Kazhdan
quantization morphism  to a formality map $F$ can be constructed {\em inductively}, i.e.\ the lifting itself does {\em not}\,  involve transcendental methods which seem to be unavoidable in the construction of Drinfeld associators and hence of the universal quantizations of Lie bialgebras; our proof is based, therefore, on a ``trivial" procedure provided one reformulates the problem of constructing $F$ in terms of differential graded props, and this reformulation occupies most of the present paper.

\setcounter{section-corollary}{2}
\begin{section-corollary}\label{1.3}
{ Every strongly homotopy Lie bialgebra structure, $\ga$, on a graded vector space $V$ can be deformation quantized, i.e.\
there exists a strongly homotopy bialgebra structure, $\Gamma(\ga,\hbar)$, on $\odot^\bullet V[[\hbar]]$
depending on a formal parameter $\hbar$ such that $\Gamma|_{\hbar=0}$ coincides with
the standard graded (co)commutative bialgebra structure in $\odot^\bullet V$
while $\frac{d \Ga}{d\hbar}|_{\hbar=0}$ induces in $V$ the original homotopy Lie bialgebra structure $\ga$.}
\end{section-corollary}

\setcounter{subsection}{3}
\subsection{An outline of the proof of the Formality Theorem}
Our main technical tool in proving Theorem~\ref{1.2} is the theory of differential graded (dg, for short) props which allows us to make the idea of {\em universality}\, of  quantizations of Lie bialgebras rigorous (a very different approach to prop interpretation of deformation quantizations was developed in \cite{EE}, while a very similar approach was used by the author in \cite{Me2} to describe
deformation quantization of wheeled Poisson structures).
We refer
 to \cite{EE, Ma3,Va}
for an introduction into the theory of props and properads, their minimal resolutions and
representations.
\sip

There is a prop $\sB$ whose representations, $\sB\rar \mathsf{End}_W$, in a dg vector space $W$ are precisely bialgebra structures in $W$. As $\f_V=\odot^\bullet V$ has a canonical bialgebra structure, there is associated a canonical morphism of props,
$$
\rho_0:  \sB\lon \mathsf{End}_{\f_V}.
$$
The standard machinery developed in \cite{MV} gives us a deformation complex of the map
$\rho_0$,
$$
\Def(\sB_\infty\stackrel{\rho_0}{\lon}\mathsf{End}_{\f_V})\simeq \prod_{m,n\geq 1\atop m+n\geq 3}\Hom(\f_V^{\ot m},
\f_V^{\ot n})[2-m-n],
$$
which comes equipped with a filtered $L_\infty$ algebra structure (depending on the choice
of a minimal resolution $\sB_\infty$ of $\sB$) whose Maurer-Cartan elements $\ga$
are in one-to-one correspondence with morphisms of dg props $\rho_0+\ga:  \sB\lon \mathsf{End}_{\f_V}$, that is, with strongly
homotopy bialgebra structures in ${\f_V}$. The differential $\mu_1$ of this $\caL_\infty$ structure is determined by the map $\rho_0$ and was proven in \cite{MV} to coincide precisely with the differential $d_{\fg\fs}$ introduced by Gerstenhaber and Schack in \cite{GS}. More details can be found in \S 3.

\sip

Similarly, there is a prop(erad) $\Lieb$ of Lie bialgebras which is Koszul and hence admits a simple minimal resolution
$\Lieb_\infty$, see \S 2 for details. Deformation complex of the zero morphism
$$
\Def(\Lieb_\infty\stackrel{0}{\lon} \mathsf{End}_V)\simeq  {\prod_{m,n\geq 1\atop m+n\geq 3}}\odot^m ( V[-1])\, {\ot}\,  \odot^{n}( V^*[-1])
$$
comes equipped, in accordance with \cite{MV}, with a graded Lie algebra structure
whose Maurer-Cartan elements $\nu$ are in one-to-one correspondence with morphisms of dg props $\nu:  \sB\lon \mathsf{End}_{V}$, that is, with strongly
homotopy bialgebra structures in ${V}$. The Lie brackets are precisely the ones $\{\ ,\ \}$ described in \S 1.1.

\sip

The idea of the proof is to relate somehow the dg props $\sB_\infty$ and $\Lieb_\infty$
to each other rather than to work with their representations. For this we have to resolve
a number of problems.

\subsubsection{}  The complex $\Def(\sB_\infty\stackrel{\rho_0}{\lon}\mathsf{End}_{\f_V})$
does not coincide with the Gerstenhaber-Schack complex (\ref{1: GS complex}); more unpleasantly,  its cohomology with respect to $d_{\fg\fs}$ does not equal the space $\Def(\Lieb_\infty\stackrel{0}{\lon} \mathsf{End}_V)$. 
 This problem can, however, be ressolved easily: there is an  functor, $F^+: \cP\rar \cP^+$
in the category of dg props introduced in \S 4.1 below: for any dg prop $\cP$, the associated prop $\cP^+$ can be uniquely characterized in terms of $\cP$ as follows --- there is a 1-1 correspondence between
representations of $\cP^+$ in a graded vector space $V$ equipped with the zero differential,
and representations of $\cP$ in the same vector space $V$ but equipped with an {\em arbitrary} (i.e.\ not fixed a priori) differential. Then we recover the desired complexes using the standard
deformation theory of morphisms of props,
$$
\Def(\sB_\infty^+\stackrel{\rho_0}{\rar}\mathsf{End}_{\f_V})=\fg\fs(\f_V,\f_V)\ \ \ \mbox{and}\ \ \
\Def(\Lieb_\infty^{_+}\stackrel{0}{\rar} \mathsf{End}_V)=\fl_V.
$$

\subsubsection{}
 The props $\sB_\infty^+$ and $\Lieb_\infty^+$ can not be related directly as the first one
has natural representations in $\f_V$ while the second one in $V$. We solve this problem in two steps:
\Bi
\item[-] first we notice that there is a sub-prop, $\mathsf{End^{poly}_{\f_{\mathit V}}}$, of the endomorphism prop $\mathsf{End}_{\f_V}$
spanned by so-called {\em polydifferential}\, operators $\Phi: \f_V^{\ot m}\rar \f_V^{\ot n}$;
 the image of the standard representation $\rho_0:\sB^+\rar \mathsf{End}_{\f_V}$
lands in   $\mathsf{End^{poly}_{\f_{\mathit V}}}$. It was shown in \cite{Me3} that the associated
$L_\infty$ (sub)algebra
$$
\mathfrak{poly}({\f_V}, \f_V):=\Def(\sB_\infty^+ \stackrel{\rho_0}{\rar} \mathsf{End^{poly}_{\f_{\mathit V}}})
\subset \fg\fs(\f_V,\f_V)
$$
is quasi-isomorphic to $\fg\fs(\f_V,\f_V)$. Thus to prove the Main Theorem it is enough to show
existence of a $\caL_\infty$ quasi-isomorphism
\Beq\label{1: Formality map F poly}
F: \fl_V[2] \lon  \mathfrak{poly}(\f_V,\f_V).
\Eeq
\item[-] There is a dg free prop $\DefQ^+$ uniquely characterized by the following property:
there is a one-to-one correspondence between polydifferential representations of $\sB_\infty^+$ in $\f_V$ and ordinary representations of $\DefQ^+$ in the vector space $V$. Its generators and the differential are described in \S 4.
\Ei

\subsubsection{} We construct a highly non-trivial continuous morphism (in fact, a quasi-isomorphism) of topological dg props
\Beq\label{1: map of props cF}
\cF^+: \DefQ^+ \lon  \widehat{\Lieb}_\infty^{_+}
\Eeq
where  $\widehat{\Lieb}_\infty^{_+}$ is the vertex $+$ genus completion\footnote{The fact that the morphism $\cF^+$ takes values in the {\em completed}\, prop is quite expected. Indeed, all the homogeneous components
 $F_k: \wedge^k(\fl_V[2]) \rar \mathsf{poly}(\f_V,\f_V)$ of (\ref{1: Formality map F poly}) are, in general, non-trivial, and every such a component $F_k$ corresponds precisely to the $k$-vertex summands in the values of $\cF^+$ on the generators of
 $\DefQ^+$.} of the prop  ${\Lieb}_\infty^{_+}$,
as follows. It was shown in \cite{EK} that, for any choice of the Drinfeld associator, there
is a universal quantization of Lie bialgebras. A differential graded version of this result was settled in the Appendix \S 8 of \cite{GH}. In our context, existence
of a universal quantization of an arbitrary (possibly, infinite-dimensional) dg Lie bialgebra
is equivalent to existence of a morphism of dg props
$$
\cE\cK^+: \DefQ^+ \lon  \widehat{\Lieb}^{_+}
$$
satisfying a certain non-triviality condition (see \S 5 for details.). Then, using the fact that the natural projection  $\widehat{\Lieb}_\infty^{_+}\rar \widehat{\Lieb}^{_+}$ is a quasi-isomorphism, the required morphism (\ref{1: map of props cF}) is inductively constructed as a lifting of
$\cE\cK^+$ making the diagram
\[
 \xymatrix{
  & \widehat{\Lieb}^{_+}_\infty
  \ar[d]^{qis} \\
 \DefQ^+ \ar[ur]^{{\cF^+}}\ar[r]_{\cE\cK^+} &
 \widehat{\Lieb}^{_+}
 }
\]
commutative. Full details are given in \S 5.4.

\subsubsection{} The final step is to use the morphism (\ref{1: map of props cF}) in order to show a proof of the Formality Theorem 2.1. This is done in two different ways in \S 5.

\subsection{Some notation}
For an $\bS$-bimodule $E=\{E(m,n)\}_{m,n\geq 1}$ the associated free prop is denoted
by $\Gamma\langle E \rangle$. The endomorphism prop of a graded vector space $V$ is denoted by
$\End_V$. The one-dimensional sign representation of the permutation group $\bS_n$
is denoted by $sgn_n$ while the trivial representation by $\id_n$.
All our differentials have degree $+1$.
 If $V=\oplus_{i\in \Z} V^i$ is a graded vector space, then
$V[k]$ is a graded vector space with $V[k]^i:=V^{i+k}$; for $a\in V^i\subset V$ we write $|a|=i$.
All our graphs are directed with flow implicitly assumed to go from bottom to top.
We work throughout in the category of $\Z$-{\em graded}\, vector spaces over a field $\K$ of characteristic 0 with its standard Koszul sign conventions; for example, the value of an element $f_1\ot f_2\in V^*\ot V^*$ on an element $a\ot b\in V\ot V$ is equal to
$(-1)^{|f_2||a|} f_1(a)f_2(b)$, where $V^*:=\Hom(V,\K)$.


\bip

\section{Prop of strongly homotopy Lie bialgebras}

\subsection{Lie bialgebras.}\label{2.1}
 A {\em  Lie bialgebra}\,
  is, by definition \cite{D}, a graded vector space $V$ together with two linear maps,
$$
\Ba{rccc}
\vartriangle: &  V & \lon & \wedge^2 V \\
       & a    & \lon & \sum a_{1}\wedge a_{2}
\Ea
\ \ \ \ \  , \ \ \ \
\Ba{rccc}
[\, , \, ]: & \wedge^2 V& \lon & V \\
       & a\ot b    & \lon & [a, b]
\Ea
$$
satisfying,
\Bi
\item[(i)]  the Jacobi identity:
$[a, [b, c]=[[a,b], c]] +
(-1)^{|b||a|}[b,[a, c]]$;
\item[(ii)] the co-Jacobi identity: $(\vartriangle\ot\Id)\vartriangle a
+ \tau (\vartriangle\ot\Id)\vartriangle a+ \tau^2 (\vartriangle\ot\Id)\vartriangle a =0$,
where $\tau$ is the cyclic permutation $(123)$ represented
naturally in $V\ot V \ot V$;
\item[(iii)]  the Drinfeld compatibility condition:
 $\vartriangle [a, b]=\sum a_1\wedge [a_2, b] - (-1)^{|a_1||a_2|} a_2\wedge
[a_1, b]
+  [a, b_1]\wedge b_2 - (-1)^{|b_1||b_2|}[a, b_2]\wedge b_1$.
\Ei
for any $a,b, c\in V$.

\subsection{Prop of Lie bialgebras}\label{Lieb}
 It is easy to construct a prop, $\Lieb$, whose
representations,
$$
\rho: \Lieb \lon \End\langle V\rangle,
$$
in a graded vector space $V$ are in one-to-one correspondence with Lie bialgebra structures
in $V$. With an association in mind,
$$
\vartriangle \leftrightarrow
 \begin{xy}
 <0mm,-0.55mm>*{};<0mm,-2.5mm>*{}**@{-},
 <0.5mm,0.5mm>*{};<2.2mm,2.2mm>*{}**@{-},
 <-0.48mm,0.48mm>*{};<-2.2mm,2.2mm>*{}**@{-},
 <0mm,0mm>*{\circ};<0mm,0mm>*{}**@{},
 \end{xy}
 \ \ \  , \ \ \
[\, , \, ] \leftrightarrow
 \begin{xy}
 <0mm,0.66mm>*{};<0mm,3mm>*{}**@{-},
 <0.39mm,-0.39mm>*{};<2.2mm,-2.2mm>*{}**@{-},
 <-0.35mm,-0.35mm>*{};<-2.2mm,-2.2mm>*{}**@{-},
 <0mm,0mm>*{\circ};<0mm,0mm>*{}**@{},
 \end{xy},
$$
one can define it as a quotient,
$$
\Lieb:= {\Gamma\langle L\rangle}/(R)
$$
of the free prop, $\Gamma\langle L \rangle$, generated by the $\bS$-bimodule $L=\{L(m,n)\}$,
\[
L(m,n):=\left\{
\Ba{rr}
sgn_2\ot \id_1\equiv\mbox{span}\left\langle
\begin{xy}
 <0mm,-0.55mm>*{};<0mm,-2.5mm>*{}**@{-},
 <0.5mm,0.5mm>*{};<2.2mm,2.2mm>*{}**@{-},
 <-0.48mm,0.48mm>*{};<-2.2mm,2.2mm>*{}**@{-},
 <0mm,0mm>*{\circ};<0mm,0mm>*{}**@{},
 <0mm,-0.55mm>*{};<0mm,-3.8mm>*{_1}**@{},
 <0.5mm,0.5mm>*{};<2.7mm,2.8mm>*{^2}**@{},
 <-0.48mm,0.48mm>*{};<-2.7mm,2.8mm>*{^1}**@{},
 \end{xy}
=-
\begin{xy}
 <0mm,-0.55mm>*{};<0mm,-2.5mm>*{}**@{-},
 <0.5mm,0.5mm>*{};<2.2mm,2.2mm>*{}**@{-},
 <-0.48mm,0.48mm>*{};<-2.2mm,2.2mm>*{}**@{-},
 <0mm,0mm>*{\circ};<0mm,0mm>*{}**@{},
 <0mm,-0.55mm>*{};<0mm,-3.8mm>*{_1}**@{},
 <0.5mm,0.5mm>*{};<2.7mm,2.8mm>*{^1}**@{},
 <-0.48mm,0.48mm>*{};<-2.7mm,2.8mm>*{^2}**@{},
 \end{xy}
   \right\rangle  & \mbox{if}\ m=2, n=1,\vspace{3mm}\\
\id_1\ot sgn_2\equiv
\mbox{span}\left\langle
\begin{xy}
 <0mm,0.66mm>*{};<0mm,3mm>*{}**@{-},
 <0.39mm,-0.39mm>*{};<2.2mm,-2.2mm>*{}**@{-},
 <-0.35mm,-0.35mm>*{};<-2.2mm,-2.2mm>*{}**@{-},
 <0mm,0mm>*{\circ};<0mm,0mm>*{}**@{},
   <0mm,0.66mm>*{};<0mm,3.4mm>*{^1}**@{},
   <0.39mm,-0.39mm>*{};<2.9mm,-4mm>*{^2}**@{},
   <-0.35mm,-0.35mm>*{};<-2.8mm,-4mm>*{^1}**@{},
\end{xy}=-
\begin{xy}
 <0mm,0.66mm>*{};<0mm,3mm>*{}**@{-},
 <0.39mm,-0.39mm>*{};<2.2mm,-2.2mm>*{}**@{-},
 <-0.35mm,-0.35mm>*{};<-2.2mm,-2.2mm>*{}**@{-},
 <0mm,0mm>*{\circ};<0mm,0mm>*{}**@{},
   <0mm,0.66mm>*{};<0mm,3.4mm>*{^1}**@{},
   <0.39mm,-0.39mm>*{};<2.9mm,-4mm>*{^1}**@{},
   <-0.35mm,-0.35mm>*{};<-2.8mm,-4mm>*{^2}**@{},
\end{xy}
\right\rangle
\ & \mbox{if}\ m=1, n=2, \vspace{3mm}\\
0 & \mbox{otherwise}
\Ea
\right.
\]
modulo the ideal generated by relations
\Beq\label{R}
R:\left\{
\Ba{r}
\begin{xy}
 <0mm,0mm>*{\circ};<0mm,0mm>*{}**@{},
 <0mm,-0.49mm>*{};<0mm,-3.0mm>*{}**@{-},
 <0.49mm,0.49mm>*{};<1.9mm,1.9mm>*{}**@{-},
 <-0.5mm,0.5mm>*{};<-1.9mm,1.9mm>*{}**@{-},
 <-2.3mm,2.3mm>*{\circ};<-2.3mm,2.3mm>*{}**@{},
 <-1.8mm,2.8mm>*{};<0mm,4.9mm>*{}**@{-},
 <-2.8mm,2.9mm>*{};<-4.6mm,4.9mm>*{}**@{-},
   <0.49mm,0.49mm>*{};<2.7mm,2.3mm>*{^3}**@{},
   <-1.8mm,2.8mm>*{};<0.4mm,5.3mm>*{^2}**@{},
   <-2.8mm,2.9mm>*{};<-5.1mm,5.3mm>*{^1}**@{},
 \end{xy}
\ + \
\begin{xy}
 <0mm,0mm>*{\circ};<0mm,0mm>*{}**@{},
 <0mm,-0.49mm>*{};<0mm,-3.0mm>*{}**@{-},
 <0.49mm,0.49mm>*{};<1.9mm,1.9mm>*{}**@{-},
 <-0.5mm,0.5mm>*{};<-1.9mm,1.9mm>*{}**@{-},
 <-2.3mm,2.3mm>*{\circ};<-2.3mm,2.3mm>*{}**@{},
 <-1.8mm,2.8mm>*{};<0mm,4.9mm>*{}**@{-},
 <-2.8mm,2.9mm>*{};<-4.6mm,4.9mm>*{}**@{-},
   <0.49mm,0.49mm>*{};<2.7mm,2.3mm>*{^2}**@{},
   <-1.8mm,2.8mm>*{};<0.4mm,5.3mm>*{^1}**@{},
   <-2.8mm,2.9mm>*{};<-5.1mm,5.3mm>*{^3}**@{},
 \end{xy}
\ + \
\begin{xy}
 <0mm,0mm>*{\circ};<0mm,0mm>*{}**@{},
 <0mm,-0.49mm>*{};<0mm,-3.0mm>*{}**@{-},
 <0.49mm,0.49mm>*{};<1.9mm,1.9mm>*{}**@{-},
 <-0.5mm,0.5mm>*{};<-1.9mm,1.9mm>*{}**@{-},
 <-2.3mm,2.3mm>*{\circ};<-2.3mm,2.3mm>*{}**@{},
 <-1.8mm,2.8mm>*{};<0mm,4.9mm>*{}**@{-},
 <-2.8mm,2.9mm>*{};<-4.6mm,4.9mm>*{}**@{-},
   <0.49mm,0.49mm>*{};<2.7mm,2.3mm>*{^1}**@{},
   <-1.8mm,2.8mm>*{};<0.4mm,5.3mm>*{^3}**@{},
   <-2.8mm,2.9mm>*{};<-5.1mm,5.3mm>*{^2}**@{},
 \end{xy}
\ \ \ \in \Gamma \langle E\rangle(3,1) \vspace{3mm}\\
 \begin{xy}
 <0mm,0mm>*{\circ};<0mm,0mm>*{}**@{},
 <0mm,0.69mm>*{};<0mm,3.0mm>*{}**@{-},
 <0.39mm,-0.39mm>*{};<2.4mm,-2.4mm>*{}**@{-},
 <-0.35mm,-0.35mm>*{};<-1.9mm,-1.9mm>*{}**@{-},
 <-2.4mm,-2.4mm>*{\circ};<-2.4mm,-2.4mm>*{}**@{},
 <-2.0mm,-2.8mm>*{};<0mm,-4.9mm>*{}**@{-},
 <-2.8mm,-2.9mm>*{};<-4.7mm,-4.9mm>*{}**@{-},
    <0.39mm,-0.39mm>*{};<3.3mm,-4.0mm>*{^3}**@{},
    <-2.0mm,-2.8mm>*{};<0.5mm,-6.7mm>*{^2}**@{},
    <-2.8mm,-2.9mm>*{};<-5.2mm,-6.7mm>*{^1}**@{},
 \end{xy}
\ + \
 \begin{xy}
 <0mm,0mm>*{\circ};<0mm,0mm>*{}**@{},
 <0mm,0.69mm>*{};<0mm,3.0mm>*{}**@{-},
 <0.39mm,-0.39mm>*{};<2.4mm,-2.4mm>*{}**@{-},
 <-0.35mm,-0.35mm>*{};<-1.9mm,-1.9mm>*{}**@{-},
 <-2.4mm,-2.4mm>*{\circ};<-2.4mm,-2.4mm>*{}**@{},
 <-2.0mm,-2.8mm>*{};<0mm,-4.9mm>*{}**@{-},
 <-2.8mm,-2.9mm>*{};<-4.7mm,-4.9mm>*{}**@{-},
    <0.39mm,-0.39mm>*{};<3.3mm,-4.0mm>*{^2}**@{},
    <-2.0mm,-2.8mm>*{};<0.5mm,-6.7mm>*{^1}**@{},
    <-2.8mm,-2.9mm>*{};<-5.2mm,-6.7mm>*{^3}**@{},
 \end{xy}
\ + \
 \begin{xy}
 <0mm,0mm>*{\circ};<0mm,0mm>*{}**@{},
 <0mm,0.69mm>*{};<0mm,3.0mm>*{}**@{-},
 <0.39mm,-0.39mm>*{};<2.4mm,-2.4mm>*{}**@{-},
 <-0.35mm,-0.35mm>*{};<-1.9mm,-1.9mm>*{}**@{-},
 <-2.4mm,-2.4mm>*{\circ};<-2.4mm,-2.4mm>*{}**@{},
 <-2.0mm,-2.8mm>*{};<0mm,-4.9mm>*{}**@{-},
 <-2.8mm,-2.9mm>*{};<-4.7mm,-4.9mm>*{}**@{-},
    <0.39mm,-0.39mm>*{};<3.3mm,-4.0mm>*{^1}**@{},
    <-2.0mm,-2.8mm>*{};<0.5mm,-6.7mm>*{^3}**@{},
    <-2.8mm,-2.9mm>*{};<-5.2mm,-6.7mm>*{^2}**@{},
 \end{xy}
\  \ \ \in  \Gamma \langle E\rangle(1,3) \vspace{2mm} \\
 \begin{xy}
 <0mm,2.47mm>*{};<0mm,0.12mm>*{}**@{-},
 <0.5mm,3.5mm>*{};<2.2mm,5.2mm>*{}**@{-},
 <-0.48mm,3.48mm>*{};<-2.2mm,5.2mm>*{}**@{-},
 <0mm,3mm>*{\circ};<0mm,3mm>*{}**@{},
  <0mm,-0.8mm>*{\circ};<0mm,-0.8mm>*{}**@{},
<-0.39mm,-1.2mm>*{};<-2.2mm,-3.5mm>*{}**@{-},
 <0.39mm,-1.2mm>*{};<2.2mm,-3.5mm>*{}**@{-},
     <0.5mm,3.5mm>*{};<2.8mm,5.7mm>*{^2}**@{},
     <-0.48mm,3.48mm>*{};<-2.8mm,5.7mm>*{^1}**@{},
   <0mm,-0.8mm>*{};<-2.7mm,-5.2mm>*{^1}**@{},
   <0mm,-0.8mm>*{};<2.7mm,-5.2mm>*{^2}**@{},
\end{xy}
\  - \
\begin{xy}
 <0mm,-1.3mm>*{};<0mm,-3.5mm>*{}**@{-},
 <0.38mm,-0.2mm>*{};<2.0mm,2.0mm>*{}**@{-},
 <-0.38mm,-0.2mm>*{};<-2.2mm,2.2mm>*{}**@{-},
<0mm,-0.8mm>*{\circ};<0mm,0.8mm>*{}**@{},
 <2.4mm,2.4mm>*{\circ};<2.4mm,2.4mm>*{}**@{},
 <2.77mm,2.0mm>*{};<4.4mm,-0.8mm>*{}**@{-},
 <2.4mm,3mm>*{};<2.4mm,5.2mm>*{}**@{-},
     <0mm,-1.3mm>*{};<0mm,-5.3mm>*{^1}**@{},
     <2.5mm,2.3mm>*{};<5.1mm,-2.6mm>*{^2}**@{},
    <2.4mm,2.5mm>*{};<2.4mm,5.7mm>*{^2}**@{},
    <-0.38mm,-0.2mm>*{};<-2.8mm,2.5mm>*{^1}**@{},
    \end{xy}
\  + \
\begin{xy}
 <0mm,-1.3mm>*{};<0mm,-3.5mm>*{}**@{-},
 <0.38mm,-0.2mm>*{};<2.0mm,2.0mm>*{}**@{-},
 <-0.38mm,-0.2mm>*{};<-2.2mm,2.2mm>*{}**@{-},
<0mm,-0.8mm>*{\circ};<0mm,0.8mm>*{}**@{},
 <2.4mm,2.4mm>*{\circ};<2.4mm,2.4mm>*{}**@{},
 <2.77mm,2.0mm>*{};<4.4mm,-0.8mm>*{}**@{-},
 <2.4mm,3mm>*{};<2.4mm,5.2mm>*{}**@{-},
     <0mm,-1.3mm>*{};<0mm,-5.3mm>*{^2}**@{},
     <2.5mm,2.3mm>*{};<5.1mm,-2.6mm>*{^1}**@{},
    <2.4mm,2.5mm>*{};<2.4mm,5.7mm>*{^2}**@{},
    <-0.38mm,-0.2mm>*{};<-2.8mm,2.5mm>*{^1}**@{},
    \end{xy}
\  - \
\begin{xy}
 <0mm,-1.3mm>*{};<0mm,-3.5mm>*{}**@{-},
 <0.38mm,-0.2mm>*{};<2.0mm,2.0mm>*{}**@{-},
 <-0.38mm,-0.2mm>*{};<-2.2mm,2.2mm>*{}**@{-},
<0mm,-0.8mm>*{\circ};<0mm,0.8mm>*{}**@{},
 <2.4mm,2.4mm>*{\circ};<2.4mm,2.4mm>*{}**@{},
 <2.77mm,2.0mm>*{};<4.4mm,-0.8mm>*{}**@{-},
 <2.4mm,3mm>*{};<2.4mm,5.2mm>*{}**@{-},
     <0mm,-1.3mm>*{};<0mm,-5.3mm>*{^2}**@{},
     <2.5mm,2.3mm>*{};<5.1mm,-2.6mm>*{^1}**@{},
    <2.4mm,2.5mm>*{};<2.4mm,5.7mm>*{^1}**@{},
    <-0.38mm,-0.2mm>*{};<-2.8mm,2.5mm>*{^2}**@{},
    \end{xy}
\ + \
\begin{xy}
 <0mm,-1.3mm>*{};<0mm,-3.5mm>*{}**@{-},
 <0.38mm,-0.2mm>*{};<2.0mm,2.0mm>*{}**@{-},
 <-0.38mm,-0.2mm>*{};<-2.2mm,2.2mm>*{}**@{-},
<0mm,-0.8mm>*{\circ};<0mm,0.8mm>*{}**@{},
 <2.4mm,2.4mm>*{\circ};<2.4mm,2.4mm>*{}**@{},
 <2.77mm,2.0mm>*{};<4.4mm,-0.8mm>*{}**@{-},
 <2.4mm,3mm>*{};<2.4mm,5.2mm>*{}**@{-},
     <0mm,-1.3mm>*{};<0mm,-5.3mm>*{^1}**@{},
     <2.5mm,2.3mm>*{};<5.1mm,-2.6mm>*{^2}**@{},
    <2.4mm,2.5mm>*{};<2.4mm,5.7mm>*{^1}**@{},
    <-0.38mm,-0.2mm>*{};<-2.8mm,2.5mm>*{^2}**@{},
    \end{xy}
\ \ \ \in  \Gamma \langle E\rangle(2,2)
\Ea
\right.
\Eeq
As the ideal is generated by two-vertex graphs, the properad behind $\Lieb$ is quadratic.

\subsection{Minimal resolution of $\Lieb$} It was shown in \cite{Gan} that  $\Lieb$ is Koszul at the dioperadic level. The extension of this result to the level of props turned out to be highly non-trivial \cite{Ko0,MaVo}. The minimal prop resolution, $\Lieb_\infty$, of $\Lieb$ is a dg free prop,
$$
\Lieb_\infty=\Gamma \langle\mathsf L\rangle,
$$
generated by the $\bS$--bimodule $\mathsf L=\{\mathsf L(m,n)\}_{m,n\geq 1, m+n\geq 3}$,
\Beq\label{L}
\mathsf L(m,n):= sgn_m\ot sgn_n[m+n-3]=\mbox{span}\left\langle
\resizebox{14mm}{!}{\begin{xy}
 <0mm,0mm>*{\circ};<0mm,0mm>*{}**@{},
 <-0.6mm,0.44mm>*{};<-8mm,5mm>*{}**@{-},
 <-0.4mm,0.7mm>*{};<-4.5mm,5mm>*{}**@{-},
 <0mm,0mm>*{};<-1mm,5mm>*{\ldots}**@{},
 <0.4mm,0.7mm>*{};<4.5mm,5mm>*{}**@{-},
 <0.6mm,0.44mm>*{};<8mm,5mm>*{}**@{-},
   <0mm,0mm>*{};<-8.5mm,5.5mm>*{^1}**@{},
   <0mm,0mm>*{};<-5mm,5.5mm>*{^2}**@{},
   <0mm,0mm>*{};<4.5mm,5.5mm>*{^{m\hspace{-0.5mm}-\hspace{-0.5mm}1}}**@{},
   <0mm,0mm>*{};<9.0mm,5.5mm>*{^m}**@{},
 <-0.6mm,-0.44mm>*{};<-8mm,-5mm>*{}**@{-},
 <-0.4mm,-0.7mm>*{};<-4.5mm,-5mm>*{}**@{-},
 <0mm,0mm>*{};<-1mm,-5mm>*{\ldots}**@{},
 <0.4mm,-0.7mm>*{};<4.5mm,-5mm>*{}**@{-},
 <0.6mm,-0.44mm>*{};<8mm,-5mm>*{}**@{-},
   <0mm,0mm>*{};<-8.5mm,-6.9mm>*{^1}**@{},
   <0mm,0mm>*{};<-5mm,-6.9mm>*{^2}**@{},
   <0mm,0mm>*{};<4.5mm,-6.9mm>*{^{n\hspace{-0.5mm}-\hspace{-0.5mm}1}}**@{},
   <0mm,0mm>*{};<9.0mm,-6.9mm>*{^n}**@{},
 \end{xy}}
\right\rangle,
\Eeq
and with the differential of the form
$$
\delta
\resizebox{14mm}{!}{\begin{xy}
 <0mm,0mm>*{\circ};<0mm,0mm>*{}**@{},
 <-0.6mm,0.44mm>*{};<-8mm,5mm>*{}**@{-},
 <-0.4mm,0.7mm>*{};<-4.5mm,5mm>*{}**@{-},
 <0mm,0mm>*{};<-1mm,5mm>*{\ldots}**@{},
 <0.4mm,0.7mm>*{};<4.5mm,5mm>*{}**@{-},
 <0.6mm,0.44mm>*{};<8mm,5mm>*{}**@{-},
   <0mm,0mm>*{};<-8.5mm,5.5mm>*{^1}**@{},
   <0mm,0mm>*{};<-5mm,5.5mm>*{^2}**@{},
   <0mm,0mm>*{};<4.5mm,5.5mm>*{^{m\hspace{-0.5mm}-\hspace{-0.5mm}1}}**@{},
   <0mm,0mm>*{};<9.0mm,5.5mm>*{^m}**@{},
 <-0.6mm,-0.44mm>*{};<-8mm,-5mm>*{}**@{-},
 <-0.4mm,-0.7mm>*{};<-4.5mm,-5mm>*{}**@{-},
 <0mm,0mm>*{};<-1mm,-5mm>*{\ldots}**@{},
 <0.4mm,-0.7mm>*{};<4.5mm,-5mm>*{}**@{-},
 <0.6mm,-0.44mm>*{};<8mm,-5mm>*{}**@{-},
   <0mm,0mm>*{};<-8.5mm,-6.9mm>*{^1}**@{},
   <0mm,0mm>*{};<-5mm,-6.9mm>*{^2}**@{},
   <0mm,0mm>*{};<4.5mm,-6.9mm>*{^{n\hspace{-0.5mm}-\hspace{-0.5mm}1}}**@{},
   <0mm,0mm>*{};<9.0mm,-6.9mm>*{^n}**@{},
 \end{xy}}
\ \ = \ \
 \sum_{[1,\ldots,m]=I_1\sqcup I_2\atop
 {|I_1|\geq 0, |I_2|\geq 1}}
 \sum_{[1,\ldots,n]=J_1\sqcup J_2\atop
 {|J_1|\geq 1, |J_2|\geq 1}
}\hspace{0mm}
\pm\ \resizebox{24mm}{!}{ \begin{xy}
 <0mm,0mm>*{\circ};<0mm,0mm>*{}**@{},
 <-0.6mm,0.44mm>*{};<-8mm,5mm>*{}**@{-},
 <-0.4mm,0.7mm>*{};<-4.5mm,5mm>*{}**@{-},
 <0mm,0mm>*{};<0mm,5mm>*{\ldots}**@{},
 <0.4mm,0.7mm>*{};<4.5mm,5mm>*{}**@{-},
 <0.6mm,0.44mm>*{};<12.4mm,4.8mm>*{}**@{-},
     <0mm,0mm>*{};<-2mm,7mm>*{\overbrace{\ \ \ \ \ \ \ \ \ \ \ \ }}**@{},
     <0mm,0mm>*{};<-2mm,9mm>*{^{I_1}}**@{},
 <-0.6mm,-0.44mm>*{};<-8mm,-5mm>*{}**@{-},
 <-0.4mm,-0.7mm>*{};<-4.5mm,-5mm>*{}**@{-},
 <0mm,0mm>*{};<-1mm,-5mm>*{\ldots}**@{},
 <0.4mm,-0.7mm>*{};<4.5mm,-5mm>*{}**@{-},
 <0.6mm,-0.44mm>*{};<8mm,-5mm>*{}**@{-},
      <0mm,0mm>*{};<0mm,-7mm>*{\underbrace{\ \ \ \ \ \ \ \ \ \ \ \ \ \ \
      }}**@{},
      <0mm,0mm>*{};<0mm,-10.6mm>*{_{J_1}}**@{},
 <13mm,5mm>*{};<13mm,5mm>*{\circ}**@{},
 <12.6mm,5.44mm>*{};<5mm,10mm>*{}**@{-},
 <12.6mm,5.7mm>*{};<8.5mm,10mm>*{}**@{-},
 <13mm,5mm>*{};<13mm,10mm>*{\ldots}**@{},
 <13.4mm,5.7mm>*{};<16.5mm,10mm>*{}**@{-},
 <13.6mm,5.44mm>*{};<20mm,10mm>*{}**@{-},
      <13mm,5mm>*{};<13mm,12mm>*{\overbrace{\ \ \ \ \ \ \ \ \ \ \ \ \ \ }}**@{},
      <13mm,5mm>*{};<13mm,14mm>*{^{I_2}}**@{},
 <12.4mm,4.3mm>*{};<8mm,0mm>*{}**@{-},
 <12.6mm,4.3mm>*{};<12mm,0mm>*{\ldots}**@{},
 <13.4mm,4.5mm>*{};<16.5mm,0mm>*{}**@{-},
 <13.6mm,4.8mm>*{};<20mm,0mm>*{}**@{-},
     <13mm,5mm>*{};<14.3mm,-2mm>*{\underbrace{\ \ \ \ \ \ \ \ \ \ \ }}**@{},
     <13mm,5mm>*{};<14.3mm,-4.5mm>*{_{J_2}}**@{},
 \end{xy}}
$$
where $\sigma(I_1\sqcup I_2)$ and $\sigma(J_1\sqcup J_2)$ are the signs of the shuffles
$[1,\ldots,m]=I_1\sqcup I_2$ and, respectively, $[1,\ldots,n]=J_1\sqcup J_2$. For example,
\[
\delta
\begin{xy}
 <0.5mm,0.5mm>*{};<2.2mm,2.2mm>*{}**@{-},
 <-0.48mm,0.48mm>*{};<-2.2mm,2.2mm>*{}**@{-},
 <0mm,0mm>*{\circ};<0mm,0mm>*{}**@{},
 <0.5mm,0.5mm>*{};<2.7mm,2.8mm>*{^2}**@{},
 <-0.48mm,0.48mm>*{};<-2.7mm,2.8mm>*{^1}**@{},
 <0.39mm,-0.39mm>*{};<2.2mm,-2.2mm>*{}**@{-},
 <-0.35mm,-0.35mm>*{};<-2.2mm,-2.2mm>*{}**@{-},
 <0mm,0mm>*{\circ};<0mm,0mm>*{}**@{},
   <0.39mm,-0.39mm>*{};<2.9mm,-4mm>*{^2}**@{},
   <-0.35mm,-0.35mm>*{};<-2.8mm,-4mm>*{^1}**@{},
\end{xy}
=
 \begin{xy}
 <0mm,2.47mm>*{};<0mm,0.12mm>*{}**@{-},
 <0.5mm,3.5mm>*{};<2.2mm,5.2mm>*{}**@{-},
 <-0.48mm,3.48mm>*{};<-2.2mm,5.2mm>*{}**@{-},
 <0mm,3mm>*{\circ};<0mm,3mm>*{}**@{},
  <0mm,-0.8mm>*{\circ};<0mm,-0.8mm>*{}**@{},
<-0.39mm,-1.2mm>*{};<-2.2mm,-3.5mm>*{}**@{-},
 <0.39mm,-1.2mm>*{};<2.2mm,-3.5mm>*{}**@{-},
     <0.5mm,3.5mm>*{};<2.8mm,5.7mm>*{^2}**@{},
     <-0.48mm,3.48mm>*{};<-2.8mm,5.7mm>*{^1}**@{},
   <0mm,-0.8mm>*{};<-2.7mm,-5.2mm>*{^1}**@{},
   <0mm,-0.8mm>*{};<2.7mm,-5.2mm>*{^2}**@{},
\end{xy}
\  - \
\begin{xy}
 <0mm,-1.3mm>*{};<0mm,-3.5mm>*{}**@{-},
 <0.38mm,-0.2mm>*{};<2.0mm,2.0mm>*{}**@{-},
 <-0.38mm,-0.2mm>*{};<-2.2mm,2.2mm>*{}**@{-},
<0mm,-0.8mm>*{\circ};<0mm,0.8mm>*{}**@{},
 <2.4mm,2.4mm>*{\circ};<2.4mm,2.4mm>*{}**@{},
 <2.77mm,2.0mm>*{};<4.4mm,-0.8mm>*{}**@{-},
 <2.4mm,3mm>*{};<2.4mm,5.2mm>*{}**@{-},
     <0mm,-1.3mm>*{};<0mm,-5.3mm>*{^1}**@{},
     <2.5mm,2.3mm>*{};<5.1mm,-2.6mm>*{^2}**@{},
    <2.4mm,2.5mm>*{};<2.4mm,5.7mm>*{^2}**@{},
    <-0.38mm,-0.2mm>*{};<-2.8mm,2.5mm>*{^1}**@{},
    \end{xy}
\  + \
\begin{xy}
 <0mm,-1.3mm>*{};<0mm,-3.5mm>*{}**@{-},
 <0.38mm,-0.2mm>*{};<2.0mm,2.0mm>*{}**@{-},
 <-0.38mm,-0.2mm>*{};<-2.2mm,2.2mm>*{}**@{-},
<0mm,-0.8mm>*{\circ};<0mm,0.8mm>*{}**@{},
 <2.4mm,2.4mm>*{\circ};<2.4mm,2.4mm>*{}**@{},
 <2.77mm,2.0mm>*{};<4.4mm,-0.8mm>*{}**@{-},
 <2.4mm,3mm>*{};<2.4mm,5.2mm>*{}**@{-},
     <0mm,-1.3mm>*{};<0mm,-5.3mm>*{^2}**@{},
     <2.5mm,2.3mm>*{};<5.1mm,-2.6mm>*{^1}**@{},
    <2.4mm,2.5mm>*{};<2.4mm,5.7mm>*{^2}**@{},
    <-0.38mm,-0.2mm>*{};<-2.8mm,2.5mm>*{^1}**@{},
    \end{xy}
\  - \
\begin{xy}
 <0mm,-1.3mm>*{};<0mm,-3.5mm>*{}**@{-},
 <0.38mm,-0.2mm>*{};<2.0mm,2.0mm>*{}**@{-},
 <-0.38mm,-0.2mm>*{};<-2.2mm,2.2mm>*{}**@{-},
<0mm,-0.8mm>*{\circ};<0mm,0.8mm>*{}**@{},
 <2.4mm,2.4mm>*{\circ};<2.4mm,2.4mm>*{}**@{},
 <2.77mm,2.0mm>*{};<4.4mm,-0.8mm>*{}**@{-},
 <2.4mm,3mm>*{};<2.4mm,5.2mm>*{}**@{-},
     <0mm,-1.3mm>*{};<0mm,-5.3mm>*{^2}**@{},
     <2.5mm,2.3mm>*{};<5.1mm,-2.6mm>*{^1}**@{},
    <2.4mm,2.5mm>*{};<2.4mm,5.7mm>*{^1}**@{},
    <-0.38mm,-0.2mm>*{};<-2.8mm,2.5mm>*{^2}**@{},
    \end{xy}
\ + \
\begin{xy}
 <0mm,-1.3mm>*{};<0mm,-3.5mm>*{}**@{-},
 <0.38mm,-0.2mm>*{};<2.0mm,2.0mm>*{}**@{-},
 <-0.38mm,-0.2mm>*{};<-2.2mm,2.2mm>*{}**@{-},
<0mm,-0.8mm>*{\circ};<0mm,0.8mm>*{}**@{},
 <2.4mm,2.4mm>*{\circ};<2.4mm,2.4mm>*{}**@{},
 <2.77mm,2.0mm>*{};<4.4mm,-0.8mm>*{}**@{-},
 <2.4mm,3mm>*{};<2.4mm,5.2mm>*{}**@{-},
     <0mm,-1.3mm>*{};<0mm,-5.3mm>*{^1}**@{},
     <2.5mm,2.3mm>*{};<5.1mm,-2.6mm>*{^2}**@{},
    <2.4mm,2.5mm>*{};<2.4mm,5.7mm>*{^1}**@{},
    <-0.38mm,-0.2mm>*{};<-2.8mm,2.5mm>*{^2}**@{},
    \end{xy}.
\]

\subsubsection{\bf Fact}
Representations,
$$
\rho: \Lieb_\infty \lon \End\left\langle V\right\rangle,
$$
of the dg prop $(\Lieb_\infty, \delta)$ in a dg space $(V,d)$ are in
one-to-one correspondence with degree $3$ elements, $\ga$, in the Poisson algebra
$(\fl_V, \{\, , \, \})$ satisfying the equation,
$\{\ga, \ga\}=0$ (see Corollary 5.1 in \cite{Me1}).

Indeed, using natural degree $m+n$ isomorphisms,
$$
s_m^n: \Hom\left(\wedge^n V, \wedge^m V\right) \lon \odot^n(V[1])\ot \odot^m(V^*[1]),
$$
we define a degree 3 element,
\[
\ga:= s_1^1(d) + \sum_{m,n\geq 1\atop m+n \geq 3}
s_n^m\circ \rho\left(
\resizebox{14mm}{!}{\begin{xy}
 <0mm,0mm>*{\circ};<0mm,0mm>*{}**@{},
 <-0.6mm,0.44mm>*{};<-8mm,5mm>*{}**@{-},
 <-0.4mm,0.7mm>*{};<-4.5mm,5mm>*{}**@{-},
 <0mm,0mm>*{};<-1mm,5mm>*{\ldots}**@{},
 <0.4mm,0.7mm>*{};<4.5mm,5mm>*{}**@{-},
 <0.6mm,0.44mm>*{};<8mm,5mm>*{}**@{-},
   <0mm,0mm>*{};<-8.5mm,5.5mm>*{^1}**@{},
   <0mm,0mm>*{};<-5mm,5.5mm>*{^2}**@{},
   <0mm,0mm>*{};<4.5mm,5.5mm>*{^{m\hspace{-0.5mm}-\hspace{-0.5mm}1}}**@{},
   <0mm,0mm>*{};<9.0mm,5.5mm>*{^m}**@{},
 <-0.6mm,-0.44mm>*{};<-8mm,-5mm>*{}**@{-},
 <-0.4mm,-0.7mm>*{};<-4.5mm,-5mm>*{}**@{-},
 <0mm,0mm>*{};<-1mm,-5mm>*{\ldots}**@{},
 <0.4mm,-0.7mm>*{};<4.5mm,-5mm>*{}**@{-},
 <0.6mm,-0.44mm>*{};<8mm,-5mm>*{}**@{-},
   <0mm,0mm>*{};<-8.5mm,-6.9mm>*{^1}**@{},
   <0mm,0mm>*{};<-5mm,-6.9mm>*{^2}**@{},
   <0mm,0mm>*{};<4.5mm,-6.9mm>*{^{n\hspace{-0.5mm}-\hspace{-0.5mm}1}}**@{},
   <0mm,0mm>*{};<9.0mm,-6.9mm>*{^n}**@{},
 \end{xy}}
\right)\in \f_\cV,
\]
and then check that the  equation $d\circ \rho=\rho\circ \delta$ translates precisely into the Maurer-Cartan
equation $\{\ga,\ga\}=0$.
The Lie bialgebra structures on  $V$  get identified with homogeneous polynomials of order 3,
$$
\ga_3\in \odot^2(V[-1])\ot V^*[-1]\ \  \oplus\ \  V[-1]\ot \odot^2(V^*[-1]),
 $$
 satisfying the equation $\{\ga_3, \ga_3\}=0$. Indeed, every such a polynomial is equivalent
to a pair $(\vartriangle\in \Hom(V, \wedge^2V),\ [\, ,\, ]\in \Hom(\wedge^2V,V))$, in terms of which the single
 equation  $\{\ga_3,\ga_3\}=0$ decomposes into relations (i)-(iii) of \S \ref{2.1}.


\section{Prop of strongly homotopy bialgebras}

\subsection{Associative bialgebras.} A {\em  bialgebra}\,
 is, by definition, a graded vector space $V$ equipped with two degree zero linear
maps,
$$
\Ba{rccc}
\mu: &  V\ot V& \lon & V \\
       & a\ot b    & \lon & ab
\Ea
\ \ \ \ \  , \ \ \ \
\Ba{rccc}
\Delta: & V& \lon & V\ot V \\
       & a    & \lon &  \sum a_{1}\ot a_{2}
\Ea
$$
satisfying,
\Bi
\item[(i)] the associativity identity:
$(ab)c=a(bc)$;
\item[(ii)] the coassociativity identity: $(\Delta\ot\Id)\Delta a=
(\Id\ot \Delta)\Delta a$;
\item[(iii)] the compatibility identity: $\Delta$ is a morphism of algebras,
i.e.\ $\Delta(ab)=\sum (-1)^{a_2b_1} a_1b_1\ot a_2b_2$,
\Ei
for any $a,b, c\in V$. We often abbreviate ``associative bialgebra" to simply ``bialgebra".

\subsection{Prop of bialgebras} There exists  a prop, $\sB$, whose
representations,
$$
\rho: \sB \lon \End\langle V\rangle,
$$
in a graded vector space $V$ are in one-to-one correspondence with the bialgebra structures
in $V$ \cite{EE}. With an association in mind,
$$
\Delta \leftrightarrow
 \begin{xy}
 <0mm,-0.55mm>*{};<0mm,-2.5mm>*{}**@{-},
 <0.5mm,0.5mm>*{};<2.2mm,2.2mm>*{}**@{-},
 <-0.48mm,0.48mm>*{};<-2.2mm,2.2mm>*{}**@{-},
 <0mm,0mm>*{\bullet};<0mm,0mm>*{}**@{},
 \end{xy}
 \ \ \  , \ \ \
\mu \leftrightarrow
 \begin{xy}
 <0mm,0.66mm>*{};<0mm,3mm>*{}**@{-},
 <0.39mm,-0.39mm>*{};<2.2mm,-2.2mm>*{}**@{-},
 <-0.35mm,-0.35mm>*{};<-2.2mm,-2.2mm>*{}**@{-},
 <0mm,0mm>*{\bullet};<0mm,0mm>*{}**@{},
 \end{xy}
$$
one can define it as a quotient,
$$
\sB:= {\Gamma\langle E \rangle}/(R)
$$
of the free prop, $\Gamma\langle E \rangle$, generated by the $\bS$-bimodule $E=\{E(m,n)\}$,
\[
E(m,n):=\left\{
\Ba{rr}
\K[\bS_2]\ot \id_1\equiv\mbox{span}\left\langle
\begin{xy}
 <0mm,-0.55mm>*{};<0mm,-2.5mm>*{}**@{-},
 <0.5mm,0.5mm>*{};<2.2mm,2.2mm>*{}**@{-},
 <-0.48mm,0.48mm>*{};<-2.2mm,2.2mm>*{}**@{-},
 <0mm,0mm>*{\bullet};<0mm,0mm>*{}**@{},
 <0mm,-0.55mm>*{};<0mm,-3.8mm>*{_1}**@{},
 <0.5mm,0.5mm>*{};<2.7mm,2.8mm>*{^2}**@{},
 <-0.48mm,0.48mm>*{};<-2.7mm,2.8mm>*{^1}**@{},
 \end{xy}
\,
,\,
\begin{xy}
 <0mm,-0.55mm>*{};<0mm,-2.5mm>*{}**@{-},
 <0.5mm,0.5mm>*{};<2.2mm,2.2mm>*{}**@{-},
 <-0.48mm,0.48mm>*{};<-2.2mm,2.2mm>*{}**@{-},
 <0mm,0mm>*{\bullet};<0mm,0mm>*{}**@{},
 <0mm,-0.55mm>*{};<0mm,-3.8mm>*{_1}**@{},
 <0.5mm,0.5mm>*{};<2.7mm,2.8mm>*{^1}**@{},
 <-0.48mm,0.48mm>*{};<-2.7mm,2.8mm>*{^2}**@{},
 \end{xy}
   \right\rangle  & \mbox{if}\ m=2, n=1,\vspace{3mm}\\
\id_1\ot \K[\bS_2]\equiv
\mbox{span}\left\langle
\begin{xy}
 <0mm,0.66mm>*{};<0mm,3mm>*{}**@{-},
 <0.39mm,-0.39mm>*{};<2.2mm,-2.2mm>*{}**@{-},
 <-0.35mm,-0.35mm>*{};<-2.2mm,-2.2mm>*{}**@{-},
 <0mm,0mm>*{\bullet};<0mm,0mm>*{}**@{},
   <0mm,0.66mm>*{};<0mm,3.4mm>*{^1}**@{},
   <0.39mm,-0.39mm>*{};<2.9mm,-4mm>*{^2}**@{},
   <-0.35mm,-0.35mm>*{};<-2.8mm,-4mm>*{^1}**@{},
\end{xy}
\,
,\,
\begin{xy}
 <0mm,0.66mm>*{};<0mm,3mm>*{}**@{-},
 <0.39mm,-0.39mm>*{};<2.2mm,-2.2mm>*{}**@{-},
 <-0.35mm,-0.35mm>*{};<-2.2mm,-2.2mm>*{}**@{-},
 <0mm,0mm>*{\bullet};<0mm,0mm>*{}**@{},
   <0mm,0.66mm>*{};<0mm,3.4mm>*{^1}**@{},
   <0.39mm,-0.39mm>*{};<2.9mm,-4mm>*{^1}**@{},
   <-0.35mm,-0.35mm>*{};<-2.8mm,-4mm>*{^2}**@{},
\end{xy}
\right\rangle
\ & \mbox{if}\ m=1, n=2, \vspace{3mm}\\
0 & \mbox{otherwise}
\Ea
\right.
\]
modulo the ideal generated by relations
\[
R:\left\{
\Ba{r}
\resizebox{8mm}{!}{\begin{xy}
 <0mm,0mm>*{\bullet};<0mm,0mm>*{}**@{},
 <0mm,-0.49mm>*{};<0mm,-3.0mm>*{}**@{-},
 <0.49mm,0.49mm>*{};<1.9mm,1.9mm>*{}**@{-},
 <-0.5mm,0.5mm>*{};<-1.9mm,1.9mm>*{}**@{-},
 <-2.3mm,2.3mm>*{\bullet};<-2.3mm,2.3mm>*{}**@{},
 <-1.8mm,2.8mm>*{};<0mm,4.9mm>*{}**@{-},
 <-2.8mm,2.9mm>*{};<-4.6mm,4.9mm>*{}**@{-},
   <0.49mm,0.49mm>*{};<2.7mm,2.3mm>*{^3}**@{},
   <-1.8mm,2.8mm>*{};<0.4mm,5.3mm>*{^2}**@{},
   <-2.8mm,2.9mm>*{};<-5.1mm,5.3mm>*{^1}**@{},
 \end{xy}}
\ - \
\resizebox{8mm}{!}{\begin{xy}
 <0mm,0mm>*{\bullet};<0mm,0mm>*{}**@{},
 <0mm,-0.49mm>*{};<0mm,-3.0mm>*{}**@{-},
 <0.49mm,0.49mm>*{};<1.9mm,1.9mm>*{}**@{-},
 <-0.5mm,0.5mm>*{};<-1.9mm,1.9mm>*{}**@{-},
 <2.3mm,2.3mm>*{\bullet};<-2.3mm,2.3mm>*{}**@{},
 <1.8mm,2.8mm>*{};<0mm,4.9mm>*{}**@{-},
 <2.8mm,2.9mm>*{};<4.6mm,4.9mm>*{}**@{-},
   <0.49mm,0.49mm>*{};<-2.7mm,2.3mm>*{^1}**@{},
   <-1.8mm,2.8mm>*{};<0mm,5.3mm>*{^2}**@{},
   <-2.8mm,2.9mm>*{};<5.1mm,5.3mm>*{^3}**@{},
 \end{xy}}
\ \ \ \in \Gamma \langle E\rangle(3,1) \vspace{3mm}\\
 \resizebox{8mm}{!}{\begin{xy}
 <0mm,0mm>*{\bullet};<0mm,0mm>*{}**@{},
 <0mm,0.69mm>*{};<0mm,3.0mm>*{}**@{-},
 <0.39mm,-0.39mm>*{};<2.4mm,-2.4mm>*{}**@{-},
 <-0.35mm,-0.35mm>*{};<-1.9mm,-1.9mm>*{}**@{-},
 <-2.4mm,-2.4mm>*{\bullet};<-2.4mm,-2.4mm>*{}**@{},
 <-2.0mm,-2.8mm>*{};<0mm,-4.9mm>*{}**@{-},
 <-2.8mm,-2.9mm>*{};<-4.7mm,-4.9mm>*{}**@{-},
    <0.39mm,-0.39mm>*{};<3.3mm,-4.0mm>*{^3}**@{},
    <-2.0mm,-2.8mm>*{};<0.5mm,-6.7mm>*{^2}**@{},
    <-2.8mm,-2.9mm>*{};<-5.2mm,-6.7mm>*{^1}**@{},
 \end{xy}}
\ - \
\resizebox{8mm}{!}{ \begin{xy}
 <0mm,0mm>*{\bullet};<0mm,0mm>*{}**@{},
 <0mm,0.69mm>*{};<0mm,3.0mm>*{}**@{-},
 <0.39mm,-0.39mm>*{};<2.4mm,-2.4mm>*{}**@{-},
 <-0.35mm,-0.35mm>*{};<-1.9mm,-1.9mm>*{}**@{-},
 <2.4mm,-2.4mm>*{\bullet};<-2.4mm,-2.4mm>*{}**@{},
 <2.0mm,-2.8mm>*{};<0mm,-4.9mm>*{}**@{-},
 <2.8mm,-2.9mm>*{};<4.7mm,-4.9mm>*{}**@{-},
    <0.39mm,-0.39mm>*{};<-3mm,-4.0mm>*{^1}**@{},
    <-2.0mm,-2.8mm>*{};<0mm,-6.7mm>*{^2}**@{},
    <-2.8mm,-2.9mm>*{};<5.2mm,-6.7mm>*{^3}**@{},
 \end{xy}}
\  \ \ \in  \Gamma \langle E\rangle(1,3) \vspace{2mm} \\
\resizebox{8mm}{!}{ \begin{xy}
 <0mm,2.47mm>*{};<0mm,-0.5mm>*{}**@{-},
 <0.5mm,3.5mm>*{};<2.2mm,5.2mm>*{}**@{-},
 <-0.48mm,3.48mm>*{};<-2.2mm,5.2mm>*{}**@{-},
 <0mm,3mm>*{\bullet};<0mm,3mm>*{}**@{},
  <0mm,-0.8mm>*{\bullet};<0mm,-0.8mm>*{}**@{},
<0mm,-0.8mm>*{};<-2.2mm,-3.5mm>*{}**@{-},
 <0mm,-0.8mm>*{};<2.2mm,-3.5mm>*{}**@{-},
     <0.5mm,3.5mm>*{};<2.8mm,5.7mm>*{^2}**@{},
     <-0.48mm,3.48mm>*{};<-2.8mm,5.7mm>*{^1}**@{},
   <0mm,-0.8mm>*{};<-2.7mm,-5.2mm>*{^1}**@{},
   <0mm,-0.8mm>*{};<2.7mm,-5.2mm>*{^2}**@{},
\end{xy}}
\ - \
\resizebox{11mm}{!}{\begin{xy}
 <0mm,0mm>*{\bullet};<0mm,0mm>*{}**@{},
 <0mm,-0.49mm>*{};<0mm,-3.0mm>*{}**@{-},
 <-0.5mm,0.5mm>*{};<-3mm,2mm>*{}**@{-},
 <-3mm,2mm>*{};<0mm,4mm>*{}**@{-},
 <0mm,4mm>*{\bullet};<-2.3mm,2.3mm>*{}**@{},
 <0mm,4mm>*{};<0mm,7.4mm>*{}**@{-},
<0mm,0mm>*{};<2.2mm,1.5mm>*{}**@{-},
 <6mm,0mm>*{\bullet};<0mm,0mm>*{}**@{},
 <6mm,4mm>*{};<3.8mm,2.5mm>*{}**@{-},
 <6mm,4mm>*{};<6mm,7.4mm>*{}**@{-},
 <6mm,4mm>*{\bullet};<-2.3mm,2.3mm>*{}**@{},
 <0mm,4mm>*{};<6mm,0mm>*{}**@{-},
<6mm,4mm>*{};<9mm,2mm>*{}**@{-},
<6mm,0mm>*{};<9mm,2mm>*{}**@{-},
<6mm,0mm>*{};<6mm,-3mm>*{}**@{-},
   <-1.8mm,2.8mm>*{};<0mm,7.8mm>*{^1}**@{},
   <-2.8mm,2.9mm>*{};<0mm,-4.3mm>*{_1}**@{},
<-1.8mm,2.8mm>*{};<6mm,7.8mm>*{^2}**@{},
   <-2.8mm,2.9mm>*{};<6mm,-4.3mm>*{_2}**@{},
 \end{xy}}
\ \ \ \in  \Gamma \langle E\rangle(2,2)
\Ea
\right.
\]
which are not quadratic in the properadic sense \cite{Va}.

\subsection{A minimal resolution of $\sB$}\label{min}
Existence of a minimal resolution of the prop of bialgebras,
$$
(\sB_\infty=\Gamma \langle\mathsf E\rangle, \delta) \stackrel{\pi}{\lon} (\sB, 0)
$$
was proven by M.\ Markl \cite{Ma} who showed that $\sB_\infty$
is freely generated by a relatively small $\bS$-bimodule
 $\mathsf E=\{\mathsf E(m,n)\}_{m,n\geq 1, m+n\geq 3}$, with
\[
\mathsf E(m,n):= \K[\bS_m]\ot \K[\bS_n][m+n-3]=\mbox{span}\left\langle
\resizebox{15mm}{!}{\begin{xy}
 <0mm,0mm>*{\bullet};<0mm,0mm>*{}**@{},
 <0mm,0mm>*{};<-8mm,5mm>*{}**@{-},
 <0mm,0mm>*{};<-4.5mm,5mm>*{}**@{-},
 <0mm,0mm>*{};<-1mm,5mm>*{\ldots}**@{},
 <0mm,0mm>*{};<4.5mm,5mm>*{}**@{-},
 <0mm,0mm>*{};<8mm,5mm>*{}**@{-},
   <0mm,0mm>*{};<-8.5mm,5.5mm>*{^1}**@{},
   <0mm,0mm>*{};<-5mm,5.5mm>*{^2}**@{},
   <0mm,0mm>*{};<4.5mm,5.5mm>*{^{m\hspace{-0.5mm}-\hspace{-0.5mm}1}}**@{},
   <0mm,0mm>*{};<9.0mm,5.5mm>*{^m}**@{},
 <0mm,0mm>*{};<-8mm,-5mm>*{}**@{-},
 <0mm,0mm>*{};<-4.5mm,-5mm>*{}**@{-},
 <0mm,0mm>*{};<-1mm,-5mm>*{\ldots}**@{},
 <0mm,0mm>*{};<4.5mm,-5mm>*{}**@{-},
 <0mm,0mm>*{};<8mm,-5mm>*{}**@{-},
   <0mm,0mm>*{};<-8.5mm,-6.9mm>*{^1}**@{},
   <0mm,0mm>*{};<-5mm,-6.9mm>*{^2}**@{},
   <0mm,0mm>*{};<4.5mm,-6.9mm>*{^{n\hspace{-0.5mm}-\hspace{-0.5mm}1}}**@{},
   <0mm,0mm>*{};<9.0mm,-6.9mm>*{^n}**@{},
 \end{xy}}
\right\rangle
\]
The differential $\delta$ is not
quadratic, and its explicit value on the generic $(m,n)$-corolla is not
known at present.
Rather surprisingly,  just existence of $(\sB_\infty, \delta)$
 and a theorem on Gerstenhaber-Schack differential from \cite{MV}
(recalled in \S~\ref{fact-GS} below)
are enough for our purposes.

\subsection{Deformation theory of dg morphisms.}\label{deftheory}
Let $(\Gamma \langle \mathsf A \rangle, \delta)$ be a dg  free prop generated by an $\bS$-bimodule $\mathsf A $, and
$(\mathsf P,d)$ an arbitrary dg prop. It was shown in \cite{MV} that the deformation complex of the zero morphism,
$$
\fg:=\Def(\Gamma \langle \mathsf A \rangle \stackrel{0}{\rar} \sP)\simeq  \Hom_{\bS}(\mathsf A, \mathsf P)[-1],
$$
has an induced (from the differentials $\delta$ and $d$) filtered $L_\infty$-structure,
\[
\left\{ \mu_n: \odot^n(\fg[1]) \lon \fg[1]   \right\}_{n\geq 1},
\]
whose Maurer-Cartan elements, that is, degree 1 elements $\ga$ in $\fg$ satisfying a well-defined equation
$$
\mu_1(\ga) + \frac{1}{2!}\mu_2(\ga, \ga) + \frac{1}{3!}\mu_2(\ga, \ga, \ga) + \ldots =0,
$$
are in one-to-one correspondence with morphisms, $\ga:(\Gamma \langle \mathsf A \rangle, \delta)\rar (\mathsf P, d)$,
of dg props.

\sip

If $\ga: (\Gamma \langle \mathsf A \rangle, \delta)\rar (\mathsf P,d)$ is  any particular morphism of dg props,
then $\fg$ has a
canonical $\ga$-{\em twisted}\, $L_\infty$-structure,
\[
\left\{ \mu_n^\ga: \odot^n(\fg[1]) \lon \fg[1]   \right\}_{n\geq 1},
\]
which controls deformation theory of the morphism $\gamma$. The associated Maurer-Cartan elements, $\Ga$,
$$
\mu_1^\ga(\Ga) + \frac{1}{2!}\mu_2^\ga(\Ga, \Ga) + \frac{1}{3!}\mu_2^\ga(\Ga, \Ga, \Ga) + \ldots =0,
$$
are in one-to-one correspondence with those morphisms of dg props,
$(\Gamma \langle \mathsf A \rangle, \delta)\rar (\mathsf P,d)$
 whose values on
generators are given by  $\ga+\Ga$.

\subsubsection{\mbox{\bf Deformation theory of bialgebras}.}\label{fact-GS} If $\ga: \sB \rar \End_ W$ is a bialgebra structure,
 then the differential $\mu_1^\ga$ of the induced $L_\infty$ structure in
$$
\fg\fs(W,W):=\Def(\sB\stackrel{\ga}{\rar}\mathsf{End}_V)\simeq \Hom_{\bS}(\mathsf E, \End(V))[-1]=\bigoplus_{m,n\geq 1\atop m+n\geq 3} \Hom(V^{\ot n}, V^{\ot m})[m+n-2],
$$
is precisely \cite{MV}
the Gerstenhaber-Schack differential \cite{GS},
$$
d_{\fg\fs}=d_1 \oplus d_2:  \Hom(V^{\ot n}, V^{\ot m}) \lon  \Hom(V^{\ot n+1}, V^{\ot m}) \oplus
\Hom(V^{\ot n}, V^{\ot m+1}),
$$
with $d_1$ given on an arbitrary $f\in \Hom(V^{\ot n}, V^{\ot m})$  by
\Beqrn
(d_1f)(v_0, v_1, \ldots, v_n)&:=&\Delta^{m-1}(v_0)\cdot f(v_1, v_2, \ldots, v_n) -
\sum_{i=0}^{n-1} (-1)^if(v_1, \dots, v_iv_{i+1}, \ldots, v_n)\\
&& + (-1)^{n+1}f(v_1, v_2, \ldots, v_n)\cdot\Delta^{m-1}(v_n)
\forall\ \ v_0, v_1,\ldots, v_n\in V,
\Eeqrn
where multiplication in $V$ is denoted by juxtaposition, the induced multiplication in the algebra $V^{\ot m}$ by
$\cdot$, the comultiplication in $V$ by $\Delta$, and
$$
\Delta^{m-1}: (\Delta\ot \Id^{\ot m-2 })\circ (\Delta\ot \Id^{\ot m-3})\circ \ldots \circ \Delta: V \rar V^{\ot m},
$$
for $m\geq 2$ while $\Delta^0:=\Id$.
The expression for $d_2$ is an obvious ``dual" analogue of the one for $d_1$.

\bip

\section{The Gerstenhaber-Schack complex  of polydifferential operators\\ and a dg
prop $\DefQ^+$
}

\bip

\subsection{An endofunctor in the category of prop(erad)s} Consider an endofunctor,
$^+: \sP \rar \sP^+$, on the (sub)category of dg {\em free}\, props, $\sP=(\Ga\langle E\rangle, \delta)$, $E=\{E(m,n)\}$ being an
$\bS$-bimodule and $\delta$ a differential. Define a new $\bS$-module,
$$
E^+(m,n):=\left\{ \Ba{rc} E(1)\oplus \K[-1] & \mbox{if}\ m=n=1,\\
   E(m,n) & \mbox{otherwise}.
\Ea\right.
$$
If we denote pictorially a generator of the summand $\K[-1]\subset E^+(1,1)$ by $\begin{xy}
 <0mm,-0.55mm>*{};<0mm,-3mm>*{}**@{-},
 <0mm,0.5mm>*{};<0mm,3mm>*{}**@{-},
 <0mm,0mm>*{\bullet};<0mm,0mm>*{}**@{},
 \end{xy}$,
and a generator (of, say, homological degree $a$) of $E(m,n)$ by an $(m,n)$-corolla
$$
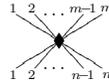
\resizebox{14mm}{!}{\begin{xy}
 <0mm,0mm>*{\blacklozenge};<0mm,0mm>*{}**@{},
 <0mm,0mm>*{};<-8mm,5mm>*{}**@{-},
 <0mm,0mm>*{};<-4.5mm,5mm>*{}**@{-},
 <0mm,0mm>*{};<-1mm,5mm>*{\ldots}**@{},
 <0mm,0mm>*{};<4.5mm,5mm>*{}**@{-},
 <0mm,0mm>*{};<8mm,5mm>*{}**@{-},
   <0mm,0mm>*{};<-8.5mm,5.5mm>*{^1}**@{},
   <0mm,0mm>*{};<-5mm,5.5mm>*{^2}**@{},
   <0mm,0mm>*{};<4.5mm,5.5mm>*{^{m\hspace{-0.5mm}-\hspace{-0.5mm}1}}**@{},
   <0mm,0mm>*{};<9.0mm,5.5mm>*{^m}**@{},
 <0mm,0mm>*{};<-8mm,-5mm>*{}**@{-},
 <0mm,0mm>*{};<-4.5mm,-5mm>*{}**@{-},
 <0mm,0mm>*{};<-1mm,-5mm>*{\ldots}**@{},
 <0mm,0mm>*{};<4.5mm,-5mm>*{}**@{-},
 <0mm,0mm>*{};<8mm,-5mm>*{}**@{-},
   <0mm,0mm>*{};<-8.5mm,-6.9mm>*{^1}**@{},
   <0mm,0mm>*{};<-5mm,-6.9mm>*{^2}**@{},
   <0mm,0mm>*{};<4.5mm,-6.9mm>*{^{n\hspace{-0.5mm}-\hspace{-0.5mm}1}}**@{},
   <0mm,0mm>*{};<9.0mm,-6.9mm>*{^n}**@{},
 \end{xy}}
$$
then $P^+$ is defined as the free prop $\Ga\langle E^+\rangle$ equipped with the following differential,
\Beqrn
\delta^+ \begin{xy}
 <0mm,-0.55mm>*{};<0mm,-3mm>*{}**@{-},
 <0mm,0.5mm>*{};<0mm,3mm>*{}**@{-},
 <0mm,0mm>*{\bullet};<0mm,0mm>*{}**@{},
 \end{xy}  &=& \Ba{c}\begin{xy}
 <0mm,0mm>*{};<0mm,-3mm>*{}**@{-},
 <0mm,0mm>*{};<0mm,6mm>*{}**@{-},
 <0mm,0mm>*{\bullet};
 <0mm,3mm>*{\bullet};
 \end{xy}\Ea\\
\delta^+
\resizebox{15mm}{!}{ \begin{xy}
 <0mm,0mm>*{\blacklozenge};<0mm,0mm>*{}**@{},
 <0mm,0mm>*{};<-8mm,5mm>*{}**@{-},
 <0mm,0mm>*{};<-4.5mm,5mm>*{}**@{-},
 <0mm,0mm>*{};<-1mm,5mm>*{\ldots}**@{},
 <0mm,0mm>*{};<4.5mm,5mm>*{}**@{-},
 <0mm,0mm>*{};<8mm,5mm>*{}**@{-},
   <0mm,0mm>*{};<-8.5mm,5.5mm>*{^1}**@{},
   <0mm,0mm>*{};<-5mm,5.5mm>*{^2}**@{},
   <0mm,0mm>*{};<4.5mm,5.5mm>*{^{m\hspace{-0.5mm}-\hspace{-0.5mm}1}}**@{},
   <0mm,0mm>*{};<9.0mm,5.5mm>*{^m}**@{},
 <0mm,0mm>*{};<-8mm,-5mm>*{}**@{-},
 <0mm,0mm>*{};<-4.5mm,-5mm>*{}**@{-},
 <0mm,0mm>*{};<-1mm,-5mm>*{\ldots}**@{},
 <0mm,0mm>*{};<4.5mm,-5mm>*{}**@{-},
 <0mm,0mm>*{};<8mm,-5mm>*{}**@{-},
   <0mm,0mm>*{};<-8.5mm,-6.9mm>*{^1}**@{},
   <0mm,0mm>*{};<-5mm,-6.9mm>*{^2}**@{},
   <0mm,0mm>*{};<4.5mm,-6.9mm>*{^{n\hspace{-0.5mm}-\hspace{-0.5mm}1}}**@{},
   <0mm,0mm>*{};<9.0mm,-6.9mm>*{^n}**@{},
 \end{xy}}
 &=&
\delta\resizebox{15mm}{!}{ \begin{xy}
 <0mm,0mm>*{\blacklozenge};<0mm,0mm>*{}**@{},
 <0mm,0mm>*{};<-8mm,5mm>*{}**@{-},
 <0mm,0mm>*{};<-4.5mm,5mm>*{}**@{-},
 <0mm,0mm>*{};<-1mm,5mm>*{\ldots}**@{},
 <0mm,0mm>*{};<4.5mm,5mm>*{}**@{-},
 <0mm,0mm>*{};<8mm,5mm>*{}**@{-},
   <0mm,0mm>*{};<-8.5mm,5.5mm>*{^1}**@{},
   <0mm,0mm>*{};<-5mm,5.5mm>*{^2}**@{},
   <0mm,0mm>*{};<4.5mm,5.5mm>*{^{m\hspace{-0.5mm}-\hspace{-0.5mm}1}}**@{},
   <0mm,0mm>*{};<9.0mm,5.5mm>*{^m}**@{},
 <0mm,0mm>*{};<-8mm,-5mm>*{}**@{-},
 <0mm,0mm>*{};<-4.5mm,-5mm>*{}**@{-},
 <0mm,0mm>*{};<-1mm,-5mm>*{\ldots}**@{},
 <0mm,0mm>*{};<4.5mm,-5mm>*{}**@{-},
 <0mm,0mm>*{};<8mm,-5mm>*{}**@{-},
   <0mm,0mm>*{};<-8.5mm,-6.9mm>*{^1}**@{},
   <0mm,0mm>*{};<-5mm,-6.9mm>*{^2}**@{},
   <0mm,0mm>*{};<4.5mm,-6.9mm>*{^{n\hspace{-0.5mm}-\hspace{-0.5mm}1}}**@{},
   <0mm,0mm>*{};<9.0mm,-6.9mm>*{^n}**@{},
 \end{xy}}
-
\overset{m-1}{\underset{i=0}{\sum}}(-1)^{a}
\resizebox{15mm}{!}{ \begin{xy}
 <0mm,0mm>*{\blacklozenge};<0mm,0mm>*{}**@{},
 <0mm,0mm>*{};<-8mm,5mm>*{}**@{-},
 <0mm,0mm>*{};<-3.5mm,5mm>*{}**@{-},
 <0mm,0mm>*{};<-6mm,5mm>*{..}**@{},
 <0mm,0mm>*{};<0mm,5mm>*{}**@{-},
  <0mm,5mm>*{\bullet};
  <0mm,5mm>*{};<0mm,8mm>*{}**@{-},
  <0mm,5mm>*{};<0mm,9mm>*{^{i\hspace{-0.2mm}+\hspace{-0.5mm}1}}**@{},
<0mm,0mm>*{};<8mm,5mm>*{}**@{-},
<0mm,0mm>*{};<3.5mm,5mm>*{}**@{-},
 <0mm,0mm>*{};<6mm,5mm>*{..}**@{},
   <0mm,0mm>*{};<-8.5mm,5.5mm>*{^1}**@{},
   <0mm,0mm>*{};<-4mm,5.5mm>*{^i}**@{},
   <0mm,0mm>*{};<9.0mm,5.5mm>*{^m}**@{},
 <0mm,0mm>*{};<-8mm,-5mm>*{}**@{-},
 <0mm,0mm>*{};<-4.5mm,-5mm>*{}**@{-},
 <0mm,0mm>*{};<-1mm,-5mm>*{\ldots}**@{},
 <0mm,0mm>*{};<4.5mm,-5mm>*{}**@{-},
 <0mm,0mm>*{};<8mm,-5mm>*{}**@{-},
   <0mm,0mm>*{};<-8.5mm,-6.9mm>*{^1}**@{},
   <0mm,0mm>*{};<-5mm,-6.9mm>*{^2}**@{},
   <0mm,0mm>*{};<4.5mm,-6.9mm>*{^{n\hspace{-0.5mm}-\hspace{-0.5mm}1}}**@{},
   <0mm,0mm>*{};<9.0mm,-6.9mm>*{^n}**@{},
 \end{xy}}
 + \
\overset{n-1}{\underset{i=0}{\sum}}\ \
\resizebox{15mm}{!}{  \begin{xy}
 <0mm,0mm>*{\blacklozenge};<0mm,0mm>*{}**@{},
 <0mm,0mm>*{};<-8mm,-5mm>*{}**@{-},
 <0mm,0mm>*{};<-3.5mm,-5mm>*{}**@{-},
 <0mm,0mm>*{};<-6mm,-5mm>*{..}**@{},
 <0mm,0mm>*{};<0mm,-5mm>*{}**@{-},
  <0mm,-5mm>*{\bullet};
  <0mm,-5mm>*{};<0mm,-8mm>*{}**@{-},
  <0mm,-5mm>*{};<0mm,-10mm>*{^{i\hspace{-0.2mm}+\hspace{-0.5mm}1}}**@{},
<0mm,0mm>*{};<8mm,-5mm>*{}**@{-},
<0mm,0mm>*{};<3.5mm,-5mm>*{}**@{-},
 <0mm,0mm>*{};<6mm,-5mm>*{..}**@{},
   <0mm,0mm>*{};<-8.5mm,-6.9mm>*{^1}**@{},
   <0mm,0mm>*{};<-4mm,-6.9mm>*{^i}**@{},
   <0mm,0mm>*{};<9.0mm,-6.9mm>*{^n}**@{},
 <0mm,0mm>*{};<-8mm,5mm>*{}**@{-},
 <0mm,0mm>*{};<-4.5mm,5mm>*{}**@{-},
 <0mm,0mm>*{};<-1mm,5mm>*{\ldots}**@{},
 <0mm,0mm>*{};<4.5mm,5mm>*{}**@{-},
 <0mm,0mm>*{};<8mm,5mm>*{}**@{-},
   <0mm,0mm>*{};<-8.5mm,5.5mm>*{^1}**@{},
   <0mm,0mm>*{};<-5mm,5.5mm>*{^2}**@{},
   <0mm,0mm>*{};<4.5mm,5.5mm>*{^{m\hspace{-0.5mm}-\hspace{-0.5mm}1}}**@{},
   <0mm,0mm>*{};<9.0mm,5.5mm>*{^m}**@{},
 \end{xy}}
\Eeqrn

Next, let $\sP=(\Ga\langle E\rangle/I, \delta)$ be a dg prop with generators $E$ and relations described by an ideal\footnote{We tacitly assume in this paper that the ideal $I$ always satisfies the condition that $I$ belongs to the subspace of the free prop $\Ga\langle E\rangle$ spanned by graphs with at most $k$ vertices for some $k\in \N$; in this case it makes sense to talk about the vertex completion of $\sP$ which is denoted by $\widehat{\sP}$. This assumption covers all the particular props we consider in this paper}  $I$. Then we set $\sP^+=(\Ga\langle E^+\rangle/I,\delta^+)$,
the point is that the ideal $I$ is respected automatically by the differential $\delta^+$ if it is respected by $\delta$. There is obviously a 1-1 correspondence between representations,
$$
\rho^+: \sP^+ \lon \mathsf{E nd}_V,
$$
in a dg vector space $(V,d)$,
and representations
$$
\rho: \sP \lon \mathsf{End}_V,
$$
in a dg vector space $V$ equipped with the deformed differential $d - \rho^+(\begin{xy}
 <0mm,-0.55mm>*{};<0mm,-3mm>*{}**@{-},
 <0mm,0.5mm>*{};<0mm,3mm>*{}**@{-},
 <0mm,0mm>*{\bullet};<0mm,0mm>*{}**@{},
 \end{xy} )$.  Note that there is   a natural forgetful map,
 $
 \sP^+ \rar \sP,
 $
of dg props so that every representation of $\sP$ is automatically a  representation of
$\sP^+$.

\sip

As an example, we have an identification of Lie algebras,
$
\fl_V= \Def(\caL ie\cB_\infty^+ \stackrel{0}{\rar} \cE nd_V)$,
where $\fl_V$ was introduced in \S 1.
 The moral is that the functor $^+$ takes care about deformations of the differential in a representation space $V$.

\sip

However the functor $^+$ kills cohomology (at least in the case of vertex completed props). If  $\sP=(\Ga\langle E\rangle, \delta)$ is a dg (free) prop
and $\widehat{\sP}$ is its vertex completion, then the associated dg prop $\widehat{\sP}^+$
is always acyclic.
Indeed, call the generator(=decorated corolla)
$\begin{xy}
 <0mm,-0.55mm>*{};<0mm,-3mm>*{}**@{-},
 <0mm,0.5mm>*{};<0mm,3mm>*{}**@{-},
 <0mm,0mm>*{\bullet};<0mm,0mm>*{}**@{},
 \end{xy}$ of $\widehat{\sP}^+$ {\em special}, all other generators {\em non-special},
and consider a filtration of $\widehat{\sP}^+$
by the number of non-special corollas. The differential in the associated graded
is $\delta_{+}$ which acts non-trivially  only the special corollas. Call a special $(1,1)$-vertex of a graph  $\Ga$ from $\sP^+$ {\em very special}\, if it belongs to the path connecting the labelled by $1$ input leg
 of $\Ga$ to its first non-special vertex   (if there are any). Consider next the filtration of $\widehat{\sP}^+$ by the number of vertices which are {\em not}\, very special. The associated graded is then isomorphic to
$A\ot C$, where $A$ is the trivial complex and the complex $(C=\oplus_{k\geq 0} C^k, d)$ is defined as follows:
 $C^k$ is  the one-dimensional vector space spanned  by the following graph with $k$ vertices ,
$$
\xy
(0,0)*{}="0",
(4,0)*{\bullet}="1",
(8,0)*{\bullet}="2",
(12,0)*{}="3",
(15,0)*{\ldots},
(18,0)*{}="4",
(22,0)*{\bullet}="5",
(26,0)*{\bullet}="6",
(30,0)*{}="7",
\ar @{-} "1";"2" <0pt>
\ar @{-} "2";"3" <0pt>
\ar @{-} "0";"1" <0pt>
\ar @{-} "4";"5" <0pt>
\ar @{-} "5";"6" <0pt>
\ar @{-} "6";"7" <0pt>
\endxy
$$
 and the differential $d$ is given by the canonical isomorphism $d: C^k \rar C^{k+1}$ for any
 $k=0,2,4,\ldots $, and by zero map for $k=1,3,\ldots$, e.g.
 $
 d(\xy
(0,0)*{}="0",
(5,0)*{}="1",
\ar @{-} "0";"1" <0pt>
\endxy)=
\xy
(0,0)*{}="0",
(4,0)*{\bullet}="1",
(8,0)*{}="2",
\ar @{-} "1";"0" <0pt>
\ar @{-} "2";"1" <0pt>
 \endxy
 $.
It is clear that $H^\bullet(C)=0$ making the cohomology of $\widehat{\sP}^+$ trivial (as all the spectral sequences involved in the argument converge due to the completeness of the filtration by the number of vertices).

\subsection{Polydifferential operators.} For an arbitrary vector space $V$ we consider a subspace,
$$
\Hom_{poly}(\f_V^{\ot \bullet}, \f_V^{\ot \bullet}) \subset
 \Hom(\f_V^{\ot \bullet}, \f_V^{\ot \bullet}),
$$
spanned by polydifferential operators,
$$
\Ba{rccc}
\Phi: & \f_V^{\ot m} &  \lon &  \f_V^{\ot n}\\
& f_1\ot \ldots \ot f_m & \lon & \Gamma(f_1, \ldots, f_m),
\Ea
$$
of the form,
\Beq\label{polyoperator}
\Phi(f_1, \ldots, f_m)= 
x^{J_1} \ot \ldots \ot x^{J_n}\cdot \Delta^{n-1}\left( \frac{\p^{|I_1|}f_1}{\p x^{I_1}}\right)\cdot\ldots \cdot
 \Delta^{n-1}\left( \frac{\p^{|I_m|}f_m}{\p x^{I_m}}\right).
\Eeq
where $x^i$ are some linear coordinates in $V$, $I$ and $J$ stand for multi-indices, say, $i_1i_2\ldots i_p$,
and $j_1j_2\ldots j_q$,
$$
x^J:= x^{j_1}x^{j_2}\ldots x^{j_q}, \ \ \ \ \ \
 \frac{\p^{|I|}}{\p x^I}:= \frac{\p^{p}}{\p x^{i_1}\p x^{i_2}\ldots \p x^{i_p}}.
$$
and $\Delta$ is the standard graded cocommutative coproduct in $\f_V$.
\sip

As $\f_V$ is naturally a bialgebra, there is  an associated Gerstenhaber-Schack complex \cite{GS},
\Beq\label{GS-complex}
\left(\bigoplus_{m,n\geq 1} \Hom(\f_V^{\ot m}, \f_V^{\ot n})[m+n-2], d_{\fg\fs}\right),
\Eeq
with the differential, $d_{\fg\fs}$,  given as in \S~\ref{fact-GS}.

\begin{lemma}\label{GS-poly}
The  subspace,
\Beq\label{4: inclusion of polydiff}
\bigoplus_{m,n\geq 1} \Hom_{poly}(\f_V^{\ot m}, \f_V^{\ot n})[m+n-2]\subset
\bigoplus_{m,n\geq 1} \Hom(\f_V^{\ot m}, \f_V^{\ot n})[m+n-2],
\Eeq
 is a subcomplex of the Gerstenhaber-Schack complex of the bialgebra $\f_V$, with the
Gerstenhaber-Schack differential given explicitly on polydifferential operators (\ref{polyoperator}) by
\Beqrn
d_{\fg\fs}\Phi &=&
\sum_{i=1}^{n}(-1)^{i+1}
x^{J_1} \ot  \mbox{...} \ot \bar{\Delta}(x^{J_i})\ot   \mbox{...}  \ot x^{J_n}\cdot \Delta^{n}
\left( \frac{\p^{|I_1|}}{\p x^{I_1}}\right)\cdot\ldots \cdot
 \Delta^{n}\left( \frac{\p^{|I_m|}}{\p x^{I_m}}\right)
\nonumber \\
&&
+\ 
\sum_{i=1}^{m}(-1)^{i+1}\sum_{I_i=I_i'\sqcup I_i''}
\nonumber \\
&&
\ \ \ \
x^{J_1} \ot \mbox{...} \ot x^{J_m}\cdot \Delta^{n-1}\left( \frac{\p^{|I_1|}}{\p x^{I_1}}\right)\cdot\mbox{...} \cdot
 \Delta^{n-1}\left( \frac{\p^{|I'_i|}}{\p x^{I'_i}}\right)\cdot  \Delta^{n-1}\left( \frac{\p^{|I_i''|}}{\p x^{I''_i}}\right)
\cdot\mbox{...}\cdot
 \Delta^{n-1}\left( \frac{\p^{|I_m|}}{\p x^{I_m}}\right),
\Eeqrn
where $\bar{\Delta}: \f_V\rar \f_V\ot \f_V$ is the reduced diagonal, and the second summation goes over all
possible splittings, ${I_i=I_i'\sqcup I_i''}$, into non-empty disjoint subsets.

\sip

Moreover, the inclusion (\ref{4: inclusion of polydiff}) is a quasi-isomorphism.
\end{lemma}

\noindent Proof is a straightforward calculation (full details are shown in \cite{Me3}).

\subsubsection{\bf Normalized polydifferential deformation complex} Note that polydifferential operators (\ref{polyoperator}) are allowed to have sets of multi-indices $I_i$ and $J_i$ with
zero cardinalities $|I_i|=0$ and/or $|J_i|=0$ (so that the standard multiplication and comultiplication in $\f_V$
belong to the class of polydifferential operators and it makes sense to consider a $L_\infty$
algebra $\Def(\sB_\infty^+\stackrel{\rho_0}{\rar}
\mathsf{End^{poly}_{\f_{\mathit V}}})$). However, the complex in the Lemma above  is quasi-isomorphic (cf.\ \cite{Lo})
to its {\em normalized}\, version, i.e.\ to the one which is spanned by operators which vanish if at least one input
function is constant, and never take constant value in each tensor factor. Therefore
we can assume from now on that the $L_\infty$-algebra
$$
\mathsf{poly}(\f_V,\f_V):=\Def(\sB_\infty^+\stackrel{\rho_0}{\rar}
\mathsf{End^{poly}_{\f_{\mathit V}}})
$$
 is spanned
by the {\em normalized}\, polydifferential operators (\ref{polyoperator}), i.e.\ the ones
which have a $|I_i|\geq 1$ and  $|J_j|\geq 1$ for all $i\in [m]$ and $j\in [n]$. Let us denote
by
$$
\mu_k: \wedge^k \mathsf{poly}(\f_V,\f_V) \lon \mathsf{poly}(\f_V,\f_V)[2-k]
$$
the $L_\infty$ structure in $\Def(\sB_\infty^+\stackrel{\rho_0}{\rar}
\mathsf{End^{poly}_{\f_{\mathit V}}})$ induced by the differential in $\sB_\infty^+$
and the standard representation $\rho_0: \sB\rar \mathsf{End^{poly}_{\f_{\mathit V}}}$.
The set of the associated Maurer-Cartan elements $\Ga\in  \mathsf{poly}(\f_V,\f_V)$,
$$
\mu_1(\Gamma) + \frac{1}{2!}\mu_2(\Gamma, \Gamma) +
\frac{1}{3!}\mu_3(\Gamma, \Gamma, \Gamma) + \ldots =0,\ \ \  |\Ga|=1,
$$
is in one-to-one correspondence with polydifferential representations of the form
$\rho_0+ \Ga:  \sB\rar \mathsf{End^{poly}_{\f_{\mathit V}}}$, where $\Ga$ is given by the {\em normalized}\, operators; we call such a representation $\rho_0+ \Ga$ a {\em normalized polydifferential}\, representation of $\sB_\infty^+$ in $\f_V$.

\subsection{Prop of quantum strongly homotopy bialgebra structures}
In his section we construct a dg free prop, $(\DefQ^+, \sd)$, whose representations,
$
                 \DefQ^+ \lon \mathsf{End}_V,
$
in a dg vector space $(V,d)$ are in one-to-one correspondence with normalized polydifferential
 representations of $\sB_\infty^+$ in $\f_V$, or, equivalently
 with Maurer-Cartan elements, $\Gamma$,
of the $L_\infty$-algebra $(\mathsf{poly}(\f_V,\f_V), \mu_\bullet)$.

\sip

An arbitrary degree 1 element $\Gamma$ in $\mathsf{poly}(\f_V,\f_V)$ has the form,
\Beq\label{Gamma}
\Gamma=\sum_{m,n\geq 1} \sum_{I_1, \ldots I_n\atop J_1,\ldots J_m}\Gamma^{I_1\ldots I_m}_{J_1\ldots J_n}
x^{J_1} \ot \ldots \ot x^{J_n}\cdot \Delta^{n-1}\left( \frac{\p^{|I_1|}}{\p x^{I_1}}\right)\cdot\ldots \cdot
 \Delta^{n-1}\left( \frac{\p^{|I_1|}}{\p x^{I_m}}\right),
\Eeq
where $\Gamma^{I_1\ldots I_m}_{J_1\ldots J_n}$ are the coordinate components of some linear map,
$$
V^{\odot{|J_1|}}\ot\ldots \ot  V^{\odot{|J_n|}} \lon   V^{\odot{|I_1|}}\ot\ldots \ot  V^{\odot{|I_m|}},
$$
 of degree $3-m-n$.
Thus there is a one-to-one correspondence between degree 1 elements $\Gamma\in \mathsf{poly}(\f_V,\f_V)$
and representations in $\mathsf{End}_V$ (not in  $\mathsf{End}_{\f_V}$!) of the free prop,
$\DefQ^+:= \Gamma\langle \sB \rangle$, generated by the $\bS$-bimodule,
$$
\sC=\{\sB(p,q):=\sC(m)\ot \sC(n)[p+q-3]\}_{p,q\geq 1},
$$
where
$$
\sB(p):=\bigoplus_{m\geq 1}\bigoplus_{ [p]=I_1\sqcup
\ldots \sqcup I_m \atop |I_1|,\ldots,|I_k|\geq 1} {\rm
Ind}^{\bS_m}_{\bS_{|I_1|}\times \ldots \times \bS_{|I_m|}} {\bf
1}_{|I_1|}\ot \ldots \ot {\bf 1}_{|I_m|}.
$$
Thus the generators of $\DefQ$ can be pictorially represented by degree $3-m-n$
directed planar corollas of the form,
\Beq\label{corolla}
\resizebox{17mm}{!}{ \xy
 <0mm,0mm>*{\mbox{$\xy *=<20mm,3mm>\txt{}*\frm{-}\endxy$}};<0mm,0mm>*{}**@{},
  <-10mm,1.5mm>*{};<-12mm,7mm>*{}**@{-},
  <-10mm,1.5mm>*{};<-11mm,7mm>*{}**@{-},
  <-10mm,1.5mm>*{};<-9.5mm,6mm>*{}**@{-},
  <-10mm,1.5mm>*{};<-8mm,7mm>*{}**@{-},
 <-10mm,1.5mm>*{};<-9.5mm,6.6mm>*{.\hspace{-0.4mm}.\hspace{-0.4mm}.}**@{},
 <0mm,0mm>*{};<-6.5mm,3.6mm>*{.\hspace{-0.1mm}.\hspace{-0.1mm}.}**@{},
  <-3mm,1.5mm>*{};<-5mm,7mm>*{}**@{-},
  <-3mm,1.5mm>*{};<-4mm,7mm>*{}**@{-},
  <-3mm,1.5mm>*{};<-2.5mm,6mm>*{}**@{-},
  <-3mm,1.5mm>*{};<-1mm,7mm>*{}**@{-},
 <-3mm,1.5mm>*{};<-2.5mm,6.6mm>*{.\hspace{-0.4mm}.\hspace{-0.4mm}.}**@{},
  <2mm,1.5mm>*{};<0mm,7mm>*{}**@{-},
  <2mm,1.5mm>*{};<1mm,7mm>*{}**@{-},
  <2mm,1.5mm>*{};<2.5mm,6mm>*{}**@{-},
  <2mm,1.5mm>*{};<4mm,7mm>*{}**@{-},
 <2mm,1.5mm>*{};<2.5mm,6.6mm>*{.\hspace{-0.4mm}.\hspace{-0.4mm}.}**@{},
 <0mm,0mm>*{};<6mm,3.6mm>*{.\hspace{-0.1mm}.\hspace{-0.1mm}.}**@{},
<10mm,1.5mm>*{};<8mm,7mm>*{}**@{-},
  <10mm,1.5mm>*{};<9mm,7mm>*{}**@{-},
  <10mm,1.5mm>*{};<10.5mm,6mm>*{}**@{-},
  <10mm,1.5mm>*{};<12mm,7mm>*{}**@{-},
 <10mm,1.5mm>*{};<10.5mm,6.6mm>*{.\hspace{-0.4mm}.\hspace{-0.4mm}.}**@{},
%
<0mm,0mm>*{};<-9.5mm,8.2mm>*{^{I_{ 1}}}**@{},
<0mm,0mm>*{};<-3mm,8.2mm>*{^{I_{ i}}}**@{},
<0mm,0mm>*{};<2mm,8.2mm>*{^{I_{ i+1}}}**@{},
<0mm,0mm>*{};<10mm,8.2mm>*{^{I_{ m}}}**@{},
<-10mm,-1.5mm>*{};<-12mm,-7mm>*{}**@{-},
  <-10mm,-1.5mm>*{};<-11mm,-7mm>*{}**@{-},
  <-10mm,-1.5mm>*{};<-9.5mm,-6mm>*{}**@{-},
  <-10mm,-1.5mm>*{};<-8mm,-7mm>*{}**@{-},
 <-10mm,-1.5mm>*{};<-9.5mm,-6.6mm>*{.\hspace{-0.4mm}.\hspace{-0.4mm}.}**@{},
 <0mm,0mm>*{};<-6.5mm,-3.6mm>*{.\hspace{-0.1mm}.\hspace{-0.1mm}.}**@{},
  <-3mm,-1.5mm>*{};<-5mm,-7mm>*{}**@{-},
  <-3mm,-1.5mm>*{};<-4mm,-7mm>*{}**@{-},
  <-3mm,-1.5mm>*{};<-2.5mm,-6mm>*{}**@{-},
  <-3mm,-1.5mm>*{};<-1mm,-7mm>*{}**@{-},
 <-3mm,-1.5mm>*{};<-2.5mm,-6.6mm>*{.\hspace{-0.4mm}.\hspace{-0.4mm}.}**@{},
  <2mm,-1.5mm>*{};<0mm,-7mm>*{}**@{-},
  <2mm,-1.5mm>*{};<1mm,-7mm>*{}**@{-},
  <2mm,-1.5mm>*{};<2.5mm,-6mm>*{}**@{-},
  <2mm,-1.5mm>*{};<4mm,-7mm>*{}**@{-},
 <2mm,-1.5mm>*{};<2.5mm,-6.6mm>*{.\hspace{-0.4mm}.\hspace{-0.4mm}.}**@{},
 <0mm,0mm>*{};<6mm,-3.6mm>*{.\hspace{-0.1mm}.\hspace{-0.1mm}.}**@{},
<10mm,-1.5mm>*{};<8mm,-7mm>*{}**@{-},
  <10mm,-1.5mm>*{};<9mm,-7mm>*{}**@{-},
  <10mm,-1.5mm>*{};<10.5mm,-6mm>*{}**@{-},
  <10mm,-1.5mm>*{};<12mm,-7mm>*{}**@{-},
 <10mm,-1.5mm>*{};<10.5mm,-6.6mm>*{.\hspace{-0.4mm}.\hspace{-0.4mm}.}**@{},
%
<0mm,0mm>*{};<-9.5mm,-9.2mm>*{^{J_{ 1}}}**@{},
<0mm,0mm>*{};<-3mm,-9.2mm>*{^{J_{j}}}**@{},
<0mm,0mm>*{};<2mm,-9.2mm>*{^{J_{ j+1}}}**@{},
<0mm,0mm>*{};<10mm,-9.2mm>*{^{J_{ n}}}**@{},
\endxy}
\ \ \ \simeq \ \ \ \Gamma^{I_1,\ldots, I_m}_{J_1,\ldots, J_n}
\Eeq
where
\Bi
\item the input legs  are labeled by the set $[q]=[1,2,\ldots,q]$
partitioned into $n$ disjoint non-empty subsets,
$$
[q]=J_1\sqcup \ldots\sqcup J_j\sqcup J_{j+1}\sqcup\ldots \sqcup J_n,
$$
and legs in each $J_j$-bunch are symmetric  (so that it does not matter
how labels from the set $J_j$ are distributed over legs in $j$-th bunch);
\item the output legs  are labeled by the set $[p]$
partitioned into $m$ disjoint non-empty subsets,
$$
[p]=I_1\sqcup \ldots\sqcup I_i\sqcup I_{i+1}\sqcup\ldots \sqcup I_k,
$$
and legs in each $I_i$-bunch are symmetric.
\Ei

%
%


Our next step is to encode of the $L_\infty$-structure $\{\mu_{\bullet}\}$ in $\mathsf{poly}(\f_V,\f_V)$
into a differential $d$ in $\DefQ^+$. The idea is simple \cite{Me2}: we replace coefficients
$\Gamma^{I_1,\ldots, I_m}_{J_1,\ldots, J_n}$ in the formal series (\ref{Gamma}) by the associated generating
corollas (\ref{corolla}) getting thereby a well-defined element
$$
\bar{\Gamma}:=
\sum_{m,n\geq 1} \sum_{I_1, \ldots I_n\atop J_1,\ldots J_m}
\resizebox{17mm}{!}{ \xy
 <0mm,0mm>*{\mbox{$\xy *=<20mm,3mm>\txt{}*\frm{-}\endxy$}};<0mm,0mm>*{}**@{},
  <-10mm,1.5mm>*{};<-12mm,7mm>*{}**@{-},
  <-10mm,1.5mm>*{};<-11mm,7mm>*{}**@{-},
  <-10mm,1.5mm>*{};<-9.5mm,6mm>*{}**@{-},
  <-10mm,1.5mm>*{};<-8mm,7mm>*{}**@{-},
 <-10mm,1.5mm>*{};<-9.5mm,6.6mm>*{.\hspace{-0.4mm}.\hspace{-0.4mm}.}**@{},
 <0mm,0mm>*{};<-6.5mm,3.6mm>*{.\hspace{-0.1mm}.\hspace{-0.1mm}.}**@{},
  <-3mm,1.5mm>*{};<-5mm,7mm>*{}**@{-},
  <-3mm,1.5mm>*{};<-4mm,7mm>*{}**@{-},
  <-3mm,1.5mm>*{};<-2.5mm,6mm>*{}**@{-},
  <-3mm,1.5mm>*{};<-1mm,7mm>*{}**@{-},
 <-3mm,1.5mm>*{};<-2.5mm,6.6mm>*{.\hspace{-0.4mm}.\hspace{-0.4mm}.}**@{},
  <2mm,1.5mm>*{};<0mm,7mm>*{}**@{-},
  <2mm,1.5mm>*{};<1mm,7mm>*{}**@{-},
  <2mm,1.5mm>*{};<2.5mm,6mm>*{}**@{-},
  <2mm,1.5mm>*{};<4mm,7mm>*{}**@{-},
 <2mm,1.5mm>*{};<2.5mm,6.6mm>*{.\hspace{-0.4mm}.\hspace{-0.4mm}.}**@{},
 <0mm,0mm>*{};<6mm,3.6mm>*{.\hspace{-0.1mm}.\hspace{-0.1mm}.}**@{},
<10mm,1.5mm>*{};<8mm,7mm>*{}**@{-},
  <10mm,1.5mm>*{};<9mm,7mm>*{}**@{-},
  <10mm,1.5mm>*{};<10.5mm,6mm>*{}**@{-},
  <10mm,1.5mm>*{};<12mm,7mm>*{}**@{-},
 <10mm,1.5mm>*{};<10.5mm,6.6mm>*{.\hspace{-0.4mm}.\hspace{-0.4mm}.}**@{},
%
<0mm,0mm>*{};<-9.5mm,8.2mm>*{^{I_{ 1}}}**@{},
<0mm,0mm>*{};<-3mm,8.2mm>*{^{I_{ i}}}**@{},
<0mm,0mm>*{};<2mm,8.2mm>*{^{I_{ i+1}}}**@{},
<0mm,0mm>*{};<10mm,8.2mm>*{^{I_{ m}}}**@{},
<-10mm,-1.5mm>*{};<-12mm,-7mm>*{}**@{-},
  <-10mm,-1.5mm>*{};<-11mm,-7mm>*{}**@{-},
  <-10mm,-1.5mm>*{};<-9.5mm,-6mm>*{}**@{-},
  <-10mm,-1.5mm>*{};<-8mm,-7mm>*{}**@{-},
 <-10mm,-1.5mm>*{};<-9.5mm,-6.6mm>*{.\hspace{-0.4mm}.\hspace{-0.4mm}.}**@{},
 <0mm,0mm>*{};<-6.5mm,-3.6mm>*{.\hspace{-0.1mm}.\hspace{-0.1mm}.}**@{},
  <-3mm,-1.5mm>*{};<-5mm,-7mm>*{}**@{-},
  <-3mm,-1.5mm>*{};<-4mm,-7mm>*{}**@{-},
  <-3mm,-1.5mm>*{};<-2.5mm,-6mm>*{}**@{-},
  <-3mm,-1.5mm>*{};<-1mm,-7mm>*{}**@{-},
 <-3mm,-1.5mm>*{};<-2.5mm,-6.6mm>*{.\hspace{-0.4mm}.\hspace{-0.4mm}.}**@{},
  <2mm,-1.5mm>*{};<0mm,-7mm>*{}**@{-},
  <2mm,-1.5mm>*{};<1mm,-7mm>*{}**@{-},
  <2mm,-1.5mm>*{};<2.5mm,-6mm>*{}**@{-},
  <2mm,-1.5mm>*{};<4mm,-7mm>*{}**@{-},
 <2mm,-1.5mm>*{};<2.5mm,-6.6mm>*{.\hspace{-0.4mm}.\hspace{-0.4mm}.}**@{},
 <0mm,0mm>*{};<6mm,-3.6mm>*{.\hspace{-0.1mm}.\hspace{-0.1mm}.}**@{},
<10mm,-1.5mm>*{};<8mm,-7mm>*{}**@{-},
  <10mm,-1.5mm>*{};<9mm,-7mm>*{}**@{-},
  <10mm,-1.5mm>*{};<10.5mm,-6mm>*{}**@{-},
  <10mm,-1.5mm>*{};<12mm,-7mm>*{}**@{-},
 <10mm,-1.5mm>*{};<10.5mm,-6.6mm>*{.\hspace{-0.4mm}.\hspace{-0.4mm}.}**@{},
%
<0mm,0mm>*{};<-9.5mm,-9.2mm>*{^{J_{ 1}}}**@{},
<0mm,0mm>*{};<-3mm,-9.2mm>*{^{J_{j}}}**@{},
<0mm,0mm>*{};<2mm,-9.2mm>*{^{J_{ j+1}}}**@{},
<0mm,0mm>*{};<10mm,-9.2mm>*{^{J_{ n}}}**@{},
\endxy}
\ot
x^{J_1} \ot \ldots \ot x^{J_n}\cdot \Delta^{n-1}\left( \frac{\p^{|I_1|}}{\p x^{I_1}}\right)\cdot\ldots \cdot
 \Delta^{n-1}\left( \frac{\p^{|I_1|}}{\p x^{I_m}}\right),
$$
in the tensor product $\DefQ^+\ot \mathsf{poly}(\f_V,\f_V)$. We extend $L_\infty$-operations
 $\{\mu_\bullet\}$ in $\mathsf{poly}(\f_V,\f_V)$  to the tensor product
 $\DefQ^+\ot \fg\fs^\bullet(\f_V)$ by replacing contractions of the coefficients $\Gamma^{I_1\ldots I_m}_{J_1\ldots J_n}$ by the associated  prop compositions in $\DefQ^+$ and then define a derivation, $\sd$, in $\DefQ^+$ by comparing coefficients
of the monomials,
$ x^{J_1} \ot \ldots \ot x^{J_n}\cdot \Delta^{n-1}\left( \frac{\p^{|I_1|}}{\p x^{I_1}}\right)\cdot\ldots \cdot
 \Delta^{n-1}\left( \frac{\p^{|I_1|}}{\p x^{I_m}}\right)$, in the equality,
\Beq\label{sd}
\sd\bar{\Gamma}= \sum_{n=1}^\infty\frac{1}{n!} \mu_n(\bar{\Gamma}, \bar{\Gamma}, \ldots, \bar{\Gamma})
\Eeq
The quadratic $L_\infty$ equations which hold for $\{\mu_\bullet\}$ imply then
$\sd^2=0$ \footnote{See also \S 2.5 and \S 2.5.1 in \cite{Me2} for a detailed explanation of this implication
 in the case when $\mu_\bullet$ is a dg Lie algebra;
below in \S~\ref{examples} we illustrate definition
(\ref{sd}) with some explicit computations.}.
 This completes the construction of the dg prop $(\DefQ^+,\sd)$ we  call from now on the dg {\em prop of quantum strongly homotopy bialgebra structures}.

\subsection{Examples.}\label{examples} The generators
$
\resizebox{16mm}{!}{
\xy
 <0mm,0mm>*{\mbox{$\xy *=<20mm,3mm>\txt{}*\frm{-}\endxy$}};<0mm,0mm>*{}**@{},
  <-10mm,1.5mm>*{};<-12mm,7mm>*{}**@{-},
  <-10mm,1.5mm>*{};<-11mm,7mm>*{}**@{-},
  <-10mm,1.5mm>*{};<-9.5mm,6mm>*{}**@{-},
  <-10mm,1.5mm>*{};<-8mm,7mm>*{}**@{-},
 <-10mm,1.5mm>*{};<-9.5mm,6.6mm>*{.\hspace{-0.4mm}.\hspace{-0.4mm}.}**@{},
 <0mm,0mm>*{};<-6.5mm,3.6mm>*{.\hspace{-0.1mm}.\hspace{-0.1mm}.}**@{},
  <-3mm,1.5mm>*{};<-5mm,7mm>*{}**@{-},
  <-3mm,1.5mm>*{};<-4mm,7mm>*{}**@{-},
  <-3mm,1.5mm>*{};<-2.5mm,6mm>*{}**@{-},
  <-3mm,1.5mm>*{};<-1mm,7mm>*{}**@{-},
 <-3mm,1.5mm>*{};<-2.5mm,6.6mm>*{.\hspace{-0.4mm}.\hspace{-0.4mm}.}**@{},
  <2mm,1.5mm>*{};<0mm,7mm>*{}**@{-},
  <2mm,1.5mm>*{};<1mm,7mm>*{}**@{-},
  <2mm,1.5mm>*{};<2.5mm,6mm>*{}**@{-},
  <2mm,1.5mm>*{};<4mm,7mm>*{}**@{-},
 <2mm,1.5mm>*{};<2.5mm,6.6mm>*{.\hspace{-0.4mm}.\hspace{-0.4mm}.}**@{},
 <0mm,0mm>*{};<6mm,3.6mm>*{.\hspace{-0.1mm}.\hspace{-0.1mm}.}**@{},
<10mm,1.5mm>*{};<8mm,7mm>*{}**@{-},
  <10mm,1.5mm>*{};<9mm,7mm>*{}**@{-},
  <10mm,1.5mm>*{};<10.5mm,6mm>*{}**@{-},
  <10mm,1.5mm>*{};<12mm,7mm>*{}**@{-},
 <10mm,1.5mm>*{};<10.5mm,6.6mm>*{.\hspace{-0.4mm}.\hspace{-0.4mm}.}**@{},
%
<0mm,0mm>*{};<-9.5mm,8.2mm>*{^{I_{ 1}}}**@{},
<0mm,0mm>*{};<-3mm,8.2mm>*{^{I_{ i}}}**@{},
<0mm,0mm>*{};<2mm,8.2mm>*{^{I_{ i+1}}}**@{},
<0mm,0mm>*{};<10mm,8.2mm>*{^{I_{ m}}}**@{},
<-10mm,-1.5mm>*{};<-12mm,-7mm>*{}**@{-},
  <-10mm,-1.5mm>*{};<-11mm,-7mm>*{}**@{-},
  <-10mm,-1.5mm>*{};<-9.5mm,-6mm>*{}**@{-},
  <-10mm,-1.5mm>*{};<-8mm,-7mm>*{}**@{-},
 <-10mm,-1.5mm>*{};<-9.5mm,-6.6mm>*{.\hspace{-0.4mm}.\hspace{-0.4mm}.}**@{},
 <0mm,0mm>*{};<-6.5mm,-3.6mm>*{.\hspace{-0.1mm}.\hspace{-0.1mm}.}**@{},
  <-3mm,-1.5mm>*{};<-5mm,-7mm>*{}**@{-},
  <-3mm,-1.5mm>*{};<-4mm,-7mm>*{}**@{-},
  <-3mm,-1.5mm>*{};<-2.5mm,-6mm>*{}**@{-},
  <-3mm,-1.5mm>*{};<-1mm,-7mm>*{}**@{-},
 <-3mm,-1.5mm>*{};<-2.5mm,-6.6mm>*{.\hspace{-0.4mm}.\hspace{-0.4mm}.}**@{},
  <2mm,-1.5mm>*{};<0mm,-7mm>*{}**@{-},
  <2mm,-1.5mm>*{};<1mm,-7mm>*{}**@{-},
  <2mm,-1.5mm>*{};<2.5mm,-6mm>*{}**@{-},
  <2mm,-1.5mm>*{};<4mm,-7mm>*{}**@{-},
 <2mm,-1.5mm>*{};<2.5mm,-6.6mm>*{.\hspace{-0.4mm}.\hspace{-0.4mm}.}**@{},
 <0mm,0mm>*{};<6mm,-3.6mm>*{.\hspace{-0.1mm}.\hspace{-0.1mm}.}**@{},
<10mm,-1.5mm>*{};<8mm,-7mm>*{}**@{-},
  <10mm,-1.5mm>*{};<9mm,-7mm>*{}**@{-},
  <10mm,-1.5mm>*{};<10.5mm,-6mm>*{}**@{-},
  <10mm,-1.5mm>*{};<12mm,-7mm>*{}**@{-},
 <10mm,-1.5mm>*{};<10.5mm,-6.6mm>*{.\hspace{-0.4mm}.\hspace{-0.4mm}.}**@{},
%
<0mm,0mm>*{};<-9.5mm,-9.2mm>*{^{J_{ 1}}}**@{},
<0mm,0mm>*{};<-3mm,-9.2mm>*{^{J_{ j}}}**@{},
<0mm,0mm>*{};<2mm,-9.2mm>*{^{J_{ j+1}}}**@{},
<0mm,0mm>*{};<10mm,-9.2mm>*{^{J_{ n}}}**@{},
\endxy}
$
are kind of ``thickenings" of the corresponding generators
$
\Ba{c}\resizebox{11mm}{!}{\begin{xy}
 <0mm,0mm>*{\bullet};<0mm,0mm>*{}**@{},
 <0mm,0mm>*{};<-8mm,5mm>*{}**@{-},
 <0mm,0mm>*{};<-4.5mm,5mm>*{}**@{-},
 <0mm,0mm>*{};<-1mm,5mm>*{\ldots}**@{},
 <0mm,0mm>*{};<4.5mm,5mm>*{}**@{-},
 <0mm,0mm>*{};<8mm,5mm>*{}**@{-},
   <0mm,0mm>*{};<-8.5mm,5.5mm>*{^1}**@{},
   <0mm,0mm>*{};<-5mm,5.5mm>*{^2}**@{},
   <0mm,0mm>*{};<4.5mm,5.5mm>*{^{m\hspace{-0.5mm}-\hspace{-0.5mm}1}}**@{},
   <0mm,0mm>*{};<9.0mm,5.5mm>*{^m}**@{},
 <0mm,0mm>*{};<-8mm,-5mm>*{}**@{-},
 <0mm,0mm>*{};<-4.5mm,-5mm>*{}**@{-},
 <0mm,0mm>*{};<-1mm,-5mm>*{\ldots}**@{},
 <0mm,0mm>*{};<4.5mm,-5mm>*{}**@{-},
 <0mm,0mm>*{};<8mm,-5mm>*{}**@{-},
   <0mm,0mm>*{};<-8.5mm,-6.9mm>*{^1}**@{},
   <0mm,0mm>*{};<-5mm,-6.9mm>*{^2}**@{},
   <0mm,0mm>*{};<4.5mm,-6.9mm>*{^{n\hspace{-0.5mm}-\hspace{-0.5mm}1}}**@{},
   <0mm,0mm>*{};<9.0mm,-6.9mm>*{^n}**@{},
 \end{xy}}\Ea
$
of $\hsB_\infty$ (mimicking thickening of $V$ into $\f_V$), and the differential $\sd$ in $\DefQ$ is also a kind
of ``thickening" (and twisting by $\rho_0$) of the differential $\hdelta$ in $\hsB_\infty$. Once an explicit
formula for $\hdelta$ is known it is very straightforward to compute
the associated expression for $\sd$ using (\ref{sd}). Below we show such expressions for those components
of $\hdelta$ which we know explicitly.

\subsubsection{} As $\mu_1$ is given explicitly by Lemma~\ref{GS-poly},
 we immediately get from
(\ref{sd})  the following expression for $\sd$ modulo
terms with number of vertices $\geq 2$,
\Beq\label{graph-GS}
\sd\left(
\resizebox{17mm}{!}{ \xy
 <0mm,0mm>*{\mbox{$\xy *=<20mm,3mm>\txt{}*\frm{-}\endxy$}};<0mm,0mm>*{}**@{},
  <-10mm,1.5mm>*{};<-12mm,7mm>*{}**@{-},
  <-10mm,1.5mm>*{};<-11mm,7mm>*{}**@{-},
  <-10mm,1.5mm>*{};<-9.5mm,6mm>*{}**@{-},
  <-10mm,1.5mm>*{};<-8mm,7mm>*{}**@{-},
 <-10mm,1.5mm>*{};<-9.5mm,6.6mm>*{.\hspace{-0.4mm}.\hspace{-0.4mm}.}**@{},
 <0mm,0mm>*{};<-6.5mm,3.6mm>*{.\hspace{-0.1mm}.\hspace{-0.1mm}.}**@{},
  <-3mm,1.5mm>*{};<-5mm,7mm>*{}**@{-},
  <-3mm,1.5mm>*{};<-4mm,7mm>*{}**@{-},
  <-3mm,1.5mm>*{};<-2.5mm,6mm>*{}**@{-},
  <-3mm,1.5mm>*{};<-1mm,7mm>*{}**@{-},
 <-3mm,1.5mm>*{};<-2.5mm,6.6mm>*{.\hspace{-0.4mm}.\hspace{-0.4mm}.}**@{},
  <2mm,1.5mm>*{};<0mm,7mm>*{}**@{-},
  <2mm,1.5mm>*{};<1mm,7mm>*{}**@{-},
  <2mm,1.5mm>*{};<2.5mm,6mm>*{}**@{-},
  <2mm,1.5mm>*{};<4mm,7mm>*{}**@{-},
 <2mm,1.5mm>*{};<2.5mm,6.6mm>*{.\hspace{-0.4mm}.\hspace{-0.4mm}.}**@{},
 <0mm,0mm>*{};<6mm,3.6mm>*{.\hspace{-0.1mm}.\hspace{-0.1mm}.}**@{},
<10mm,1.5mm>*{};<8mm,7mm>*{}**@{-},
  <10mm,1.5mm>*{};<9mm,7mm>*{}**@{-},
  <10mm,1.5mm>*{};<10.5mm,6mm>*{}**@{-},
  <10mm,1.5mm>*{};<12mm,7mm>*{}**@{-},
 <10mm,1.5mm>*{};<10.5mm,6.6mm>*{.\hspace{-0.4mm}.\hspace{-0.4mm}.}**@{},
%
<0mm,0mm>*{};<-9.5mm,8.2mm>*{^{I_{ 1}}}**@{},
<0mm,0mm>*{};<-3mm,8.2mm>*{^{I_{ i}}}**@{},
<0mm,0mm>*{};<2mm,8.2mm>*{^{I_{ i+1}}}**@{},
<0mm,0mm>*{};<10mm,8.2mm>*{^{I_{ m}}}**@{},
<-10mm,-1.5mm>*{};<-12mm,-7mm>*{}**@{-},
  <-10mm,-1.5mm>*{};<-11mm,-7mm>*{}**@{-},
  <-10mm,-1.5mm>*{};<-9.5mm,-6mm>*{}**@{-},
  <-10mm,-1.5mm>*{};<-8mm,-7mm>*{}**@{-},
 <-10mm,-1.5mm>*{};<-9.5mm,-6.6mm>*{.\hspace{-0.4mm}.\hspace{-0.4mm}.}**@{},
 <0mm,0mm>*{};<-6.5mm,-3.6mm>*{.\hspace{-0.1mm}.\hspace{-0.1mm}.}**@{},
  <-3mm,-1.5mm>*{};<-5mm,-7mm>*{}**@{-},
  <-3mm,-1.5mm>*{};<-4mm,-7mm>*{}**@{-},
  <-3mm,-1.5mm>*{};<-2.5mm,-6mm>*{}**@{-},
  <-3mm,-1.5mm>*{};<-1mm,-7mm>*{}**@{-},
 <-3mm,-1.5mm>*{};<-2.5mm,-6.6mm>*{.\hspace{-0.4mm}.\hspace{-0.4mm}.}**@{},
  <2mm,-1.5mm>*{};<0mm,-7mm>*{}**@{-},
  <2mm,-1.5mm>*{};<1mm,-7mm>*{}**@{-},
  <2mm,-1.5mm>*{};<2.5mm,-6mm>*{}**@{-},
  <2mm,-1.5mm>*{};<4mm,-7mm>*{}**@{-},
 <2mm,-1.5mm>*{};<2.5mm,-6.6mm>*{.\hspace{-0.4mm}.\hspace{-0.4mm}.}**@{},
 <0mm,0mm>*{};<6mm,-3.6mm>*{.\hspace{-0.1mm}.\hspace{-0.1mm}.}**@{},
<10mm,-1.5mm>*{};<8mm,-7mm>*{}**@{-},
  <10mm,-1.5mm>*{};<9mm,-7mm>*{}**@{-},
  <10mm,-1.5mm>*{};<10.5mm,-6mm>*{}**@{-},
  <10mm,-1.5mm>*{};<12mm,-7mm>*{}**@{-},
 <10mm,-1.5mm>*{};<10.5mm,-6.6mm>*{.\hspace{-0.4mm}.\hspace{-0.4mm}.}**@{},
%
<0mm,0mm>*{};<-9.5mm,-9.2mm>*{^{J_{ 1}}}**@{},
<0mm,0mm>*{};<-3mm,-9.2mm>*{^{J_{ j}}}**@{},
<0mm,0mm>*{};<2mm,-9.2mm>*{^{J_{ j+1}}}**@{},
<0mm,0mm>*{};<10mm,-9.2mm>*{^{J_{ n}}}**@{},
\endxy}
\right)
=\sum_{i=1}^{m-1}(-1)^i
\resizebox{17mm}{!}{ \xy
 <0mm,0mm>*{\mbox{$\xy *=<20mm,3mm>\txt{}*\frm{-}\endxy$}};<0mm,0mm>*{}**@{},
  <-10mm,1.5mm>*{};<-12mm,7mm>*{}**@{-},
  <-10mm,1.5mm>*{};<-11mm,7mm>*{}**@{-},
  <-10mm,1.5mm>*{};<-9.5mm,6mm>*{}**@{-},
  <-10mm,1.5mm>*{};<-8mm,7mm>*{}**@{-},
 <-10mm,1.5mm>*{};<-9.5mm,6.6mm>*{.\hspace{-0.4mm}.\hspace{-0.4mm}.}**@{},
 <0mm,0mm>*{};<-6.5mm,3.6mm>*{.\hspace{-0.1mm}.\hspace{-0.1mm}.}**@{},
  <0mm,1.5mm>*{};<-5mm,7mm>*{}**@{-},
  <0mm,1.5mm>*{};<-4mm,7mm>*{}**@{-},
  <0mm,1.5mm>*{};<-2.5mm,6mm>*{}**@{-},
  <0mm,1.5mm>*{};<-1.3mm,7mm>*{}**@{-},
 <0mm,1.5mm>*{};<-2.5mm,6.6mm>*{.\hspace{-0.4mm}.\hspace{-0.4mm}.}**@{},
  <0mm,1.5mm>*{};<1.3mm,7mm>*{}**@{-},
  <0mm,1.5mm>*{};<2.5mm,6mm>*{}**@{-},
  <0mm,1.5mm>*{};<4mm,7mm>*{}**@{-},
 <0mm,1.5mm>*{};<2.5mm,6.6mm>*{.\hspace{-0.4mm}.\hspace{-0.4mm}.}**@{},
 <0mm,0mm>*{};<6mm,3.6mm>*{.\hspace{-0.1mm}.\hspace{-0.1mm}.}**@{},
<10mm,1.5mm>*{};<8mm,7mm>*{}**@{-},
  <10mm,1.5mm>*{};<9mm,7mm>*{}**@{-},
  <10mm,1.5mm>*{};<10.5mm,6mm>*{}**@{-},
  <10mm,1.5mm>*{};<12mm,7mm>*{}**@{-},
 <10mm,1.5mm>*{};<10.5mm,6.6mm>*{.\hspace{-0.4mm}.\hspace{-0.4mm}.}**@{},
%
<0mm,0mm>*{};<-9.5mm,8.2mm>*{^{I_{ 1}}}**@{},
<0mm,0mm>*{};<-3mm,8.2mm>*{^{I_{ i}}}**@{},
<0mm,0mm>*{};<2mm,8.2mm>*{^{\sqcup I_{ i+1}}}**@{},
<0mm,0mm>*{};<10mm,8.2mm>*{^{I_{ m}}}**@{},
<-10mm,-1.5mm>*{};<-12mm,-7mm>*{}**@{-},
  <-10mm,-1.5mm>*{};<-11mm,-7mm>*{}**@{-},
  <-10mm,-1.5mm>*{};<-9.5mm,-6mm>*{}**@{-},
  <-10mm,-1.5mm>*{};<-8mm,-7mm>*{}**@{-},
 <-10mm,-1.5mm>*{};<-9.5mm,-6.6mm>*{.\hspace{-0.4mm}.\hspace{-0.4mm}.}**@{},
 <0mm,0mm>*{};<-6.5mm,-3.6mm>*{.\hspace{-0.1mm}.\hspace{-0.1mm}.}**@{},
  <-3mm,-1.5mm>*{};<-5mm,-7mm>*{}**@{-},
  <-3mm,-1.5mm>*{};<-4mm,-7mm>*{}**@{-},
  <-3mm,-1.5mm>*{};<-2.5mm,-6mm>*{}**@{-},
  <-3mm,-1.5mm>*{};<-1mm,-7mm>*{}**@{-},
 <-3mm,-1.5mm>*{};<-2.5mm,-6.6mm>*{.\hspace{-0.4mm}.\hspace{-0.4mm}.}**@{},
  <2mm,-1.5mm>*{};<0mm,-7mm>*{}**@{-},
  <2mm,-1.5mm>*{};<1mm,-7mm>*{}**@{-},
  <2mm,-1.5mm>*{};<2.5mm,-6mm>*{}**@{-},
  <2mm,-1.5mm>*{};<4mm,-7mm>*{}**@{-},
 <2mm,-1.5mm>*{};<2.5mm,-6.6mm>*{.\hspace{-0.4mm}.\hspace{-0.4mm}.}**@{},
 <0mm,0mm>*{};<6mm,-3.6mm>*{.\hspace{-0.1mm}.\hspace{-0.1mm}.}**@{},
<10mm,-1.5mm>*{};<8mm,-7mm>*{}**@{-},
  <10mm,-1.5mm>*{};<9mm,-7mm>*{}**@{-},
  <10mm,-1.5mm>*{};<10.5mm,-6mm>*{}**@{-},
  <10mm,-1.5mm>*{};<12mm,-7mm>*{}**@{-},
 <10mm,-1.5mm>*{};<10.5mm,-6.6mm>*{.\hspace{-0.4mm}.\hspace{-0.4mm}.}**@{},
%
<0mm,0mm>*{};<-9.5mm,-9.2mm>*{^{J_{ 1}}}**@{},
<0mm,0mm>*{};<-3mm,-9.2mm>*{^{J_{ j}}}**@{},
<0mm,0mm>*{};<2mm,-9.2mm>*{^{J_{ j+1}}}**@{},
<0mm,0mm>*{};<10mm,-9.2mm>*{^{J_{ n}}}**@{},
\endxy}
+
\sum_{j=1}^{n-1}(-1)^j
\resizebox{17mm}{!}{ \xy
 <0mm,0mm>*{\mbox{$\xy *=<20mm,3mm>\txt{}*\frm{-}\endxy$}};<0mm,0mm>*{}**@{},
  <-10mm,1.5mm>*{};<-12mm,7mm>*{}**@{-},
  <-10mm,1.5mm>*{};<-11mm,7mm>*{}**@{-},
  <-10mm,1.5mm>*{};<-9.5mm,6mm>*{}**@{-},
  <-10mm,1.5mm>*{};<-8mm,7mm>*{}**@{-},
 <-10mm,1.5mm>*{};<-9.5mm,6.6mm>*{.\hspace{-0.4mm}.\hspace{-0.4mm}.}**@{},
 <0mm,0mm>*{};<-6.5mm,3.6mm>*{.\hspace{-0.1mm}.\hspace{-0.1mm}.}**@{},
  <-3mm,1.5mm>*{};<-5mm,7mm>*{}**@{-},
  <-3mm,1.5mm>*{};<-4mm,7mm>*{}**@{-},
  <-3mm,1.5mm>*{};<-2.5mm,6mm>*{}**@{-},
  <-3mm,1.5mm>*{};<-1mm,7mm>*{}**@{-},
 <-3mm,1.5mm>*{};<-2.5mm,6.6mm>*{.\hspace{-0.4mm}.\hspace{-0.4mm}.}**@{},
  <2mm,1.5mm>*{};<0mm,7mm>*{}**@{-},
  <2mm,1.5mm>*{};<1mm,7mm>*{}**@{-},
  <2mm,1.5mm>*{};<2.5mm,6mm>*{}**@{-},
  <2mm,1.5mm>*{};<4mm,7mm>*{}**@{-},
 <2mm,1.5mm>*{};<2.5mm,6.6mm>*{.\hspace{-0.4mm}.\hspace{-0.4mm}.}**@{},
 <0mm,0mm>*{};<6mm,3.6mm>*{.\hspace{-0.1mm}.\hspace{-0.1mm}.}**@{},
<10mm,1.5mm>*{};<8mm,7mm>*{}**@{-},
  <10mm,1.5mm>*{};<9mm,7mm>*{}**@{-},
  <10mm,1.5mm>*{};<10.5mm,6mm>*{}**@{-},
  <10mm,1.5mm>*{};<12mm,7mm>*{}**@{-},
 <10mm,1.5mm>*{};<10.5mm,6.6mm>*{.\hspace{-0.4mm}.\hspace{-0.4mm}.}**@{},
%
<0mm,0mm>*{};<-9.5mm,8.2mm>*{^{I_{ 1}}}**@{},
<0mm,0mm>*{};<-3mm,8.2mm>*{^{I_{ i}}}**@{},
<0mm,0mm>*{};<2mm,8.2mm>*{^{I_{ i+1}}}**@{},
<0mm,0mm>*{};<10mm,8.2mm>*{^{I_{ m}}}**@{},
<-10mm,-1.5mm>*{};<-12mm,-7mm>*{}**@{-},
  <-10mm,-1.5mm>*{};<-11mm,-7mm>*{}**@{-},
  <-10mm,-1.5mm>*{};<-9.5mm,-6mm>*{}**@{-},
  <-10mm,-1.5mm>*{};<-8mm,-7mm>*{}**@{-},
 <-10mm,-1.5mm>*{};<-9.5mm,-6.6mm>*{.\hspace{-0.4mm}.\hspace{-0.4mm}.}**@{},
 <0mm,0mm>*{};<-6.5mm,-3.6mm>*{.\hspace{-0.1mm}.\hspace{-0.1mm}.}**@{},
  <0mm,-1.5mm>*{};<-5mm,-7mm>*{}**@{-},
  <0mm,-1.5mm>*{};<-4mm,-7mm>*{}**@{-},
  <0mm,-1.5mm>*{};<-2.5mm,-6mm>*{}**@{-},
  <0mm,-1.5mm>*{};<-1.3mm,-7mm>*{}**@{-},
 <0mm,-1.5mm>*{};<-2.5mm,-6.6mm>*{.\hspace{-0.4mm}.\hspace{-0.4mm}.}**@{},
  <0mm,-1.5mm>*{};<1.3mm,-7mm>*{}**@{-},
  <0mm,-1.5mm>*{};<2.5mm,-6mm>*{}**@{-},
  <0mm,-1.5mm>*{};<4mm,-7mm>*{}**@{-},
 <0mm,-1.5mm>*{};<2.5mm,-6.6mm>*{.\hspace{-0.4mm}.\hspace{-0.4mm}.}**@{},
 <0mm,0mm>*{};<6mm,-3.6mm>*{.\hspace{-0.1mm}.\hspace{-0.1mm}.}**@{},
<10mm,-1.5mm>*{};<8mm,-7mm>*{}**@{-},
  <10mm,-1.5mm>*{};<9mm,-7mm>*{}**@{-},
  <10mm,-1.5mm>*{};<10.5mm,-6mm>*{}**@{-},
  <10mm,-1.5mm>*{};<12mm,-7mm>*{}**@{-},
 <10mm,-1.5mm>*{};<10.5mm,-6.6mm>*{.\hspace{-0.4mm}.\hspace{-0.4mm}.}**@{},
%
<0mm,0mm>*{};<-9.5mm,-9.2mm>*{^{J_{ 1}}}**@{},
<0mm,0mm>*{};<-3mm,-9.2mm>*{^{J_{ j}}}**@{},
<0mm,0mm>*{};<2.2mm,-9.2mm>*{^{\sqcup J_{ j+1}}}**@{},
<0mm,0mm>*{};<10mm,-9.2mm>*{^{J_{ n}}}**@{},
\endxy}
+ O( 2).
\Eeq

\subsubsection{}\label{Hoch}
On the generators
$
\Ba{c}\resizebox{8mm}{!}{ \begin{xy}
<0mm,0mm>*{\bullet},
<0mm,-5mm>*{}**@{-},
<-5mm,5mm>*{}**@{-},
<-2mm,5mm>*{}**@{-},
<2mm,5mm>*{}**@{-},
<5mm,5mm>*{}**@{-},
<0mm,7mm>*{_{1\ \ \  \ldots\ \  \ n}},
\end{xy}}\Ea$, $n\geq 1$,
of $\hsB_\infty$ the differential is given by
$$
\hdelta
\resizebox{8mm}{!}{ \begin{xy}
<0mm,0mm>*{\bullet},
<0mm,-5mm>*{}**@{-},
<-5mm,5mm>*{}**@{-},
<-2mm,5mm>*{}**@{-},
<2mm,5mm>*{}**@{-},
<5mm,5mm>*{}**@{-},
<0mm,7mm>*{_{1\ \ \  \ldots\ \  \ n}},
\end{xy}}
=\sum_{i=0}^{n-1}\sum_{q=1}^{n-i}
(-1)^{i+l(n-i-q)+1}
\resizebox{21mm}{!}{ \begin{xy}
<0mm,0mm>*{\bullet},
<0mm,-5mm>*{}**@{-},
<4mm,7mm>*{^{1\ \dots\ i  \ \qquad\ \  i+q+1\ \dots \ n}},
<-14mm,5mm>*{}**@{-},
<-6mm,5mm>*{}**@{-},
<20mm,5mm>*{}**@{-},
<8mm,5mm>*{}**@{-},
<0mm,5mm>*{}**@{-},
<0mm,5mm>*{\bullet};
<-5mm,10mm>*{}**@{-},
<-2mm,10mm>*{}**@{-},
<2mm,10mm>*{}**@{-},
<5mm,10mm>*{}**@{-},
<0mm,12.5mm>*{_{i+1\ \dots\ i+q}},
\end{xy}}.
$$
Using (\ref{sd}) get after some tedious calculation (cf.\ \cite{Me2}) the following expression for the value of
$\sd$
on the associated ``thickened" generators of $\DefQ^+$,
$$
\sd\left(\resizebox{17mm}{!}{ \xy
 <0mm,0mm>*{\mbox{$\xy *=<20mm,3mm>\txt{}*\frm{-}\endxy$}};<0mm,0mm>*{}**@{},
  <-10mm,1.5mm>*{};<-12mm,7mm>*{}**@{-},
  <-10mm,1.5mm>*{};<-11mm,7mm>*{}**@{-},
  <-10mm,1.5mm>*{};<-9.5mm,6mm>*{}**@{-},
  <-10mm,1.5mm>*{};<-8mm,7mm>*{}**@{-},
 <-10mm,1.5mm>*{};<-9.5mm,6.6mm>*{.\hspace{-0.4mm}.\hspace{-0.4mm}.}**@{},
 <0mm,0mm>*{};<-6.4mm,3.6mm>*{.\hspace{-0.1mm}.\hspace{-0.1mm}.}**@{},
  <-3mm,1.5mm>*{};<-5mm,7mm>*{}**@{-},
  <-3mm,1.5mm>*{};<-4mm,7mm>*{}**@{-},
  <-3mm,1.5mm>*{};<-2.5mm,6mm>*{}**@{-},
  <-3mm,1.5mm>*{};<-1mm,7mm>*{}**@{-},
 <-3mm,1.5mm>*{};<-2.5mm,6.6mm>*{.\hspace{-0.4mm}.\hspace{-0.4mm}.}**@{},
  <2mm,1.5mm>*{};<0mm,7mm>*{}**@{-},
  <2mm,1.5mm>*{};<1mm,7mm>*{}**@{-},
  <2mm,1.5mm>*{};<2.5mm,6mm>*{}**@{-},
  <2mm,1.5mm>*{};<4mm,7mm>*{}**@{-},
 <2mm,1.5mm>*{};<2.5mm,6.6mm>*{.\hspace{-0.4mm}.\hspace{-0.4mm}.}**@{},
 <0mm,0mm>*{};<6mm,3.6mm>*{.\hspace{-0.1mm}.\hspace{-0.1mm}.}**@{},
<10mm,1.5mm>*{};<8mm,7mm>*{}**@{-},
  <10mm,1.5mm>*{};<9mm,7mm>*{}**@{-},
  <10mm,1.5mm>*{};<10.5mm,6mm>*{}**@{-},
  <10mm,1.5mm>*{};<12mm,7mm>*{}**@{-},
 <10mm,1.5mm>*{};<10.5mm,6.6mm>*{.\hspace{-0.4mm}.\hspace{-0.4mm}.}**@{},
%
<0mm,0mm>*{};<-9.5mm,8.2mm>*{^{I_{ 1}}}**@{},
<0mm,0mm>*{};<-3mm,8.2mm>*{^{I_{ i}}}**@{},
<0mm,0mm>*{};<2mm,8.2mm>*{^{I_{ i+1}}}**@{},
<0mm,0mm>*{};<10mm,8.2mm>*{^{I_{ n}}}**@{},
<0mm,-1.5mm>*{};<-0.9mm,-7mm>*{}**@{-},
  <0mm,-1.5mm>*{};<-2mm,-7mm>*{}**@{-},
  <0mm,-1.5mm>*{};<0.5mm,-6mm>*{}**@{-},
  <0mm,-1.5mm>*{};<2mm,-7mm>*{}**@{-},
 <0mm,-1.5mm>*{};<0.5mm,-6.6mm>*{.\hspace{-0.4mm}.\hspace{-0.4mm}.}**@{},
<0mm,0mm>*{};<0mm,-9.2mm>*{^{J}}**@{},
\endxy}
\right) = \sum_{i=1}^{n-1}(-1)^{i}
\resizebox{17mm}{!}{ \xy
 <0mm,0mm>*{\mbox{$\xy *=<20mm,3mm>\txt{}*\frm{-}\endxy$}};<0mm,0mm>*{}**@{},
  <-10mm,1.5mm>*{};<-12mm,7mm>*{}**@{-},
  <-10mm,1.5mm>*{};<-11mm,7mm>*{}**@{-},
  <-10mm,1.5mm>*{};<-9.5mm,6mm>*{}**@{-},
  <-10mm,1.5mm>*{};<-8mm,7mm>*{}**@{-},
 <-10mm,1.5mm>*{};<-9.5mm,6.6mm>*{.\hspace{-0.4mm}.\hspace{-0.4mm}.}**@{},
 <0mm,0mm>*{};<-5.5mm,3.6mm>*{.\hspace{-0.1mm}.\hspace{-0.1mm}.}**@{},
%
  <0mm,1.5mm>*{};<-4mm,7mm>*{}**@{-},
  <0mm,1.5mm>*{};<-2.0mm,6mm>*{}**@{-},
  <0mm,1.5mm>*{};<-1mm,7mm>*{}**@{-},
 <0mm,1.5mm>*{};<-2.3mm,6.6mm>*{.\hspace{-0.4mm}.\hspace{-0.4mm}.}**@{},
%
  <0mm,1.5mm>*{};<1mm,7mm>*{}**@{-},
  <0mm,1.5mm>*{};<2.0mm,6mm>*{}**@{-},
  <0mm,1.5mm>*{};<4mm,7mm>*{}**@{-},
 <0mm,1.5mm>*{};<2.3mm,6.6mm>*{.\hspace{-0.4mm}.\hspace{-0.4mm}.}**@{},
 <0mm,0mm>*{};<6mm,3.6mm>*{.\hspace{-0.1mm}.\hspace{-0.1mm}.}**@{},
<10mm,1.5mm>*{};<8mm,7mm>*{}**@{-},
  <10mm,1.5mm>*{};<9mm,7mm>*{}**@{-},
  <10mm,1.5mm>*{};<10.5mm,6mm>*{}**@{-},
  <10mm,1.5mm>*{};<12mm,7mm>*{}**@{-},
 <10mm,1.5mm>*{};<10.5mm,6.6mm>*{.\hspace{-0.4mm}.\hspace{-0.4mm}.}**@{},
%
<0mm,0mm>*{};<-9.5mm,8.2mm>*{^{I_{ 1}}}**@{},
<0mm,0mm>*{};<-2.5mm,8.2mm>*{{^{I_{ i}\sqcup}}}**@{},
<0mm,0mm>*{};<2.7mm,8.2mm>*{^{I_{ i+1}}}**@{},
<0mm,0mm>*{};<10mm,8.2mm>*{^{I_{ n}}}**@{},
%
<0mm,-1.5mm>*{};<-0.9mm,-7mm>*{}**@{-},
  <0mm,-1.5mm>*{};<-2mm,-7mm>*{}**@{-},
  <0mm,-1.5mm>*{};<0.5mm,-6mm>*{}**@{-},
  <0mm,-1.5mm>*{};<2mm,-7mm>*{}**@{-},
 <0mm,-1.5mm>*{};<0.5mm,-6.6mm>*{.\hspace{-0.4mm}.\hspace{-0.4mm}.}**@{},
<0mm,0mm>*{};<0mm,-9.2mm>*{^{J}}**@{},
\endxy}
\
 +\ \sum_{p+q=n+1\atop {0\leq i\leq p-1 \atop 0\leq s}}
\sum_{  {I_{i+1}=I_{i+1}'\sqcup I''_{i+1}\atop .......................}
\atop
{I_{i+q}=I_{i+q}'\sqcup I''_{i+q} \atop J=J_1\sqcup J_2}}
 \frac{(-1)^\var}{s!}
\
\resizebox{35mm}{!}{ \xy
 <19mm,0mm>*{\mbox{$\xy *=<58mm,3mm>\txt{}*\frm{-}\endxy$}};<0mm,0mm>*{}**@{},
  <-10mm,1.5mm>*{};<-12mm,7mm>*{}**@{-},
  <-10mm,1.5mm>*{};<-11mm,7mm>*{}**@{-},
  <-10mm,1.5mm>*{};<-9.5mm,6mm>*{}**@{-},
  <-10mm,1.5mm>*{};<-8mm,7mm>*{}**@{-},
 <-10mm,1.5mm>*{};<-9.5mm,6.6mm>*{.\hspace{-0.4mm}.\hspace{-0.4mm}.}**@{},
 <0mm,0mm>*{};<-6.5mm,3.6mm>*{.\hspace{-0.1mm}.\hspace{-0.1mm}.}**@{},
  <-3mm,1.5mm>*{};<-5mm,7mm>*{}**@{-},
  <-3mm,1.5mm>*{};<-4mm,7mm>*{}**@{-},
  <-3mm,1.5mm>*{};<-2.5mm,6mm>*{}**@{-},
  <-3mm,1.5mm>*{};<-1mm,7mm>*{}**@{-},
 <-3mm,1.5mm>*{};<-2.5mm,6.6mm>*{.\hspace{-0.4mm}.\hspace{-0.4mm}.}**@{},
  <10mm,1.5mm>*{};<0mm,7mm>*{}**@{-},
  <10mm,1.5mm>*{};<4mm,7mm>*{}**@{-},
<10mm,1.5mm>*{};<7.3mm,5.9mm>*{.\hspace{-0.0mm}.\hspace{-0.0mm}.}**@{},
   <10mm,1.5mm>*{};<3.8mm,6.0mm>*{}**@{-},
 <10mm,1.5mm>*{};<2.5mm,6.6mm>*{.\hspace{-0.4mm}.\hspace{-0.4mm}.}**@{},
%
%
<10mm,1.5mm>*{};<9mm,7mm>*{}**@{-},
  <10mm,1.5mm>*{};<10.5mm,6mm>*{}**@{-},
  <10mm,1.5mm>*{};<12mm,7mm>*{}**@{-},
 <10mm,1.5mm>*{};<10.5mm,6.6mm>*{.\hspace{-0.4mm}.\hspace{-0.4mm}.}**@{},
<0mm,0mm>*{};<-9.5mm,8.4mm>*{^{I_{1}}}**@{},
<0mm,0mm>*{};<-3mm,8.4mm>*{^{I_{ i}}}**@{},
<0mm,0mm>*{};<2mm,8.6mm>*{^{I_{ i+1}'}}**@{},
<0mm,0mm>*{};<10.5mm,8.6mm>*{^{I_{ i+q}'}}**@{},
<10mm,1.5mm>*{};<18mm,10mm>*{}**@{-},
<18mm,10mm>*{};<25mm,12mm>*{}**@{-},
<10mm,1.5mm>*{};<20.0mm,7mm>*{}**@{-},
<20mm,7mm>*{};<25mm,12mm>*{}**@{-},
<10mm,1.5mm>*{};<25mm,12mm>*{}**@{-},
<10mm,1.5mm>*{};<18mm,8.1mm>*{.\hspace{-0.4mm}.\hspace{-0.4mm}.}**@{},
<10mm,1.5mm>*{};<19mm,8.5mm>*{^s}**@{},
%
<0mm,0mm>*{};<-9.5mm,8.2mm>*{^{I_{ 1}}}**@{},
<0mm,0mm>*{};<-3mm,8.2mm>*{^{I_{ i}}}**@{},
<0mm,0mm>*{};<2mm,8.2mm>*{^{I_{ i+1}}}**@{},
<0mm,0mm>*{};<10mm,8.2mm>*{^{I_{ n}}}**@{},
<20mm,-1.5mm>*{};<19.1mm,-7mm>*{}**@{-},
  <20mm,-1.5mm>*{};<18mm,-7mm>*{}**@{-},
  <20mm,-1.5mm>*{};<20.5mm,-6mm>*{}**@{-},
  <20mm,-1.5mm>*{};<22mm,-7mm>*{}**@{-},
 <20mm,-1.5mm>*{};<20.5mm,-6.6mm>*{.\hspace{-0.4mm}.\hspace{-0.4mm}.}**@{},
<0mm,0mm>*{};<20mm,-9.2mm>*{^{J_1}}**@{},
<25mm,13.75mm>*{\mbox{$\xy *=<14mm,3mm>\txt{}*\frm{-}\endxy$}};
<0mm,0mm>*{}**@{},
 <18mm,15mm>*{};<16mm,20.5mm>*{}**@{-},
 <18mm,15mm>*{};<17mm,20.5mm>*{}**@{-},
 <18mm,15mm>*{};<18.5mm,19.6mm>*{}**@{-},
 <18mm,15mm>*{};<20mm,20.5mm>*{}**@{-},
 <18mm,15mm>*{};<18.6mm,20.3mm>*{.\hspace{-0.4mm}.\hspace{-0.4mm}.}**@{},
<22mm,15mm>*{};<25.5mm,17.7mm>*{\cdots}**@{},
 <32mm,15.2mm>*{};<30mm,20.5mm>*{}**@{-},
 <32mm,15.2mm>*{};<31mm,20.5mm>*{}**@{-},
 <32mm,15.2mm>*{};<32.5mm,19.6mm>*{}**@{-},
 <32mm,15mm>*{};<34mm,20.5mm>*{}**@{-},
 <32mm,15mm>*{};<32.3mm,20.3mm>*{.\hspace{-0.4mm}.\hspace{-0.4mm}.}**@{},
%
<0mm,0mm>*{};<18mm,22.6mm>*{^{I_{ i+1}''}}**@{},
<0mm,0mm>*{};<32.5mm,22.6mm>*{^{I_{ i+q}''}}**@{},
 <25mm,12mm>*{};<25mm,9mm>*{}**@{-},
 <25mm,12mm>*{};<26.2mm,9mm>*{}**@{-},
 <29mm,12mm>*{};<27.4mm,9.2mm>*{.\hspace{-0.1mm}.}**@{},
 <25mm,12mm>*{};<28.9mm,9mm>*{}**@{-},
 <29mm,12mm>*{};<27.2mm,7.3mm>*{_{J_2}}**@{},
<38mm,1.5mm>*{};<36mm,7mm>*{}**@{-},
<38mm,1.5mm>*{};<37mm,7mm>*{}**@{-},
<38mm,1.5mm>*{};<38.5mm,6mm>*{}**@{-},
<38mm,1.5mm>*{};<40mm,7mm>*{}**@{-},
<38mm,1.5mm>*{};<38.5mm,6.6mm>*{.\hspace{-0.4mm}.\hspace{-0.4mm}.}**@{},
<38mm,1.5mm>*{};<43mm,4mm>*{.\hspace{-0.0mm}.\hspace{-0.0mm}.}**@{},
<48mm,1.5mm>*{};<46mm,7mm>*{}**@{-},
<48mm,1.5mm>*{};<47mm,7mm>*{}**@{-},
<48mm,1.5mm>*{};<48.5mm,6mm>*{}**@{-},
<48mm,1.5mm>*{};<50mm,7mm>*{}**@{-},
<48mm,1.5mm>*{};<48.5mm,6.6mm>*{.\hspace{-0.4mm}.\hspace{-0.4mm}.}**@{},
<0mm,0mm>*{};<40.3mm,8.6mm>*{^{I_{ i+q+1}}}**@{},
<0mm,0mm>*{};<48.5mm,8.6mm>*{^{I_{ k}}}**@{},
\endxy}
$$
where $\var={i+l(n-i-q)+1}$.

\subsection{Important remark}\label{4: DefQ genus completion} The latter example
shows that the free prop $\DefQ^+$ makes sense as a {\em differential}\, prop only if it is considered as genus completed. This immediately raises a question: what is a representation of the completed prop $\DefQ^+$ in an arbitrary vector space $V$? A  morphism of dg props
$$
\rho: \DefQ^+ \lon \mathsf{End}_V
$$
is uniquely specified by its values on the generators,
$$
\rho\left(
\resizebox{17mm}{!}{\xy
 <0mm,0mm>*{\mbox{$\xy *=<20mm,3mm>\txt{}*\frm{-}\endxy$}};<0mm,0mm>*{}**@{},
  <-10mm,1.5mm>*{};<-12mm,7mm>*{}**@{-},
  <-10mm,1.5mm>*{};<-11mm,7mm>*{}**@{-},
  <-10mm,1.5mm>*{};<-9.5mm,6mm>*{}**@{-},
  <-10mm,1.5mm>*{};<-8mm,7mm>*{}**@{-},
 <-10mm,1.5mm>*{};<-9.5mm,6.6mm>*{.\hspace{-0.4mm}.\hspace{-0.4mm}.}**@{},
 <0mm,0mm>*{};<-6.5mm,3.6mm>*{.\hspace{-0.1mm}.\hspace{-0.1mm}.}**@{},
  <-3mm,1.5mm>*{};<-5mm,7mm>*{}**@{-},
  <-3mm,1.5mm>*{};<-4mm,7mm>*{}**@{-},
  <-3mm,1.5mm>*{};<-2.5mm,6mm>*{}**@{-},
  <-3mm,1.5mm>*{};<-1mm,7mm>*{}**@{-},
 <-3mm,1.5mm>*{};<-2.5mm,6.6mm>*{.\hspace{-0.4mm}.\hspace{-0.4mm}.}**@{},
  <2mm,1.5mm>*{};<0mm,7mm>*{}**@{-},
  <2mm,1.5mm>*{};<1mm,7mm>*{}**@{-},
  <2mm,1.5mm>*{};<2.5mm,6mm>*{}**@{-},
  <2mm,1.5mm>*{};<4mm,7mm>*{}**@{-},
 <2mm,1.5mm>*{};<2.5mm,6.6mm>*{.\hspace{-0.4mm}.\hspace{-0.4mm}.}**@{},
 <0mm,0mm>*{};<6mm,3.6mm>*{.\hspace{-0.1mm}.\hspace{-0.1mm}.}**@{},
<10mm,1.5mm>*{};<8mm,7mm>*{}**@{-},
  <10mm,1.5mm>*{};<9mm,7mm>*{}**@{-},
  <10mm,1.5mm>*{};<10.5mm,6mm>*{}**@{-},
  <10mm,1.5mm>*{};<12mm,7mm>*{}**@{-},
 <10mm,1.5mm>*{};<10.5mm,6.6mm>*{.\hspace{-0.4mm}.\hspace{-0.4mm}.}**@{},
%
<0mm,0mm>*{};<-9.5mm,8.2mm>*{^{I_{ 1}}}**@{},
<0mm,0mm>*{};<-3mm,8.2mm>*{^{I_{ i}}}**@{},
<0mm,0mm>*{};<2mm,8.2mm>*{^{I_{ i+1}}}**@{},
<0mm,0mm>*{};<10mm,8.2mm>*{^{I_{ m}}}**@{},
<-10mm,-1.5mm>*{};<-12mm,-7mm>*{}**@{-},
  <-10mm,-1.5mm>*{};<-11mm,-7mm>*{}**@{-},
  <-10mm,-1.5mm>*{};<-9.5mm,-6mm>*{}**@{-},
  <-10mm,-1.5mm>*{};<-8mm,-7mm>*{}**@{-},
 <-10mm,-1.5mm>*{};<-9.5mm,-6.6mm>*{.\hspace{-0.4mm}.\hspace{-0.4mm}.}**@{},
 <0mm,0mm>*{};<-6.5mm,-3.6mm>*{.\hspace{-0.1mm}.\hspace{-0.1mm}.}**@{},
  <-3mm,-1.5mm>*{};<-5mm,-7mm>*{}**@{-},
  <-3mm,-1.5mm>*{};<-4mm,-7mm>*{}**@{-},
  <-3mm,-1.5mm>*{};<-2.5mm,-6mm>*{}**@{-},
  <-3mm,-1.5mm>*{};<-1mm,-7mm>*{}**@{-},
 <-3mm,-1.5mm>*{};<-2.5mm,-6.6mm>*{.\hspace{-0.4mm}.\hspace{-0.4mm}.}**@{},
  <2mm,-1.5mm>*{};<0mm,-7mm>*{}**@{-},
  <2mm,-1.5mm>*{};<1mm,-7mm>*{}**@{-},
  <2mm,-1.5mm>*{};<2.5mm,-6mm>*{}**@{-},
  <2mm,-1.5mm>*{};<4mm,-7mm>*{}**@{-},
 <2mm,-1.5mm>*{};<2.5mm,-6.6mm>*{.\hspace{-0.4mm}.\hspace{-0.4mm}.}**@{},
 <0mm,0mm>*{};<6mm,-3.6mm>*{.\hspace{-0.1mm}.\hspace{-0.1mm}.}**@{},
<10mm,-1.5mm>*{};<8mm,-7mm>*{}**@{-},
  <10mm,-1.5mm>*{};<9mm,-7mm>*{}**@{-},
  <10mm,-1.5mm>*{};<10.5mm,-6mm>*{}**@{-},
  <10mm,-1.5mm>*{};<12mm,-7mm>*{}**@{-},
 <10mm,-1.5mm>*{};<10.5mm,-6.6mm>*{.\hspace{-0.4mm}.\hspace{-0.4mm}.}**@{},
%
<0mm,0mm>*{};<-9.5mm,-9.2mm>*{^{J_{ 1}}}**@{},
<0mm,0mm>*{};<-3mm,-9.2mm>*{^{J_{ ij}}}**@{},
<0mm,0mm>*{};<2mm,-9.2mm>*{^{J_{ j+1}}}**@{},
<0mm,0mm>*{};<10mm,-9.2mm>*{^{J_{ n}}}**@{},
\endxy}
\right)\subset \Hom(\ot^n \f_V, \ot^m \f_V) \subset \mathsf{End}_V
$$
which we assume from now on {\em to vanish for sufficiently large values
of the cardinalities $|I_i|$ and $|J_j|$}, $i\in [m]$, $j\in [n]$, $m,n\geq 1$. Under this assumption the map
$\rho$ is well-defined, and when we consider the deformation complex
$\Def(\DefQ^+\rar  \mathsf{End}_V)$ we always tacitly assume that we work within the class of the maps defined just above. Only under these assumptions on representations of $\DefQ^+$ we have an isomorphism of complexes,
$$
\Def(\DefQ^+\rar  \mathsf{End}_V)\simeq \Def(\sB_\infty^+ \rar  \mathsf{End}_{\f_V})
$$
Indeed, without the above assumption on representations $\rho$ we shall
get an element in  $\Def(\sB_\infty^+ \rar  \mathsf{End}_{\f_V})$ which is an
{\em infinite}\, sum of polydifferential operators of the form (\ref{polyoperator})
(with  order tending to infinity); such an infinite sum can not be an element of
 $\Def(\sB_\infty^+ \rar  \mathsf{End}_{\f_V})$ (and it has no sense to apply $L_\infty$
 operation $\mu_\bullet$ to such an infinite sum).

\subsection{ $\DefQ^+$ versus $\Lieb_\infty^+$}\label{cohomology-DefQ}
It is easy to check that a derivation, $\sd_1$, of $\DefQ^+$, given on generators by
\Beq\label{sd_1}
\sd_1\left(
\resizebox{17mm}{!}{\xy
 <0mm,0mm>*{\mbox{$\xy *=<20mm,3mm>\txt{}*\frm{-}\endxy$}};<0mm,0mm>*{}**@{},
  <-10mm,1.5mm>*{};<-12mm,7mm>*{}**@{-},
  <-10mm,1.5mm>*{};<-11mm,7mm>*{}**@{-},
  <-10mm,1.5mm>*{};<-9.5mm,6mm>*{}**@{-},
  <-10mm,1.5mm>*{};<-8mm,7mm>*{}**@{-},
 <-10mm,1.5mm>*{};<-9.5mm,6.6mm>*{.\hspace{-0.4mm}.\hspace{-0.4mm}.}**@{},
 <0mm,0mm>*{};<-6.5mm,3.6mm>*{.\hspace{-0.1mm}.\hspace{-0.1mm}.}**@{},
  <-3mm,1.5mm>*{};<-5mm,7mm>*{}**@{-},
  <-3mm,1.5mm>*{};<-4mm,7mm>*{}**@{-},
  <-3mm,1.5mm>*{};<-2.5mm,6mm>*{}**@{-},
  <-3mm,1.5mm>*{};<-1mm,7mm>*{}**@{-},
 <-3mm,1.5mm>*{};<-2.5mm,6.6mm>*{.\hspace{-0.4mm}.\hspace{-0.4mm}.}**@{},
  <2mm,1.5mm>*{};<0mm,7mm>*{}**@{-},
  <2mm,1.5mm>*{};<1mm,7mm>*{}**@{-},
  <2mm,1.5mm>*{};<2.5mm,6mm>*{}**@{-},
  <2mm,1.5mm>*{};<4mm,7mm>*{}**@{-},
 <2mm,1.5mm>*{};<2.5mm,6.6mm>*{.\hspace{-0.4mm}.\hspace{-0.4mm}.}**@{},
 <0mm,0mm>*{};<6mm,3.6mm>*{.\hspace{-0.1mm}.\hspace{-0.1mm}.}**@{},
<10mm,1.5mm>*{};<8mm,7mm>*{}**@{-},
  <10mm,1.5mm>*{};<9mm,7mm>*{}**@{-},
  <10mm,1.5mm>*{};<10.5mm,6mm>*{}**@{-},
  <10mm,1.5mm>*{};<12mm,7mm>*{}**@{-},
 <10mm,1.5mm>*{};<10.5mm,6.6mm>*{.\hspace{-0.4mm}.\hspace{-0.4mm}.}**@{},
%
<0mm,0mm>*{};<-9.5mm,8.2mm>*{^{I_{ 1}}}**@{},
<0mm,0mm>*{};<-3mm,8.2mm>*{^{I_{ i}}}**@{},
<0mm,0mm>*{};<2mm,8.2mm>*{^{I_{ i+1}}}**@{},
<0mm,0mm>*{};<10mm,8.2mm>*{^{I_{ m}}}**@{},
<-10mm,-1.5mm>*{};<-12mm,-7mm>*{}**@{-},
  <-10mm,-1.5mm>*{};<-11mm,-7mm>*{}**@{-},
  <-10mm,-1.5mm>*{};<-9.5mm,-6mm>*{}**@{-},
  <-10mm,-1.5mm>*{};<-8mm,-7mm>*{}**@{-},
 <-10mm,-1.5mm>*{};<-9.5mm,-6.6mm>*{.\hspace{-0.4mm}.\hspace{-0.4mm}.}**@{},
 <0mm,0mm>*{};<-6.5mm,-3.6mm>*{.\hspace{-0.1mm}.\hspace{-0.1mm}.}**@{},
  <-3mm,-1.5mm>*{};<-5mm,-7mm>*{}**@{-},
  <-3mm,-1.5mm>*{};<-4mm,-7mm>*{}**@{-},
  <-3mm,-1.5mm>*{};<-2.5mm,-6mm>*{}**@{-},
  <-3mm,-1.5mm>*{};<-1mm,-7mm>*{}**@{-},
 <-3mm,-1.5mm>*{};<-2.5mm,-6.6mm>*{.\hspace{-0.4mm}.\hspace{-0.4mm}.}**@{},
  <2mm,-1.5mm>*{};<0mm,-7mm>*{}**@{-},
  <2mm,-1.5mm>*{};<1mm,-7mm>*{}**@{-},
  <2mm,-1.5mm>*{};<2.5mm,-6mm>*{}**@{-},
  <2mm,-1.5mm>*{};<4mm,-7mm>*{}**@{-},
 <2mm,-1.5mm>*{};<2.5mm,-6.6mm>*{.\hspace{-0.4mm}.\hspace{-0.4mm}.}**@{},
 <0mm,0mm>*{};<6mm,-3.6mm>*{.\hspace{-0.1mm}.\hspace{-0.1mm}.}**@{},
<10mm,-1.5mm>*{};<8mm,-7mm>*{}**@{-},
  <10mm,-1.5mm>*{};<9mm,-7mm>*{}**@{-},
  <10mm,-1.5mm>*{};<10.5mm,-6mm>*{}**@{-},
  <10mm,-1.5mm>*{};<12mm,-7mm>*{}**@{-},
 <10mm,-1.5mm>*{};<10.5mm,-6.6mm>*{.\hspace{-0.4mm}.\hspace{-0.4mm}.}**@{},
%
<0mm,0mm>*{};<-9.5mm,-9.2mm>*{^{J_{ 1}}}**@{},
<0mm,0mm>*{};<-3mm,-9.2mm>*{^{J_{ ij}}}**@{},
<0mm,0mm>*{};<2mm,-9.2mm>*{^{J_{ j+1}}}**@{},
<0mm,0mm>*{};<10mm,-9.2mm>*{^{J_{ n}}}**@{},
\endxy}
\right)
=\sum_{i=1}^{m-1}(-1)^i
\resizebox{17mm}{!}{\xy
 <0mm,0mm>*{\mbox{$\xy *=<20mm,3mm>\txt{}*\frm{-}\endxy$}};<0mm,0mm>*{}**@{},
  <-10mm,1.5mm>*{};<-12mm,7mm>*{}**@{-},
  <-10mm,1.5mm>*{};<-11mm,7mm>*{}**@{-},
  <-10mm,1.5mm>*{};<-9.5mm,6mm>*{}**@{-},
  <-10mm,1.5mm>*{};<-8mm,7mm>*{}**@{-},
 <-10mm,1.5mm>*{};<-9.5mm,6.6mm>*{.\hspace{-0.4mm}.\hspace{-0.4mm}.}**@{},
 <0mm,0mm>*{};<-6.5mm,3.6mm>*{.\hspace{-0.1mm}.\hspace{-0.1mm}.}**@{},
  <0mm,1.5mm>*{};<-5mm,7mm>*{}**@{-},
  <0mm,1.5mm>*{};<-4mm,7mm>*{}**@{-},
  <0mm,1.5mm>*{};<-2.5mm,6mm>*{}**@{-},
  <0mm,1.5mm>*{};<-1.3mm,7mm>*{}**@{-},
 <0mm,1.5mm>*{};<-2.5mm,6.6mm>*{.\hspace{-0.4mm}.\hspace{-0.4mm}.}**@{},
  <0mm,1.5mm>*{};<1.3mm,7mm>*{}**@{-},
  <0mm,1.5mm>*{};<2.5mm,6mm>*{}**@{-},
  <0mm,1.5mm>*{};<4mm,7mm>*{}**@{-},
 <0mm,1.5mm>*{};<2.5mm,6.6mm>*{.\hspace{-0.4mm}.\hspace{-0.4mm}.}**@{},
 <0mm,0mm>*{};<6mm,3.6mm>*{.\hspace{-0.1mm}.\hspace{-0.1mm}.}**@{},
<10mm,1.5mm>*{};<8mm,7mm>*{}**@{-},
  <10mm,1.5mm>*{};<9mm,7mm>*{}**@{-},
  <10mm,1.5mm>*{};<10.5mm,6mm>*{}**@{-},
  <10mm,1.5mm>*{};<12mm,7mm>*{}**@{-},
 <10mm,1.5mm>*{};<10.5mm,6.6mm>*{.\hspace{-0.4mm}.\hspace{-0.4mm}.}**@{},
%
<0mm,0mm>*{};<-9.5mm,8.2mm>*{^{I_{ 1}}}**@{},
<0mm,0mm>*{};<-3mm,8.2mm>*{^{I_{ i}}}**@{},
<0mm,0mm>*{};<2mm,8.2mm>*{^{\sqcup I_{ i+1}}}**@{},
<0mm,0mm>*{};<10mm,8.2mm>*{^{I_{ m}}}**@{},
<-10mm,-1.5mm>*{};<-12mm,-7mm>*{}**@{-},
  <-10mm,-1.5mm>*{};<-11mm,-7mm>*{}**@{-},
  <-10mm,-1.5mm>*{};<-9.5mm,-6mm>*{}**@{-},
  <-10mm,-1.5mm>*{};<-8mm,-7mm>*{}**@{-},
 <-10mm,-1.5mm>*{};<-9.5mm,-6.6mm>*{.\hspace{-0.4mm}.\hspace{-0.4mm}.}**@{},
 <0mm,0mm>*{};<-6.5mm,-3.6mm>*{.\hspace{-0.1mm}.\hspace{-0.1mm}.}**@{},
  <-3mm,-1.5mm>*{};<-5mm,-7mm>*{}**@{-},
  <-3mm,-1.5mm>*{};<-4mm,-7mm>*{}**@{-},
  <-3mm,-1.5mm>*{};<-2.5mm,-6mm>*{}**@{-},
  <-3mm,-1.5mm>*{};<-1mm,-7mm>*{}**@{-},
 <-3mm,-1.5mm>*{};<-2.5mm,-6.6mm>*{.\hspace{-0.4mm}.\hspace{-0.4mm}.}**@{},
  <2mm,-1.5mm>*{};<0mm,-7mm>*{}**@{-},
  <2mm,-1.5mm>*{};<1mm,-7mm>*{}**@{-},
  <2mm,-1.5mm>*{};<2.5mm,-6mm>*{}**@{-},
  <2mm,-1.5mm>*{};<4mm,-7mm>*{}**@{-},
 <2mm,-1.5mm>*{};<2.5mm,-6.6mm>*{.\hspace{-0.4mm}.\hspace{-0.4mm}.}**@{},
 <0mm,0mm>*{};<6mm,-3.6mm>*{.\hspace{-0.1mm}.\hspace{-0.1mm}.}**@{},
<10mm,-1.5mm>*{};<8mm,-7mm>*{}**@{-},
  <10mm,-1.5mm>*{};<9mm,-7mm>*{}**@{-},
  <10mm,-1.5mm>*{};<10.5mm,-6mm>*{}**@{-},
  <10mm,-1.5mm>*{};<12mm,-7mm>*{}**@{-},
 <10mm,-1.5mm>*{};<10.5mm,-6.6mm>*{.\hspace{-0.4mm}.\hspace{-0.4mm}.}**@{},
%
<0mm,0mm>*{};<-9.5mm,-9.2mm>*{^{J_{ 1}}}**@{},
<0mm,0mm>*{};<-3mm,-9.2mm>*{^{J_{ j}}}**@{},
<0mm,0mm>*{};<2mm,-9.2mm>*{^{J_{ j+1}}}**@{},
<0mm,0mm>*{};<10mm,-9.2mm>*{^{J_{ n}}}**@{},
\endxy}
+
\sum_{j=1}^{n-1}(-1)^j
\resizebox{17mm}{!}{\xy
 <0mm,0mm>*{\mbox{$\xy *=<20mm,3mm>\txt{}*\frm{-}\endxy$}};<0mm,0mm>*{}**@{},
  <-10mm,1.5mm>*{};<-12mm,7mm>*{}**@{-},
  <-10mm,1.5mm>*{};<-11mm,7mm>*{}**@{-},
  <-10mm,1.5mm>*{};<-9.5mm,6mm>*{}**@{-},
  <-10mm,1.5mm>*{};<-8mm,7mm>*{}**@{-},
 <-10mm,1.5mm>*{};<-9.5mm,6.6mm>*{.\hspace{-0.4mm}.\hspace{-0.4mm}.}**@{},
 <0mm,0mm>*{};<-6.5mm,3.6mm>*{.\hspace{-0.1mm}.\hspace{-0.1mm}.}**@{},
  <-3mm,1.5mm>*{};<-5mm,7mm>*{}**@{-},
  <-3mm,1.5mm>*{};<-4mm,7mm>*{}**@{-},
  <-3mm,1.5mm>*{};<-2.5mm,6mm>*{}**@{-},
  <-3mm,1.5mm>*{};<-1mm,7mm>*{}**@{-},
 <-3mm,1.5mm>*{};<-2.5mm,6.6mm>*{.\hspace{-0.4mm}.\hspace{-0.4mm}.}**@{},
  <2mm,1.5mm>*{};<0mm,7mm>*{}**@{-},
  <2mm,1.5mm>*{};<1mm,7mm>*{}**@{-},
  <2mm,1.5mm>*{};<2.5mm,6mm>*{}**@{-},
  <2mm,1.5mm>*{};<4mm,7mm>*{}**@{-},
 <2mm,1.5mm>*{};<2.5mm,6.6mm>*{.\hspace{-0.4mm}.\hspace{-0.4mm}.}**@{},
 <0mm,0mm>*{};<6mm,3.6mm>*{.\hspace{-0.1mm}.\hspace{-0.1mm}.}**@{},
<10mm,1.5mm>*{};<8mm,7mm>*{}**@{-},
  <10mm,1.5mm>*{};<9mm,7mm>*{}**@{-},
  <10mm,1.5mm>*{};<10.5mm,6mm>*{}**@{-},
  <10mm,1.5mm>*{};<12mm,7mm>*{}**@{-},
 <10mm,1.5mm>*{};<10.5mm,6.6mm>*{.\hspace{-0.4mm}.\hspace{-0.4mm}.}**@{},
%
<0mm,0mm>*{};<-9.5mm,8.2mm>*{^{I_{ 1}}}**@{},
<0mm,0mm>*{};<-3mm,8.2mm>*{^{I_{ i}}}**@{},
<0mm,0mm>*{};<2mm,8.2mm>*{^{I_{ i+1}}}**@{},
<0mm,0mm>*{};<10mm,8.2mm>*{^{I_{ m}}}**@{},
<-10mm,-1.5mm>*{};<-12mm,-7mm>*{}**@{-},
  <-10mm,-1.5mm>*{};<-11mm,-7mm>*{}**@{-},
  <-10mm,-1.5mm>*{};<-9.5mm,-6mm>*{}**@{-},
  <-10mm,-1.5mm>*{};<-8mm,-7mm>*{}**@{-},
 <-10mm,-1.5mm>*{};<-9.5mm,-6.6mm>*{.\hspace{-0.4mm}.\hspace{-0.4mm}.}**@{},
 <0mm,0mm>*{};<-6.5mm,-3.6mm>*{.\hspace{-0.1mm}.\hspace{-0.1mm}.}**@{},
  <0mm,-1.5mm>*{};<-5mm,-7mm>*{}**@{-},
  <0mm,-1.5mm>*{};<-4mm,-7mm>*{}**@{-},
  <0mm,-1.5mm>*{};<-2.5mm,-6mm>*{}**@{-},
  <0mm,-1.5mm>*{};<-1.3mm,-7mm>*{}**@{-},
 <0mm,-1.5mm>*{};<-2.5mm,-6.6mm>*{.\hspace{-0.4mm}.\hspace{-0.4mm}.}**@{},
  <0mm,-1.5mm>*{};<1.3mm,-7mm>*{}**@{-},
  <0mm,-1.5mm>*{};<2.5mm,-6mm>*{}**@{-},
  <0mm,-1.5mm>*{};<4mm,-7mm>*{}**@{-},
 <0mm,-1.5mm>*{};<2.5mm,-6.6mm>*{.\hspace{-0.4mm}.\hspace{-0.4mm}.}**@{},
 <0mm,0mm>*{};<6mm,-3.6mm>*{.\hspace{-0.1mm}.\hspace{-0.1mm}.}**@{},
<10mm,-1.5mm>*{};<8mm,-7mm>*{}**@{-},
  <10mm,-1.5mm>*{};<9mm,-7mm>*{}**@{-},
  <10mm,-1.5mm>*{};<10.5mm,-6mm>*{}**@{-},
  <10mm,-1.5mm>*{};<12mm,-7mm>*{}**@{-},
 <10mm,-1.5mm>*{};<10.5mm,-6.6mm>*{.\hspace{-0.4mm}.\hspace{-0.4mm}.}**@{},
%
<0mm,0mm>*{};<-9.5mm,-9.2mm>*{^{J_{ 1}}}**@{},
<0mm,0mm>*{};<-3mm,-9.2mm>*{^{J_{ j}}}**@{},
<0mm,0mm>*{};<2.2mm,-9.2mm>*{^{\sqcup J_{ j+1}}}**@{},
<0mm,0mm>*{};<10mm,-9.2mm>*{^{J_{ n}}}**@{},
\endxy}
\Eeq
is a differential. By (\ref{graph-GS}), the differential $\sd_1$ is a part part of the
full differential
$\sd$. As a first approximation to the formality interrelation between $(\DefQ^+, \sd)$ and $(\Lieb_\infty^+, \delta^+)$,
we notice the following result which was proven in Proposition 3.5.1 of \cite{Me3}.

\begin{theorem}\label{sd-cohomology}
The cohomology prop $H^\bullet(\DefQ^+, \sd_1)$ is isomorphic to the (non-differential) prop ${\Lieb}_\infty^+$.
\end{theorem}

\subsection{A quotient of $\DefQ^+$}\label{defq} The formula for the differential in \S~\ref{Hoch} implies
$$
\sd\, \Ba{c}
\resizebox{9mm}{!}{\xy
 <-10mm,0mm>*{\mbox{$\xy *=<3mm,3mm>\txt{}*\frm{-}\endxy$}};<0mm,0mm>*{}**@{},
  <-10mm,1.5mm>*{};<-12mm,7mm>*{}**@{-},
  <-10mm,1.5mm>*{};<-11mm,7mm>*{}**@{-},
  <-10mm,1.5mm>*{};<-9.5mm,6mm>*{}**@{-},
  <-10mm,1.5mm>*{};<-8mm,7mm>*{}**@{-},
 <-10mm,-1.5mm>*{};<-9.5mm,6.6mm>*{.\hspace{-0.4mm}.\hspace{-0.4mm}.}**@{},
<-10mm,-1.5mm>*{};<-12mm,-7mm>*{}**@{-},
  <-10mm,-1.5mm>*{};<-11mm,-7mm>*{}**@{-},
  <-10mm,-1.5mm>*{};<-9.5mm,-6mm>*{}**@{-},
  <-10mm,-1.5mm>*{};<-8mm,-7mm>*{}**@{-},
 <-10mm,-1.5mm>*{};<-9.5mm,-6.6mm>*{.\hspace{-0.4mm}.\hspace{-0.4mm}.}**@{},
<0mm,0mm>*{};<-9.5mm,-9.2mm>*{^{J}}**@{},
<0mm,0mm>*{};<-9.5mm,8.2mm>*{^{I}}**@{},
\endxy}\Ea
=\sum_{I=I'\sqcup I''\atop J=J'\sqcup J''}\sum_{s\geq 1}\frac{1}{s!}\
\resizebox{11mm}{!}{\xy
 <-10mm,-3mm>*{\mbox{$\xy *=<3mm,3mm>\txt{}*\frm{-}\endxy$}};<0mm,0mm>*{}**@{},
<-10mm,-4.5mm>*{};<-12mm,-10mm>*{}**@{-},
  <-10mm,-4.5mm>*{};<-11mm,-10mm>*{}**@{-},
  <-10mm,-4.5mm>*{};<-9.5mm,-9mm>*{}**@{-},
  <-10mm,-4.5mm>*{};<-8mm,-10mm>*{}**@{-},
 <-10mm,-4.5mm>*{};<-9.5mm,-9.6mm>*{.\hspace{-0.4mm}.\hspace{-0.4mm}.}**@{},
<-5mm,6mm>*{\mbox{$\xy *=<3mm,3mm>\txt{}*\frm{-}\endxy$}};<0mm,0mm>*{}**@{},
<-5mm,7.5mm>*{};<-7mm,13mm>*{}**@{-},
  <-5mm,7.5mm>*{};<-6mm,13mm>*{}**@{-},
  <-5mm,7.5mm>*{};<-4.5mm,12mm>*{}**@{-},
  <-5mm,7.5mm>*{};<-3mm,13mm>*{}**@{-},
 <-5mm,7.5mm>*{};<-4.5mm,12.6mm>*{.\hspace{-0.4mm}.\hspace{-0.4mm}.}**@{},
<-10mm,-1.5mm>*{};<-12.5mm,4mm>*{}**@{-},
  <-10mm,-1.5mm>*{};<-9.5mm,4mm>*{}**@{-},
 <-10mm,-1.5mm>*{};<-10.8mm,3.6mm>*{.\hspace{-0.4mm}.\hspace{-0.4mm}.}**@{},
<0mm,0mm>*{};<-9.5mm,-12.8mm>*{^{J'}}**@{},
<0mm,0mm>*{};<-2.2mm,-2.5mm>*{^{J''}}**@{},
<0mm,0mm>*{};<-10.5mm,6.2mm>*{^{I'}}**@{},
<0mm,0mm>*{};<-4.5mm,15.2mm>*{^{I''}}**@{},
<-10mm,-1.5mm>*{};<-8.3mm,2mm>*{}**@{-},
<-8.3mm,2mm>*{};<-5mm,4.5mm>*{}**@{-},
<-10mm,-1.5mm>*{};<-5.7mm,1mm>*{}**@{-},
<-5.7mm,1mm>*{};<-5mm,4.5mm>*{}**@{-},
 <0mm,0mm>*{};<-7mm,1.2mm>*{.\hspace{-0.4mm}.}**@{},
<0mm,0mm>*{};<-6.2mm,1.6mm>*{^s}**@{},
<-5mm,4.5mm>*{};<-4mm,0.5mm>*{}**@{-},
<-5mm,4.5mm>*{};<-1mm,0.5mm>*{}**@{-},
 <0mm,0mm>*{};<-2.7mm,1mm>*{.\hspace{-0.4mm}.}**@{},
\endxy}
$$
which in turn implies that the ideal $\mathsf I$ in $\DefQ^+$ generated by corollas \
$
\Ba{c}\resizebox{7mm}{!}{\xy
 <-10mm,0mm>*{\mbox{$\xy *=<3mm,3mm>\txt{}*\frm{-}\endxy$}};<0mm,0mm>*{}**@{},
  <-10mm,1.5mm>*{};<-12mm,7mm>*{}**@{-},
  <-10mm,1.5mm>*{};<-11mm,7mm>*{}**@{-},
  <-10mm,1.5mm>*{};<-9.5mm,6mm>*{}**@{-},
  <-10mm,1.5mm>*{};<-8mm,7mm>*{}**@{-},
 <-10mm,-1.5mm>*{};<-9.5mm,6.6mm>*{.\hspace{-0.4mm}.\hspace{-0.4mm}.}**@{},
<-10mm,-1.5mm>*{};<-12mm,-7mm>*{}**@{-},
  <-10mm,-1.5mm>*{};<-11mm,-7mm>*{}**@{-},
  <-10mm,-1.5mm>*{};<-9.5mm,-6mm>*{}**@{-},
  <-10mm,-1.5mm>*{};<-8mm,-7mm>*{}**@{-},
 <-10mm,-1.5mm>*{};<-9.5mm,-6.6mm>*{.\hspace{-0.4mm}.\hspace{-0.4mm}.}**@{},
<0mm,0mm>*{};<-9.5mm,-9.2mm>*{^{J}}**@{},
<0mm,0mm>*{};<-9.5mm,8.2mm>*{^{I}}**@{},
\endxy}\Ea
$
is differential. Hence the quotient,
$
\DefQ:= {\DefQ^+}/{\mathsf I},
$
is a {\em differential}\,  prop.
The induced differential we denote by $\sd$. It is a free prop generated by corollas
$
\Ba{c}\resizebox{15mm}{!}{\xy
 <0mm,0mm>*{\mbox{$\xy *=<20mm,3mm>\txt{}*\frm{-}\endxy$}};<0mm,0mm>*{}**@{},
  <-10mm,1.5mm>*{};<-12mm,7mm>*{}**@{-},
  <-10mm,1.5mm>*{};<-11mm,7mm>*{}**@{-},
  <-10mm,1.5mm>*{};<-9.5mm,6mm>*{}**@{-},
  <-10mm,1.5mm>*{};<-8mm,7mm>*{}**@{-},
 <-10mm,1.5mm>*{};<-9.5mm,6.6mm>*{.\hspace{-0.4mm}.\hspace{-0.4mm}.}**@{},
 <0mm,0mm>*{};<-6.5mm,3.6mm>*{.\hspace{-0.1mm}.\hspace{-0.1mm}.}**@{},
  <-3mm,1.5mm>*{};<-5mm,7mm>*{}**@{-},
  <-3mm,1.5mm>*{};<-4mm,7mm>*{}**@{-},
  <-3mm,1.5mm>*{};<-2.5mm,6mm>*{}**@{-},
  <-3mm,1.5mm>*{};<-1mm,7mm>*{}**@{-},
 <-3mm,1.5mm>*{};<-2.5mm,6.6mm>*{.\hspace{-0.4mm}.\hspace{-0.4mm}.}**@{},
  <2mm,1.5mm>*{};<0mm,7mm>*{}**@{-},
  <2mm,1.5mm>*{};<1mm,7mm>*{}**@{-},
  <2mm,1.5mm>*{};<2.5mm,6mm>*{}**@{-},
  <2mm,1.5mm>*{};<4mm,7mm>*{}**@{-},
 <2mm,1.5mm>*{};<2.5mm,6.6mm>*{.\hspace{-0.4mm}.\hspace{-0.4mm}.}**@{},
 <0mm,0mm>*{};<6mm,3.6mm>*{.\hspace{-0.1mm}.\hspace{-0.1mm}.}**@{},
<10mm,1.5mm>*{};<8mm,7mm>*{}**@{-},
  <10mm,1.5mm>*{};<9mm,7mm>*{}**@{-},
  <10mm,1.5mm>*{};<10.5mm,6mm>*{}**@{-},
  <10mm,1.5mm>*{};<12mm,7mm>*{}**@{-},
 <10mm,1.5mm>*{};<10.5mm,6.6mm>*{.\hspace{-0.4mm}.\hspace{-0.4mm}.}**@{},
%
<0mm,0mm>*{};<-9.5mm,8.2mm>*{^{I_{ 1}}}**@{},
<0mm,0mm>*{};<-3mm,8.2mm>*{^{I_{ i}}}**@{},
<0mm,0mm>*{};<2mm,8.2mm>*{^{I_{ i+1}}}**@{},
<0mm,0mm>*{};<10mm,8.2mm>*{^{I_{ m}}}**@{},
<-10mm,-1.5mm>*{};<-12mm,-7mm>*{}**@{-},
  <-10mm,-1.5mm>*{};<-11mm,-7mm>*{}**@{-},
  <-10mm,-1.5mm>*{};<-9.5mm,-6mm>*{}**@{-},
  <-10mm,-1.5mm>*{};<-8mm,-7mm>*{}**@{-},
 <-10mm,-1.5mm>*{};<-9.5mm,-6.6mm>*{.\hspace{-0.4mm}.\hspace{-0.4mm}.}**@{},
 <0mm,0mm>*{};<-6.5mm,-3.6mm>*{.\hspace{-0.1mm}.\hspace{-0.1mm}.}**@{},
  <-3mm,-1.5mm>*{};<-5mm,-7mm>*{}**@{-},
  <-3mm,-1.5mm>*{};<-4mm,-7mm>*{}**@{-},
  <-3mm,-1.5mm>*{};<-2.5mm,-6mm>*{}**@{-},
  <-3mm,-1.5mm>*{};<-1mm,-7mm>*{}**@{-},
 <-3mm,-1.5mm>*{};<-2.5mm,-6.6mm>*{.\hspace{-0.4mm}.\hspace{-0.4mm}.}**@{},
  <2mm,-1.5mm>*{};<0mm,-7mm>*{}**@{-},
  <2mm,-1.5mm>*{};<1mm,-7mm>*{}**@{-},
  <2mm,-1.5mm>*{};<2.5mm,-6mm>*{}**@{-},
  <2mm,-1.5mm>*{};<4mm,-7mm>*{}**@{-},
 <2mm,-1.5mm>*{};<2.5mm,-6.6mm>*{.\hspace{-0.4mm}.\hspace{-0.4mm}.}**@{},
 <0mm,0mm>*{};<6mm,-3.6mm>*{.\hspace{-0.1mm}.\hspace{-0.1mm}.}**@{},
<10mm,-1.5mm>*{};<8mm,-7mm>*{}**@{-},
  <10mm,-1.5mm>*{};<9mm,-7mm>*{}**@{-},
  <10mm,-1.5mm>*{};<10.5mm,-6mm>*{}**@{-},
  <10mm,-1.5mm>*{};<12mm,-7mm>*{}**@{-},
 <10mm,-1.5mm>*{};<10.5mm,-6.6mm>*{.\hspace{-0.4mm}.\hspace{-0.4mm}.}**@{},
%
<0mm,0mm>*{};<-9.5mm,-9.2mm>*{^{J_{ 1}}}**@{},
<0mm,0mm>*{};<-3mm,-9.2mm>*{^{J_{j}}}**@{},
<0mm,0mm>*{};<2mm,-9.2mm>*{^{J_{ j+1}}}**@{},
<0mm,0mm>*{};<10mm,-9.2mm>*{^{J_{ n}}}**@{},
\endxy}\Ea
$
with $m+n \geq 3$. It is clear that there is a one-to-one correspondence between
normalized polydifferential representations of the dg prop $\sB_\infty$ in $\f_V$ and ordinary (see \S \ref{4: DefQ genus completion}) representations of $\DefQ$ in a graded vector space $V$. The prop $\DefQ$ is non-positively graded
so that any morphism
 from $\DefQ^+$ into a non-positively graded prop $\sP$ factors through the natural projection to $\DefQ$, $\DefQ^+\rar \DefQ \rar \sP$.

\bip

\section{Formality morphism of $L_\infty$-algebras via morphism of props}

\subsection{\mbox{\bf On representations of completed props and formal parameters}}
\label{planck}
 The Etingof-Kazhdan
quantization morphism $\cE\cK$ (which is discussed in the next section, see (\ref{5: EK morphism}) below)
is a  morphism from the dg prop $\DefQ$ into a {\em completed}\, prop of Lie bialgebras. The value of $\cE\cK$ on a generator
of $\DefQ$  is a well-defined --- probably, infinite!  --- linear combination of graphs from $\widehat{\Lieb}$, but then one has to be careful on how to define a representation
of the completed prop $\widehat{\Lieb}$ in a vector space $V$ as such an infinite linear combinations of graphs is mapped in general into a divergent  infinite sum of linear operators. This the point where in the prop approach to deformation quantization the formal ``Planck constant" has to enter the story.



\sip

A bialgebra structure on $V$ is a morphism
of props,
$$
\rho: \Lieb \lon \End_V,
$$
with
$$
\rho\left(
\begin{xy}
 <0mm,-0.55mm>*{};<0mm,-2.5mm>*{}**@{-},
 <0.5mm,0.5mm>*{};<2.2mm,2.2mm>*{}**@{-},
 <-0.48mm,0.48mm>*{};<-2.2mm,2.2mm>*{}**@{-},
 <0mm,0mm>*{\bullet};<0mm,0mm>*{}**@{},
 <0mm,-0.55mm>*{};<0mm,-3.8mm>*{_1}**@{},
 <0.5mm,0.5mm>*{};<2.7mm,2.8mm>*{^1}**@{},
 <-0.48mm,0.48mm>*{};<-2.7mm,2.8mm>*{^2}**@{},
 \end{xy}
\right) =:\vartriangle\in \Hom(V, V\ot V), \ \ \
\rho\left(
\begin{xy}
 <0mm,0.66mm>*{};<0mm,3mm>*{}**@{-},
 <0.39mm,-0.39mm>*{};<2.2mm,-2.2mm>*{}**@{-},
 <-0.35mm,-0.35mm>*{};<-2.2mm,-2.2mm>*{}**@{-},
 <0mm,0mm>*{\bullet};<0mm,0mm>*{}**@{},
   <0mm,0.66mm>*{};<0mm,3.4mm>*{^1}**@{},
   <0.39mm,-0.39mm>*{};<2.9mm,-4mm>*{^2}**@{},
   <-0.35mm,-0.35mm>*{};<-2.8mm,-4mm>*{^1}**@{},
\end{xy}
\right)=:[\, ,\, ]\in \Hom(V\ot V, V)
$$
being co-Lie and, respectively, Lie brackets in $V$. However, the morphism $\rho$ can {\em not}\, be automatically
extended to a representation,
$$
\hat{\rho}: \widehat{\Lieb} \lon  \End\langle V\rangle,
$$
of the completed prop as, for example, a well-defined infinite sum of graphs in $\widehat{\Lieb}$,
$$
\sum_{n=1}^\infty\ \ \left.
\begin{xy}
 <0mm,-9mm>*{};<0mm,-11.5mm>*{}**@{-},
 <0mm,-9mm>*{};<2.2mm,-6.8mm>*{}**@{-},
 <0mm,-9mm>*{};<-2.2mm,-6.8mm>*{}**@{-},
 <0mm,-9mm>*{\bullet};
<2.2mm,-6.8mm>*{};<0mm,-4.6mm>*{}**@{-},
<-2.2mm,-6.8mm>*{};<0mm,-4.6mm>*{}**@{-},
<0mm,-4.6mm>*{\bullet};
<0mm,-4.6mm>*{};<0mm,-2.1mm>*{}**@{-},
<0mm,9mm>*{};<0mm,11.5mm>*{}**@{-},
 <0mm,9mm>*{};<2.2mm,6.8mm>*{}**@{-},
 <0mm,9mm>*{};<-2.2mm,6.8mm>*{}**@{-},
 <0mm,9mm>*{\bullet};
<2.2mm,6.8mm>*{};<0mm,4.6mm>*{}**@{-},
<-2.2mm,6.8mm>*{};<0mm,4.6mm>*{}**@{-},
<0mm,4.6mm>*{\bullet};
<0mm,4.6mm>*{};<0mm,2.1mm>*{}**@{-},
 <0mm,0mm>*{\cdot};
<0mm,1mm>*{\cdot};
<0mm,-1mm>*{\cdot};
 \end{xy}
\right\}\mbox{$n$ times}
$$
gets mapped into an infinite sum of elements of $\Hom(V,V)$ which is divergent in general.

Thus
we interpret instead a bialgebra structure in $V$ as a morphism of props,
$$
{\rho}: {\Lieb} \lon  \End_V[[\hbar]],
$$
by setting
$$
\rho\left(
\begin{xy}
 <0mm,-0.55mm>*{};<0mm,-2.5mm>*{}**@{-},
 <0.5mm,0.5mm>*{};<2.2mm,2.2mm>*{}**@{-},
 <-0.48mm,0.48mm>*{};<-2.2mm,2.2mm>*{}**@{-},
 <0mm,0mm>*{\bullet};<0mm,0mm>*{}**@{},
 <0mm,-0.55mm>*{};<0mm,-3.8mm>*{_1}**@{},
 <0.5mm,0.5mm>*{};<2.7mm,2.8mm>*{^1}**@{},
 <-0.48mm,0.48mm>*{};<-2.7mm,2.8mm>*{^2}**@{},
 \end{xy}
\right) =:\hbar\vartriangle\in \Hom(V, V\ot V)[[\hbar]], \ \ \
\rho\left(
\begin{xy}
 <0mm,0.66mm>*{};<0mm,3mm>*{}**@{-},
 <0.39mm,-0.39mm>*{};<2.2mm,-2.2mm>*{}**@{-},
 <-0.35mm,-0.35mm>*{};<-2.2mm,-2.2mm>*{}**@{-},
 <0mm,0mm>*{\bullet};<0mm,0mm>*{}**@{},
   <0mm,0.66mm>*{};<0mm,3.4mm>*{^1}**@{},
   <0.39mm,-0.39mm>*{};<2.9mm,-4mm>*{^2}**@{},
   <-0.35mm,-0.35mm>*{};<-2.8mm,-4mm>*{^1}**@{},
\end{xy}
\right)=:\hbar[\, ,\, ]\in \Hom(V\ot V, V)[[\hbar]].
$$
Here $\hbar$ is a formal parameter and, for a vector space $W$ the symbol $W[[\hbar]]$ stands for the vector space of formal power series in $\hbar$ with coefficients in $W$.
This morphism extends to a {\em continuous}\,  morphism of {\em topological}\, props
$$
\hat{\rho}: \widehat{\Lieb} \lon  \End_V[[\hbar]],
$$
with no divergency problems.

\subsection{Etingof-Kazhdan quantizations.} There are two universal quantizations of Lie bialgebras which were constructed
in \cite{EK}: the first one involves universal formulae with traces and hence is applicable only to
finite-dimensional Lie bialgebras, while the second one avoids traces and hence is applicable to
infinite-dimensional Lie bialgebras as well. We can reinterpret these results in terms of props as follows. Denote by $\mathrm{Mor}(\mathsf{P}\rar \mathsf{Q})$ the set of continuous morphisms from a topological prop $\mathsf P$ to a topological prop $\mathsf Q$. Then a Lie bialgebra structure in $V[[\hbar]]$ is the same as an element
in  $\mathrm{Mor}(\widehat{\Lieb} \rar \mathsf{End}_{V}[[\hbar]])$ while a (not necessary commutative/ cocommutative) bialgebra structure in $\f_V[[\hbar]]$ can be understood as an element
of  $\mathrm{Mor}(\sB_\infty \rar  \mathsf{End}_{\f_V}[[\hbar]])$. The
 Etingof-Kazhdan theorem
 says that, for any graded vector space $V$, there exists a non-trivial map of sets
 $$
 \mathrm{Mor}(\widehat{\Lieb} \rar \mathsf{End}_{V}[[\hbar]]) \lon \mathrm{Mor}(\sB_\infty \rar  \mathsf{End}_{\f_V}[[\hbar]]).
 $$
 By its very definition, the prop $\DefQ$ has the following property
 $$
  \mathrm{Mor}(\sB_\infty \rar  \mathsf{End}_{\f_V}[[\hbar]])= \mathrm{Mor}(\DefQ\rar  \mathsf{End}_V[[\hbar]])
 $$
 so that we can reinterpret the Etingof-Kazhdan theorem as an existence of a non-trivial
 map
 $$
\cE\cK_V:  \mathrm{Mor}(\widehat{\Lieb} \rar \mathsf{End}_{V}[[\hbar]]) \lon  \mathrm{Mor}(\DefQ\rar  \mathsf{End}_V[[\hbar]])
 $$
 which  in fact does not depend on the choice of a particular vector space $V$ nor on  a particular
 Lie bialgebra structure in $V$, i.e.\ on a particular element in the set  $\mathrm{Mor}(\widehat{\Lieb} \rar \mathsf{End}_{V}[[\hbar]])$. Put another way, the Etingof-Kazhdan theorem can be reformulated as existence of a non-trivial  (in the sense made rigorous in item (ii) of the Theorem 5.3 below) morphism of props
 \Beq\label{5: EK morphism}
\cE\cK: (\DefQ, \sd) \lon (\widehat{\Lieb}, 0),
\Eeq
which induces the above mentioned morphism of sets.
 This observation is not new --- it was first made in \cite{EE} in a closely related (but not identical) form. Using a version  \cite{GH}  of the Etingof-Kazhdan theorem  which holds true for {\em dg}\, Lie bialgebras one can also talk about existence of a non-trivial morphism of dg props
 $$
\cE\cK^+: (\DefQ^+, \sd^+) \lon (\widehat{\Lieb}^+, \delta^+),
$$

\sip

Similarly, existence of universal quantizations of {\em finite}-dimensional Lie bialgebras can be reformulated as existence of continuous morphism of topological dg props,
$$
\cE\cK^\circlearrowright: (\DefQ, \sd) \lon (\widehat{\Lieb}^{_\circlearrowright}, 0),
$$
where $\widehat{\Lieb}^{_\circlearrowright}$  is the vertex (and hence genus) completion of the {\em wheeled}\, prop $\Lieb^\circlearrowright$ of Lie bialgebras  (see, e.g., \cite{Me2} for
its precise definition). It is existence of the  map $\cE\cK$ (or, equivalently, of $\cE\cK^+$, see \S 5.4) which allows us to prove the formality theorem
via a more less elementary lifting procedure; this iterative lifting argument does {\em not}\, work for the map $\cE\cK^\circlearrowright$.

\subsection{Theorem}\label{barF}
 {\em There exists a morphism of dg props, ${\cF}$, making the diagram
\[
 \xymatrix{
  & (\widehat{\Lieb}_\infty,\delta)
  \ar[d]^{qis} \\
 (\DefQ, \delta) \ar[ur]^{{\cF}}\ar[r]_{\cE\cK} &
 (\widehat{\Lieb}, 0)
 }
\]
commutative. Moreover, such an ${\cF}$ satisfies the following conditions:
\Bi
\item[(i)] for any generating corolla $\mathsf e$ in $\DefQ$, the composition $p_k\circ {\cF}(\mathsf e)$
is a\, {\em finite} linear combination of graphs from $\Lieb_\infty$, where $p_k$ is the projection from
$\widehat{\Lieb}_\infty$ to its subspace spanned by decorated graphs with precisely $k$ vertices.
\vspace{3mm}

\item[(ii)]
$
p_1\circ
{\cF}\left(
\xy
 <0mm,0mm>*{\mbox{$\xy *=<20mm,3mm>\txt{}*\frm{-}\endxy$}};<0mm,0mm>*{}**@{},
  <-10mm,1.5mm>*{};<-12mm,7mm>*{}**@{-},
  <-10mm,1.5mm>*{};<-11mm,7mm>*{}**@{-},
  <-10mm,1.5mm>*{};<-9.5mm,6mm>*{}**@{-},
  <-10mm,1.5mm>*{};<-8mm,7mm>*{}**@{-},
 <-10mm,1.5mm>*{};<-9.5mm,6.6mm>*{.\hspace{-0.4mm}.\hspace{-0.4mm}.}**@{},
 <0mm,0mm>*{};<-6.5mm,3.6mm>*{.\hspace{-0.1mm}.\hspace{-0.1mm}.}**@{},
  <-3mm,1.5mm>*{};<-5mm,7mm>*{}**@{-},
  <-3mm,1.5mm>*{};<-4mm,7mm>*{}**@{-},
  <-3mm,1.5mm>*{};<-2.5mm,6mm>*{}**@{-},
  <-3mm,1.5mm>*{};<-1mm,7mm>*{}**@{-},
 <-3mm,1.5mm>*{};<-2.5mm,6.6mm>*{.\hspace{-0.4mm}.\hspace{-0.4mm}.}**@{},
  <2mm,1.5mm>*{};<0mm,7mm>*{}**@{-},
  <2mm,1.5mm>*{};<1mm,7mm>*{}**@{-},
  <2mm,1.5mm>*{};<2.5mm,6mm>*{}**@{-},
  <2mm,1.5mm>*{};<4mm,7mm>*{}**@{-},
 <2mm,1.5mm>*{};<2.5mm,6.6mm>*{.\hspace{-0.4mm}.\hspace{-0.4mm}.}**@{},
 <0mm,0mm>*{};<6mm,3.6mm>*{.\hspace{-0.1mm}.\hspace{-0.1mm}.}**@{},
<10mm,1.5mm>*{};<8mm,7mm>*{}**@{-},
  <10mm,1.5mm>*{};<9mm,7mm>*{}**@{-},
  <10mm,1.5mm>*{};<10.5mm,6mm>*{}**@{-},
  <10mm,1.5mm>*{};<12mm,7mm>*{}**@{-},
 <10mm,1.5mm>*{};<10.5mm,6.6mm>*{.\hspace{-0.4mm}.\hspace{-0.4mm}.}**@{},
%
<0mm,0mm>*{};<-9.5mm,8.2mm>*{^{I_{ 1}}}**@{},
<0mm,0mm>*{};<-3mm,8.2mm>*{^{I_{ i}}}**@{},
<0mm,0mm>*{};<2mm,8.2mm>*{^{I_{ i+1}}}**@{},
<0mm,0mm>*{};<10mm,8.2mm>*{^{I_{ m}}}**@{},
<-10mm,-1.5mm>*{};<-12mm,-7mm>*{}**@{-},
  <-10mm,-1.5mm>*{};<-11mm,-7mm>*{}**@{-},
  <-10mm,-1.5mm>*{};<-9.5mm,-6mm>*{}**@{-},
  <-10mm,-1.5mm>*{};<-8mm,-7mm>*{}**@{-},
 <-10mm,-1.5mm>*{};<-9.5mm,-6.6mm>*{.\hspace{-0.4mm}.\hspace{-0.4mm}.}**@{},
 <0mm,0mm>*{};<-6.5mm,-3.6mm>*{.\hspace{-0.1mm}.\hspace{-0.1mm}.}**@{},
  <-3mm,-1.5mm>*{};<-5mm,-7mm>*{}**@{-},
  <-3mm,-1.5mm>*{};<-4mm,-7mm>*{}**@{-},
  <-3mm,-1.5mm>*{};<-2.5mm,-6mm>*{}**@{-},
  <-3mm,-1.5mm>*{};<-1mm,-7mm>*{}**@{-},
 <-3mm,-1.5mm>*{};<-2.5mm,-6.6mm>*{.\hspace{-0.4mm}.\hspace{-0.4mm}.}**@{},
  <2mm,-1.5mm>*{};<0mm,-7mm>*{}**@{-},
  <2mm,-1.5mm>*{};<1mm,-7mm>*{}**@{-},
  <2mm,-1.5mm>*{};<2.5mm,-6mm>*{}**@{-},
  <2mm,-1.5mm>*{};<4mm,-7mm>*{}**@{-},
 <2mm,-1.5mm>*{};<2.5mm,-6.6mm>*{.\hspace{-0.4mm}.\hspace{-0.4mm}.}**@{},
 <0mm,0mm>*{};<6mm,-3.6mm>*{.\hspace{-0.1mm}.\hspace{-0.1mm}.}**@{},
<10mm,-1.5mm>*{};<8mm,-7mm>*{}**@{-},
  <10mm,-1.5mm>*{};<9mm,-7mm>*{}**@{-},
  <10mm,-1.5mm>*{};<10.5mm,-6mm>*{}**@{-},
  <10mm,-1.5mm>*{};<12mm,-7mm>*{}**@{-},
 <10mm,-1.5mm>*{};<10.5mm,-6.6mm>*{.\hspace{-0.4mm}.\hspace{-0.4mm}.}**@{},
%
<0mm,0mm>*{};<-9.5mm,-9.2mm>*{^{J_{ 1}}}**@{},
<0mm,0mm>*{};<-3mm,-9.2mm>*{^{J_{ij}}}**@{},
<0mm,0mm>*{};<2mm,-9.2mm>*{^{J_{ j+1}}}**@{},
<0mm,0mm>*{};<10mm,-9.2mm>*{^{J_{ n}}}**@{},
\endxy
\right)
=\left\{\Ba{cr}
\begin{xy}
 <0mm,0mm>*{\bullet};<0mm,0mm>*{}**@{},
 <0mm,0mm>*{};<-8mm,5mm>*{}**@{-},
 <0mm,0mm>*{};<-4.5mm,5mm>*{}**@{-},
 <0mm,0mm>*{};<-1mm,5mm>*{\ldots}**@{},
 <0mm,0mm>*{};<4.5mm,5mm>*{}**@{-},
 <0mm,0mm>*{};<8mm,5mm>*{}**@{-},
   <0mm,0mm>*{};<-8.5mm,5.5mm>*{^1}**@{},
   <0mm,0mm>*{};<-5mm,5.5mm>*{^2}**@{},
   <0mm,0mm>*{};<4.5mm,5.5mm>*{^{m\hspace{-0.5mm}-\hspace{-0.5mm}1}}**@{},
   <0mm,0mm>*{};<9.0mm,5.5mm>*{^m}**@{},
 <0mm,0mm>*{};<-8mm,-5mm>*{}**@{-},
 <0mm,0mm>*{};<-4.5mm,-5mm>*{}**@{-},
 <0mm,0mm>*{};<-1mm,-5mm>*{\ldots}**@{},
 <0mm,0mm>*{};<4.5mm,-5mm>*{}**@{-},
 <0mm,0mm>*{};<8mm,-5mm>*{}**@{-},
   <0mm,0mm>*{};<-8.5mm,-6.9mm>*{^1}**@{},
   <0mm,0mm>*{};<-5mm,-6.9mm>*{^2}**@{},
   <0mm,0mm>*{};<4.5mm,-6.9mm>*{^{n\hspace{-0.5mm}-\hspace{-0.5mm}1}}**@{},
   <0mm,0mm>*{};<9.0mm,-6.9mm>*{^n}**@{},
 \end{xy}
& \mbox{for}\ |I_1|=\ldots=|I_m|=|J_1|=\ldots=|J_n|=1\\
0 & \mbox{otherwise}.
\Ea
\right.
$
\Ei

}

\begin{proof}
As $\DefQ$ is a free prop, a morphism ${\cF}$ is completely determined by its values,
$$
{\cF}\left(
\xy
 <0mm,0mm>*{\mbox{$\xy *=<20mm,3mm>\txt{}*\frm{-}\endxy$}};<0mm,0mm>*{}**@{},
  <-10mm,1.5mm>*{};<-12mm,7mm>*{}**@{-},
  <-10mm,1.5mm>*{};<-11mm,7mm>*{}**@{-},
  <-10mm,1.5mm>*{};<-9.5mm,6mm>*{}**@{-},
  <-10mm,1.5mm>*{};<-8mm,7mm>*{}**@{-},
 <-10mm,1.5mm>*{};<-9.5mm,6.6mm>*{.\hspace{-0.4mm}.\hspace{-0.4mm}.}**@{},
 <0mm,0mm>*{};<-6.5mm,3.6mm>*{.\hspace{-0.1mm}.\hspace{-0.1mm}.}**@{},
  <-3mm,1.5mm>*{};<-5mm,7mm>*{}**@{-},
  <-3mm,1.5mm>*{};<-4mm,7mm>*{}**@{-},
  <-3mm,1.5mm>*{};<-2.5mm,6mm>*{}**@{-},
  <-3mm,1.5mm>*{};<-1mm,7mm>*{}**@{-},
 <-3mm,1.5mm>*{};<-2.5mm,6.6mm>*{.\hspace{-0.4mm}.\hspace{-0.4mm}.}**@{},
  <2mm,1.5mm>*{};<0mm,7mm>*{}**@{-},
  <2mm,1.5mm>*{};<1mm,7mm>*{}**@{-},
  <2mm,1.5mm>*{};<2.5mm,6mm>*{}**@{-},
  <2mm,1.5mm>*{};<4mm,7mm>*{}**@{-},
 <2mm,1.5mm>*{};<2.5mm,6.6mm>*{.\hspace{-0.4mm}.\hspace{-0.4mm}.}**@{},
 <0mm,0mm>*{};<6mm,3.6mm>*{.\hspace{-0.1mm}.\hspace{-0.1mm}.}**@{},
<10mm,1.5mm>*{};<8mm,7mm>*{}**@{-},
  <10mm,1.5mm>*{};<9mm,7mm>*{}**@{-},
  <10mm,1.5mm>*{};<10.5mm,6mm>*{}**@{-},
  <10mm,1.5mm>*{};<12mm,7mm>*{}**@{-},
 <10mm,1.5mm>*{};<10.5mm,6.6mm>*{.\hspace{-0.4mm}.\hspace{-0.4mm}.}**@{},
%
<0mm,0mm>*{};<-9.5mm,8.2mm>*{^{I_{ 1}}}**@{},
<0mm,0mm>*{};<-3mm,8.2mm>*{^{I_{ i}}}**@{},
<0mm,0mm>*{};<2mm,8.2mm>*{^{I_{ i+1}}}**@{},
<0mm,0mm>*{};<10mm,8.2mm>*{^{I_{ m}}}**@{},
<-10mm,-1.5mm>*{};<-12mm,-7mm>*{}**@{-},
  <-10mm,-1.5mm>*{};<-11mm,-7mm>*{}**@{-},
  <-10mm,-1.5mm>*{};<-9.5mm,-6mm>*{}**@{-},
  <-10mm,-1.5mm>*{};<-8mm,-7mm>*{}**@{-},
 <-10mm,-1.5mm>*{};<-9.5mm,-6.6mm>*{.\hspace{-0.4mm}.\hspace{-0.4mm}.}**@{},
 <0mm,0mm>*{};<-6.5mm,-3.6mm>*{.\hspace{-0.1mm}.\hspace{-0.1mm}.}**@{},
  <-3mm,-1.5mm>*{};<-5mm,-7mm>*{}**@{-},
  <-3mm,-1.5mm>*{};<-4mm,-7mm>*{}**@{-},
  <-3mm,-1.5mm>*{};<-2.5mm,-6mm>*{}**@{-},
  <-3mm,-1.5mm>*{};<-1mm,-7mm>*{}**@{-},
 <-3mm,-1.5mm>*{};<-2.5mm,-6.6mm>*{.\hspace{-0.4mm}.\hspace{-0.4mm}.}**@{},
  <2mm,-1.5mm>*{};<0mm,-7mm>*{}**@{-},
  <2mm,-1.5mm>*{};<1mm,-7mm>*{}**@{-},
  <2mm,-1.5mm>*{};<2.5mm,-6mm>*{}**@{-},
  <2mm,-1.5mm>*{};<4mm,-7mm>*{}**@{-},
 <2mm,-1.5mm>*{};<2.5mm,-6.6mm>*{.\hspace{-0.4mm}.\hspace{-0.4mm}.}**@{},
 <0mm,0mm>*{};<6mm,-3.6mm>*{.\hspace{-0.1mm}.\hspace{-0.1mm}.}**@{},
<10mm,-1.5mm>*{};<8mm,-7mm>*{}**@{-},
  <10mm,-1.5mm>*{};<9mm,-7mm>*{}**@{-},
  <10mm,-1.5mm>*{};<10.5mm,-6mm>*{}**@{-},
  <10mm,-1.5mm>*{};<12mm,-7mm>*{}**@{-},
 <10mm,-1.5mm>*{};<10.5mm,-6.6mm>*{.\hspace{-0.4mm}.\hspace{-0.4mm}.}**@{},
%
<0mm,0mm>*{};<-9.5mm,-9.2mm>*{^{J_{ 1}}}**@{},
<0mm,0mm>*{};<-3mm,-9.2mm>*{^{J_{ij}}}**@{},
<0mm,0mm>*{};<2mm,-9.2mm>*{^{J_{ j+1}}}**@{},
<0mm,0mm>*{};<10mm,-9.2mm>*{^{J_{ n}}}**@{},
\endxy
\right)\in \widehat{\Lieb}_\infty, \ \ \ m+n\geq 3,
$$
 on the generating corollas. We shall construct ${\cF}$
by induction\footnote{this induction is an almost literal analogue
of the Whitehead lifting trick in the theory of $CW$-complexes in algebraic topology.}
 on the ``weight", $\mathfrak w:=m+n-3$, associated to such corollas. For $\mathfrak w=0$ we set ${\cF}$ to be
an arbitrary lift along the surjection $qis$ of the Etingof-Kazhdan morphism $\cE\cK$, i.e. we begin our induction by
setting
$$
{\cF}\left(
\xy
 <-0.5mm,0mm>*{\mbox{$\xy *=<6mm,3mm>\txt{}*\frm{-}\endxy$}};<0mm,0mm>*{}**@{},
  <-3mm,1.5mm>*{};<-5mm,7mm>*{}**@{-},
  <-3mm,1.5mm>*{};<-4mm,7mm>*{}**@{-},
  <-3mm,1.5mm>*{};<-2.5mm,6mm>*{}**@{-},
  <-3mm,1.5mm>*{};<-1mm,7mm>*{}**@{-},
 <-3mm,1.5mm>*{};<-2.5mm,6.6mm>*{.\hspace{-0.4mm}.\hspace{-0.4mm}.}**@{},
  <2mm,1.5mm>*{};<0mm,7mm>*{}**@{-},
  <2mm,1.5mm>*{};<1mm,7mm>*{}**@{-},
  <2mm,1.5mm>*{};<2.5mm,6mm>*{}**@{-},
  <2mm,1.5mm>*{};<4mm,7mm>*{}**@{-},
 <2mm,1.5mm>*{};<2.5mm,6.6mm>*{.\hspace{-0.4mm}.\hspace{-0.4mm}.}**@{},
<0mm,0mm>*{};<-3mm,8.2mm>*{^{I_{1}}}**@{},
<0mm,0mm>*{};<2mm,8.2mm>*{^{I_2}}**@{},
%
  <-0.5mm,-1.5mm>*{};<-2.5mm,-7mm>*{}**@{-},
  <-0.5mm,-1.5mm>*{};<-1.5mm,-7mm>*{}**@{-},
  <-0.5mm,-1.5mm>*{};<0mm,-6mm>*{}**@{-},
  <-0.5mm,-1.5mm>*{};<1.5mm,-7mm>*{}**@{-},
 <-0.5mm,-1.5mm>*{};<0mm,-6.6mm>*{.\hspace{-0.4mm}.\hspace{-0.4mm}.}**@{},
%
<0mm,0mm>*{};<-0.5mm,-9.2mm>*{^{J}}**@{},
\endxy
\right):= qis^{-1}\circ \cE\cK\left(
\xy
 <-0.5mm,0mm>*{\mbox{$\xy *=<6mm,3mm>\txt{}*\frm{-}\endxy$}};<0mm,0mm>*{}**@{},
  <-3mm,1.5mm>*{};<-5mm,7mm>*{}**@{-},
  <-3mm,1.5mm>*{};<-4mm,7mm>*{}**@{-},
  <-3mm,1.5mm>*{};<-2.5mm,6mm>*{}**@{-},
  <-3mm,1.5mm>*{};<-1mm,7mm>*{}**@{-},
 <-3mm,1.5mm>*{};<-2.5mm,6.6mm>*{.\hspace{-0.4mm}.\hspace{-0.4mm}.}**@{},
  <2mm,1.5mm>*{};<0mm,7mm>*{}**@{-},
  <2mm,1.5mm>*{};<1mm,7mm>*{}**@{-},
  <2mm,1.5mm>*{};<2.5mm,6mm>*{}**@{-},
  <2mm,1.5mm>*{};<4mm,7mm>*{}**@{-},
 <2mm,1.5mm>*{};<2.5mm,6.6mm>*{.\hspace{-0.4mm}.\hspace{-0.4mm}.}**@{},
<0mm,0mm>*{};<-3mm,8.2mm>*{^{I_{1}}}**@{},
<0mm,0mm>*{};<2mm,8.2mm>*{^{I_2}}**@{},
%
  <-0.5mm,-1.5mm>*{};<-2.5mm,-7mm>*{}**@{-},
  <-0.5mm,-1.5mm>*{};<-1.5mm,-7mm>*{}**@{-},
  <-0.5mm,-1.5mm>*{};<0mm,-6mm>*{}**@{-},
  <-0.5mm,-1.5mm>*{};<1.5mm,-7mm>*{}**@{-},
 <-0.5mm,-1.5mm>*{};<0mm,-6.6mm>*{.\hspace{-0.4mm}.\hspace{-0.4mm}.}**@{},
%
<0mm,0mm>*{};<-0.5mm,-9.2mm>*{^{J}}**@{},
\endxy
\right),
 \ \ \ \mbox{and}\ \ \
{\cF}\left(
\xy
 <-0.5mm,0mm>*{\mbox{$\xy *=<6mm,3mm>\txt{}*\frm{-}\endxy$}};<0mm,0mm>*{}**@{},
  <-3mm,-1.5mm>*{};<-5mm,-7mm>*{}**@{-},
  <-3mm,-1.5mm>*{};<-4mm,-7mm>*{}**@{-},
  <-3mm,-1.5mm>*{};<-2.5mm,-6mm>*{}**@{-},
  <-3mm,-1.5mm>*{};<-1mm,-7mm>*{}**@{-},
 <-3mm,-1.5mm>*{};<-2.5mm,-6.6mm>*{.\hspace{-0.4mm}.\hspace{-0.4mm}.}**@{},
  <2mm,-1.5mm>*{};<0mm,-7mm>*{}**@{-},
  <2mm,-1.5mm>*{};<1mm,-7mm>*{}**@{-},
  <2mm,-1.5mm>*{};<2.5mm,-6mm>*{}**@{-},
  <2mm,-1.5mm>*{};<4mm,-7mm>*{}**@{-},
 <2mm,-1.5mm>*{};<2.5mm,-6.6mm>*{.\hspace{-0.4mm}.\hspace{-0.4mm}.}**@{},
<0mm,0mm>*{};<-3mm,-9.2mm>*{^{J_{1}}}**@{},
<0mm,0mm>*{};<2mm,-9.2mm>*{^{J_2}}**@{},
%
  <-0.5mm,1.5mm>*{};<-2.5mm,7mm>*{}**@{-},
  <-0.5mm,1.5mm>*{};<-1.5mm,7mm>*{}**@{-},
  <-0.5mm,1.5mm>*{};<0mm,6mm>*{}**@{-},
  <-0.5mm,1.5mm>*{};<1.5mm,7mm>*{}**@{-},
 <-0.5mm,1.5mm>*{};<0mm,6.6mm>*{.\hspace{-0.4mm}.\hspace{-0.4mm}.}**@{},
%
<0mm,0mm>*{};<-0.5mm,9.2mm>*{^{I}}**@{},
\endxy
\right):=
qis^{-1}\circ \cE\cK\left(
\xy
 <-0.5mm,0mm>*{\mbox{$\xy *=<6mm,3mm>\txt{}*\frm{-}\endxy$}};<0mm,0mm>*{}**@{},
  <-3mm,1.5mm>*{};<-5mm,7mm>*{}**@{-},
  <-3mm,1.5mm>*{};<-4mm,7mm>*{}**@{-},
  <-3mm,1.5mm>*{};<-2.5mm,6mm>*{}**@{-},
  <-3mm,1.5mm>*{};<-1mm,7mm>*{}**@{-},
 <-3mm,1.5mm>*{};<-2.5mm,6.6mm>*{.\hspace{-0.4mm}.\hspace{-0.4mm}.}**@{},
  <2mm,1.5mm>*{};<0mm,7mm>*{}**@{-},
  <2mm,1.5mm>*{};<1mm,7mm>*{}**@{-},
  <2mm,1.5mm>*{};<2.5mm,6mm>*{}**@{-},
  <2mm,1.5mm>*{};<4mm,7mm>*{}**@{-},
 <2mm,1.5mm>*{};<2.5mm,6.6mm>*{.\hspace{-0.4mm}.\hspace{-0.4mm}.}**@{},
<0mm,0mm>*{};<-3mm,8.2mm>*{^{I_{1}}}**@{},
<0mm,0mm>*{};<2mm,8.2mm>*{^{I_2}}**@{},
%
  <-0.5mm,-1.5mm>*{};<-2.5mm,-7mm>*{}**@{-},
  <-0.5mm,-1.5mm>*{};<-1.5mm,-7mm>*{}**@{-},
  <-0.5mm,-1.5mm>*{};<0mm,-6mm>*{}**@{-},
  <-0.5mm,-1.5mm>*{};<1.5mm,-7mm>*{}**@{-},
 <-0.5mm,-1.5mm>*{};<0mm,-6.6mm>*{.\hspace{-0.4mm}.\hspace{-0.4mm}.}**@{},
%
<0mm,0mm>*{};<-0.5mm,-9.2mm>*{^{J}}**@{},
\endxy
\right),
$$
where $qis^{-1}$ is an arbitrary section of the quasi-isomorphism $qis$, i.e.\ an arbitrary lifting
of cohomology classes into cycles.

Assume we constructed values of ${\cF}$ on all corollas of weight $\mathfrak w\leq N$.
Let $\mathsf e$ be a generating corolla of $\DefQ$ with non-zero
weight
$\mathfrak w=N+1$. Note that $\sd \mathsf e$ is a linear combination of graphs whose vertices are decorated by
corollas of weight
$\leq N$ (as $\mathfrak w$ is the precisely minus of the degree of $\mathsf e$, and $\sd$ increases degree by $+1$).
By induction, ${\cF}(\sd \mathsf e)$
is a well-defined element in $\widehat{\Lieb}_\infty$. As $\cE\cK(\mathsf e)=0$, the element,
$$
{\cF}(\sd \mathsf e)
$$
is a closed element in $\widehat{\Lieb}_\infty$ which projects under $qis$ to zero. Since the surjection $qis$ is a
quasi-isomorphism, this element must be exact. Thus there exists $\mathfrak e\in \widehat{\Lieb}_\infty$ such that
$$
\delta {\mathfrak e}= {\cF}(\sd\mathsf e).
$$
We set ${\cF}(\mathsf e):= {\mathfrak e}$ completing thereby inductive construction of ${\cF}$.

\sip

Next, if $\mathsf e$ is a generating corolla in $\Defq$ of degree $3-m-n$ then $p_k\circ {\cF}(\mathsf e)$
is a linear combination of graphs built from $k$ generating corollas in $\Lieb_\infty$ with total number of half edges
attached to these $k$ vertices being equal  $3(k-1)+m+n$. There is only a {\em finite}\, number of
graphs in $\Lieb_\infty$
satisfying these conditions. This proves Claim (i).

\sip

Claim (ii) is obvious in the part {\em otherwise}.
Let
$$
e_n^m:= \sum_{\sigma\in \bS_m} \sum_{\tau\in \bS_n}(-1)^{sgn(\sigma)+sgn(\tau)}
\xy
 <0mm,0mm>*{\mbox{$\xy *=<20mm,3mm>\txt{}*\frm{-}\endxy$}};<0mm,0mm>*{}**@{},
  <-10mm,1.5mm>*{};<-12mm,7mm>*{}**@{-},
 <0mm,0mm>*{};<-6.5mm,3.6mm>*{.\hspace{-0.1mm}.\hspace{-0.1mm}.}**@{},
  <-3mm,1.5mm>*{};<-3.5mm,7mm>*{}**@{-},
  <3mm,1.5mm>*{};<3.5mm,7mm>*{}**@{-},
 <0mm,0mm>*{};<6mm,3.6mm>*{.\hspace{-0.1mm}.\hspace{-0.1mm}.}**@{},
<10mm,1.5mm>*{};<11.5mm,7mm>*{}**@{-},
%
<0mm,0mm>*{};<-14mm,8.2mm>*{^{\sigma(1)}}**@{},
<0mm,0mm>*{};<-5mm,8.2mm>*{^{\sigma(i)}}**@{},
<0mm,0mm>*{};<4mm,8.2mm>*{^{\sigma(i\hspace{-0.1mm}+\hspace{-0.1mm}1)}}**@{},
<0mm,0mm>*{};<14mm,8.2mm>*{^{\sigma(m)}}**@{},
 <-10mm,-1.5mm>*{};<-12mm,-7mm>*{}**@{-},
 <0mm,0mm>*{};<-6.5mm,-3.6mm>*{.\hspace{-0.1mm}.\hspace{-0.1mm}.}**@{},
  <-3mm,-1.5mm>*{};<-3.5mm,-7mm>*{}**@{-},
  <3mm,-1.5mm>*{};<3.5mm,-7mm>*{}**@{-},
 <0mm,0mm>*{};<6mm,-3.6mm>*{.\hspace{-0.1mm}.\hspace{-0.1mm}.}**@{},
<10mm,-1.5mm>*{};<11.5mm,-7mm>*{}**@{-},
%
<0mm,0mm>*{};<-14mm,-9.2mm>*{^{\tau(1)}}**@{},
<0mm,0mm>*{};<-5mm,-9.2mm>*{^{\tau(i)}}**@{},
<0mm,0mm>*{};<4mm,-9.2mm>*{^{\tau(i\hspace{-0.1mm}+\hspace{-0.1mm}1)}}**@{},
<0mm,0mm>*{};<14mm,-9.2mm>*{^{\tau(m)}}**@{},
\endxy
$$
 be the skewsymmetrization  of the generating corolla in $\DefQ$
with $ |I_1|=\ldots=|I_m|=|J_1|=\ldots=|J_n|=1$. To prove Claim~(ii) it is enough to show that
$p_1\circ {\cF}(e_n^m)\neq 0$ for all $m,n\geq 1$, $m+n\geq 3$. We shall show this by induction on the weight
$\mathfrak w=m+n-3$.

Denote the composition $p_k\circ {\cF}$ by ${\cF}_k$.

If $m+n=3$, then  ${\cF}_1(e_n^m)\neq 0$ by the construction of ${\cF}$.

 Assume that   ${\cF}_1(e_p^q)\neq 0$ for
all  $e_p^q$ with weight $\mathfrak w\leq N$ and consider  $e_n^m$ with
non-zero weight $\mathfrak w=N+1$. Then
\[
{\cF}_1(e_n^m)\neq 0      \ \  \Leftrightarrow \delta\left({\cF}_1(e_n^m)\right)\neq 0.
\]
By our construction of ${\cF}$, we have
\[
 \delta({\cF}_1(e_n^m)) ={\cF}_2(\sd e_n^m)={\cF}_2(\sd_1 e_n^m) +
{\cF}_1\boxtimes {\cF}_1(\sd_2 e_n^m),
\]
where $\sd_1$ is the linear in number of vertices part of the differential $\sd$ in $\DefQ$ and is given by (\ref{sd_1}),
  $\sd_2$ stands for the quadratic (i.e.\ spanned by two-vertex graphs)
part of  $\sd$,
and
${\cF}_1\boxtimes {\cF}_1$ means the morphism ${\cF}_1$ applied to
decoration of each of the two vertices in every
graph summand of $\sd_2e_m^n$. Now $\sd_1 e_n^m=0$, while  $\sd_2 e_m^n$  contains, for example,
 the following linear combination
of graphs,
$$
 \sum_{\sigma\in \bS_m} \sum_{\tau\in \bS_n}(-1)^{sgn(\sigma)+sgn(\tau)}
\xy
 <0mm,3mm>*{\mbox{$\xy *=<20mm,3mm>\txt{}*\frm{-}\endxy$}};<0mm,0mm>*{}**@{},
  <-10mm,4.5mm>*{};<-12mm,10mm>*{}**@{-},
 <0mm,3mm>*{};<-6.5mm,6.6mm>*{.\hspace{-0.1mm}.\hspace{-0.1mm}.}**@{},
  <-3mm,4.5mm>*{};<-3.5mm,10mm>*{}**@{-},
  <3mm,4.5mm>*{};<3.5mm,10mm>*{}**@{-},
 <0mm,3mm>*{};<6mm,6.6mm>*{.\hspace{-0.1mm}.\hspace{-0.1mm}.}**@{},
<10mm,4.5mm>*{};<11.5mm,10mm>*{}**@{-},
%
<0mm,0mm>*{};<-14mm,11.2mm>*{^{\sigma(1)}}**@{},
<0mm,0mm>*{};<-5mm,11.2mm>*{^{\sigma(i)}}**@{},
<0mm,0mm>*{};<4mm,11.2mm>*{^{\sigma(i\hspace{-0.1mm}+\hspace{-0.1mm}1)}}**@{},
<0mm,0mm>*{};<14mm,11.2mm>*{^{\sigma(m)}}**@{},
 <0mm,-3mm>*{\mbox{$\xy *=<20mm,3mm>\txt{}*\frm{-}\endxy$}};<0mm,0mm>*{}**@{},
  <-10mm,-4.5mm>*{};<-12mm,-10mm>*{}**@{-},
 <0mm,-3mm>*{};<-6.5mm,-6.6mm>*{.\hspace{-0.1mm}.\hspace{-0.1mm}.}**@{},
  <-3mm,-4.5mm>*{};<-3.5mm,-10mm>*{}**@{-},
  <3mm,-4.5mm>*{};<3.5mm,-10mm>*{}**@{-},
 <0mm,-3mm>*{};<6mm,-6.6mm>*{.\hspace{-0.1mm}.\hspace{-0.1mm}.}**@{},
<10mm,-4.5mm>*{};<11.5mm,-10mm>*{}**@{-},
<0mm,-1.5mm>*{};<0mm,1.5mm>*{}**@{-},
<0mm,0mm>*{};<-14mm,-12.2mm>*{^{\tau(1)}}**@{},
<0mm,0mm>*{};<-5mm,-12.2mm>*{^{\tau(i)}}**@{},
<0mm,0mm>*{};<4mm,-12.2mm>*{^{\tau(i\hspace{-0.1mm}+\hspace{-0.1mm}1)}}**@{},
<0mm,0mm>*{};<14mm,-12.2mm>*{^{\tau(m)}}**@{},
\endxy
\ ,
$$
as irreducible summands coming from the genus zero part of the differential $\delta$ in $\sB_\infty$ (see \cite{Ma}).
Then, by the induction assumption, ${\cF}_1\boxtimes {\cF}_1(\sd_2 e_n^m)$ can not be zero
implying ${\cF}_1(e_n^m)\neq 0$. This completes the proof of Claim (ii) and hence of the Theorem.
\end{proof}

\subsection{First Proof of Formality Theorem~\ref{1.2}} We are going to construct a sequence of linear maps,
$$
F_k: \wedge^k (\fl_V[2]) \lon \mathsf{poly}(\f_V,\f_V), \ \ k\geq 1,
$$
of degree $1-k$ satisfying quadratic relations of an $L_\infty$-morphism,
\Beq\label{morphism}
\sum_{\sigma\in Sh(2)}F_{k-1}\left(\{f_{\sigma(1)}, f_{\sigma(2)}\}, f_{\sigma(3)}, \ldots, f_{\sigma(k)}\right)
=\hspace{60mm}
\Eeq
\[
\hspace{15mm}
=
\sum_{i=1}^k \sum_{k=k_1+\ldots + k_i}\sum_{\sigma\in Sh(k_1,\ldots, k_r)} \pm \mu_r(F_{k_1}\left(f_{\sigma(1)}, \ldots, f_{\sigma(k_1)}),
\dots, F_{k_r}(f_{\sigma(k-k_r+1)}, \ldots, f_{\sigma(k)}\right),
\]
 where $ Sh(k_1,\ldots, k_r)$ stands for the subgroup of $(k_1,\ldots, k_r)$-shuffles
in $\bS_k$, and $f_1, \ldots, f_k$ are arbitrary elements in $\fl_V[2]$. In a linear coordinate system $\{x^j\}$
on $V$ (and the dual coordinate system $\{p_i\}$ on $V^*$) such an element $f$ is a formal power series,
$$
f=\sum_{m,n\geq 1} f^{i_1\ldots i_m}_{j_1\ldots j_n} x^{j_1}\wedge\ldots \wedge x^{j_n}\wedge p_{j_1}\wedge\ldots
\wedge p_{j_m}.
$$
We define a degree $0$ linear map
$$
\Ba{rccc}
F_1: & \fl_V[2] & \lon & \mathsf{poly}(\f_V,\f_V)\\
     & f & \lon & F_1(f)
\Ea
$$
by setting
$$
F_1(f):=\sum_{m,n\geq 1}  f^{i_1\ldots i_m}_{j_1\ldots j_n} x^{j_1}\ot \ldots \ot x^{j_n}
\cdot \Delta^{n-1}\left( \frac{\p}{\p x^{i_1}}\right)\cdot\ldots \cdot
 \Delta^{n-1}\left( \frac{\p}{\p x^{i_m}}\right).
$$
Clearly, $d_{\fg\fs}\circ F_1=0$ so that this map sends $\fl_V[2]$ into cycles in $\mathsf{poly}(\f_V,\f_V)$.

Next we shall read off the maps $F_k$ for $k\geq 2$ from the components, ${\cF_k}$, of the morphism ${\cF}$ we
constructed in the proof of Theorem~\ref{barF}. As the morphism ${\cF}$ lands in the prop which does not contain
the $(1,1)$-corolla
$
\begin{xy}
 <0mm,-0.55mm>*{};<0mm,-3mm>*{}**@{-},
 <0mm,0.5mm>*{};<0mm,3mm>*{}**@{-},
 <0mm,0mm>*{\bullet};<0mm,0mm>*{}**@{},
 \end{xy},
$
we are forced to set
$$
F_k(f_1, \ldots, f_k)=0,  \ \ \ k\geq 2,
$$
if at least one input $f_i$ lies in $V[1]\ot V^*[1]$. Thus $F_k$ for $k\geq 2$ must  factor through the projection,
$$
\wedge^k (\fl_V[2]) \lon \wedge^k \fg \stackrel{F_k}{\lon} \mathsf{poly}(\f_V,\f_V),
$$
where
$$
\fg:= \bigoplus_{m,n\geq 1\atop m+n\geq 3} \wedge^m V\ot \wedge^n V^*[2-m-n] \subset \fl_V[2].
$$
From now on we identify (see Corollary 5.1 in \cite{Me1}) the latter with the vector space of $\bS$-equivariant linear maps,
$$
\fg\equiv \Hom(\mathsf L, \mathsf{End}_V)[-1]\simeq \Def(\Lieb_\infty \stackrel{0}{\rar} \mathsf{End}_V)
$$
where  the $\bS$-bimodule $\mathsf L$ is given by (\ref{L})).


The maps $F_k$ are defined once  the compositions,
\[
{F}_{k, |J_1|, \ldots, |J_n|}^{\ \, |I_1|, \ldots, |I_m|}: \wedge^k \fg \stackrel{F_k}{\lon} \mathsf{poly}(\f_V,\f_V)
\stackrel{proj}{\lon}\Hom_{poly}(\f_V^{\ot n}, \f_V^{\ot m})
\stackrel{proj}{\lon} \hspace{60mm}
\]
\[
\hspace{70mm}
\stackrel{proj}{\lon}\Hom\left(\odot^{|J_1|}V\ot\ldots \ot \odot^{|J_n|}V, \odot^{|I_1|}\ot\ldots \ot \odot^{|I_m|}V\right).
\]
are defined for all $k\geq 2$, $m+n\geq 3$, and $|J_1|,\ldots, |J_n|, |I_1|, \ldots, |I_m|\geq 1$.

We construct the value $F_{k,|J_1|, \ldots, |J_n|}^{\ \, |I_1|, \ldots, |I_m|}(f_1, \ldots, f_k)$ on arbitrary
$f_1, \ldots, f_k\in  \Hom(\mathsf L, \End\langle V\rangle)[-1]$ in three steps:

\noindent {\em Step 1}. Consider
$$
{\cF}_k
\left(
\xy
 <0mm,0mm>*{\mbox{$\xy *=<20mm,3mm>\txt{}*\frm{-}\endxy$}};<0mm,0mm>*{}**@{},
  <-10mm,1.5mm>*{};<-12mm,7mm>*{}**@{-},
  <-10mm,1.5mm>*{};<-11mm,7mm>*{}**@{-},
  <-10mm,1.5mm>*{};<-9.5mm,6mm>*{}**@{-},
  <-10mm,1.5mm>*{};<-8mm,7mm>*{}**@{-},
 <-10mm,1.5mm>*{};<-9.5mm,6.6mm>*{.\hspace{-0.4mm}.\hspace{-0.4mm}.}**@{},
 <0mm,0mm>*{};<-6.5mm,3.6mm>*{.\hspace{-0.1mm}.\hspace{-0.1mm}.}**@{},
  <-3mm,1.5mm>*{};<-5mm,7mm>*{}**@{-},
  <-3mm,1.5mm>*{};<-4mm,7mm>*{}**@{-},
  <-3mm,1.5mm>*{};<-2.5mm,6mm>*{}**@{-},
  <-3mm,1.5mm>*{};<-1mm,7mm>*{}**@{-},
 <-3mm,1.5mm>*{};<-2.5mm,6.6mm>*{.\hspace{-0.4mm}.\hspace{-0.4mm}.}**@{},
  <2mm,1.5mm>*{};<0mm,7mm>*{}**@{-},
  <2mm,1.5mm>*{};<1mm,7mm>*{}**@{-},
  <2mm,1.5mm>*{};<2.5mm,6mm>*{}**@{-},
  <2mm,1.5mm>*{};<4mm,7mm>*{}**@{-},
 <2mm,1.5mm>*{};<2.5mm,6.6mm>*{.\hspace{-0.4mm}.\hspace{-0.4mm}.}**@{},
 <0mm,0mm>*{};<6mm,3.6mm>*{.\hspace{-0.1mm}.\hspace{-0.1mm}.}**@{},
<10mm,1.5mm>*{};<8mm,7mm>*{}**@{-},
  <10mm,1.5mm>*{};<9mm,7mm>*{}**@{-},
  <10mm,1.5mm>*{};<10.5mm,6mm>*{}**@{-},
  <10mm,1.5mm>*{};<12mm,7mm>*{}**@{-},
 <10mm,1.5mm>*{};<10.5mm,6.6mm>*{.\hspace{-0.4mm}.\hspace{-0.4mm}.}**@{},
%
<0mm,0mm>*{};<-9.5mm,8.2mm>*{^{I_{ 1}}}**@{},
<0mm,0mm>*{};<-3mm,8.2mm>*{^{I_{ i}}}**@{},
<0mm,0mm>*{};<2mm,8.2mm>*{^{I_{ i+1}}}**@{},
<0mm,0mm>*{};<10mm,8.2mm>*{^{I_{ m}}}**@{},
<-10mm,-1.5mm>*{};<-12mm,-7mm>*{}**@{-},
  <-10mm,-1.5mm>*{};<-11mm,-7mm>*{}**@{-},
  <-10mm,-1.5mm>*{};<-9.5mm,-6mm>*{}**@{-},
  <-10mm,-1.5mm>*{};<-8mm,-7mm>*{}**@{-},
 <-10mm,-1.5mm>*{};<-9.5mm,-6.6mm>*{.\hspace{-0.4mm}.\hspace{-0.4mm}.}**@{},
 <0mm,0mm>*{};<-6.5mm,-3.6mm>*{.\hspace{-0.1mm}.\hspace{-0.1mm}.}**@{},
  <-3mm,-1.5mm>*{};<-5mm,-7mm>*{}**@{-},
  <-3mm,-1.5mm>*{};<-4mm,-7mm>*{}**@{-},
  <-3mm,-1.5mm>*{};<-2.5mm,-6mm>*{}**@{-},
  <-3mm,-1.5mm>*{};<-1mm,-7mm>*{}**@{-},
 <-3mm,-1.5mm>*{};<-2.5mm,-6.6mm>*{.\hspace{-0.4mm}.\hspace{-0.4mm}.}**@{},
  <2mm,-1.5mm>*{};<0mm,-7mm>*{}**@{-},
  <2mm,-1.5mm>*{};<1mm,-7mm>*{}**@{-},
  <2mm,-1.5mm>*{};<2.5mm,-6mm>*{}**@{-},
  <2mm,-1.5mm>*{};<4mm,-7mm>*{}**@{-},
 <2mm,-1.5mm>*{};<2.5mm,-6.6mm>*{.\hspace{-0.4mm}.\hspace{-0.4mm}.}**@{},
 <0mm,0mm>*{};<6mm,-3.6mm>*{.\hspace{-0.1mm}.\hspace{-0.1mm}.}**@{},
<10mm,-1.5mm>*{};<8mm,-7mm>*{}**@{-},
  <10mm,-1.5mm>*{};<9mm,-7mm>*{}**@{-},
  <10mm,-1.5mm>*{};<10.5mm,-6mm>*{}**@{-},
  <10mm,-1.5mm>*{};<12mm,-7mm>*{}**@{-},
 <10mm,-1.5mm>*{};<10.5mm,-6.6mm>*{.\hspace{-0.4mm}.\hspace{-0.4mm}.}**@{},
%
<0mm,0mm>*{};<-9.5mm,-9.2mm>*{^{J_{ 1}}}**@{},
<0mm,0mm>*{};<-3mm,-9.2mm>*{^{J_{ij}}}**@{},
<0mm,0mm>*{};<2mm,-9.2mm>*{^{J_{ j+1}}}**@{},
<0mm,0mm>*{};<10mm,-9.2mm>*{^{J_{ n}}}**@{},
\endxy
\right)=
\sum_{G} G\langle l_{1}, l_{2}, \ldots, l_{k}\rangle _{Aut(G)} \in \Lieb_\infty
$$
where the sun runs over a family of graphs with $k$ vertices, $Vert(G)=\{v_1, \dots, v_k\}$, decorated with
with the unique generators
$l_i\in \mathsf L(|In(v_i)|, |Out(v_i)|$, $i=1, \ldots, k$ from $\mathsf{L}$,
 where $|In(v_i)|$ (resp.\  $|Out(v_i)|$) stands for the number
of input (resp.\ output) half edges attached to the vertex $v_i$ and
$$
G\langle l_{v_1}, l_{v_2},
\ldots, l_{v_k}\rangle
:=\left(\sum_{s: [k]\rar Vert(G)} l_{s(1)}\ot \ldots \ot l_{s(k)}
\right)_{\bS_k}
$$
stands for the {\em unordered}\, tensor product
of $l_i$ over the set $Vert(G)$.

\noindent {\em Step 2}. Define next
$$
{\cF}_k(f_1,\ldots, f_k):= \sum_G\sum_{\sigma\in \bS_k} (-1)^{Koszul\ sgn}
 G\langle f_{\sigma(1)}(l_1),   f_{\sigma(2)}(l_2), \ldots,   f_{\sigma(k)}(l_k)\rangle _{Aut(G)}.
$$
This is a sum of graphs whose vertices are decorated by elements of the endomorphism prop
$\mathsf{End}_V$ (cf. \cite{Ma3,MV}). Hence we can apply to $\cF_k(f_1,\ldots, f_k)$ the vertical and horizontal
$\mathsf{End}_V$-compositions
to get a well defined element,
$$
comp_{\mathsf{End}_V}\left({\cF}_k(f_1,\ldots, f_k)\right)\in \End\langle V\rangle,
$$
which, in fact, lies by the  construction in the subspace
$$
\Hom\left(\odot^{|J_1|}V\ot\ldots \ot \odot^{|J_n|}V, \odot^{|I_1|}V\ot\ldots \ot \odot^{|I_m|}V\right)\subset
\mathsf{End}_V.
$$

{\em Step 3}. Finally  set
$$
F_{k,|J_1|, \ldots, |J_n|}^{\ \, |I_1|, \ldots, |I_m|}(f_1, \ldots, f_k):=
comp_{\End\langle V\rangle}\left({\cF}_k(l_1,\ldots, l_k)\right).
$$

It is easy to check that the constructed map $F_k$ has degree $1-k$. It is also a straightforward
untwisting of the definitions of differentials $\sd$ in $\Defq$ and $\delta$ in $\Lieb_\infty$ to show that
equations $(\ref{morphism})$ follow directly from the basic property, $\delta\circ {\cF}= {\cF}\circ \sd$,
of the morphism  ${\cF}$ (cf.\ again \cite{Ma2,MV}).
\hfill $\Box$

\subsection{Second Proof of Formality Theorem~\ref{1.2}}
There exists a commutative diagram of dg props,
\[
 \xymatrix{
  & (\widehat{\Lieb}^{_+}_\infty,\delta^+)
  \ar[d]^{qis} \\
 (\DefQ^+, \delta^+) \ar[ur]^{{\cF^+}}\ar[r]_{\cE\cK^+} &
 (\widehat{\Lieb}^{_+}, \delta^+)
 }
\]
with
\[
p_1\circ
{\cF^+}\left(
\xy
 <0mm,0mm>*{\mbox{$\xy *=<20mm,3mm>\txt{}*\frm{-}\endxy$}};<0mm,0mm>*{}**@{},
  <-10mm,1.5mm>*{};<-12mm,7mm>*{}**@{-},
  <-10mm,1.5mm>*{};<-11mm,7mm>*{}**@{-},
  <-10mm,1.5mm>*{};<-9.5mm,6mm>*{}**@{-},
  <-10mm,1.5mm>*{};<-8mm,7mm>*{}**@{-},
 <-10mm,1.5mm>*{};<-9.5mm,6.6mm>*{.\hspace{-0.4mm}.\hspace{-0.4mm}.}**@{},
 <0mm,0mm>*{};<-6.5mm,3.6mm>*{.\hspace{-0.1mm}.\hspace{-0.1mm}.}**@{},
  <-3mm,1.5mm>*{};<-5mm,7mm>*{}**@{-},
  <-3mm,1.5mm>*{};<-4mm,7mm>*{}**@{-},
  <-3mm,1.5mm>*{};<-2.5mm,6mm>*{}**@{-},
  <-3mm,1.5mm>*{};<-1mm,7mm>*{}**@{-},
 <-3mm,1.5mm>*{};<-2.5mm,6.6mm>*{.\hspace{-0.4mm}.\hspace{-0.4mm}.}**@{},
  <2mm,1.5mm>*{};<0mm,7mm>*{}**@{-},
  <2mm,1.5mm>*{};<1mm,7mm>*{}**@{-},
  <2mm,1.5mm>*{};<2.5mm,6mm>*{}**@{-},
  <2mm,1.5mm>*{};<4mm,7mm>*{}**@{-},
 <2mm,1.5mm>*{};<2.5mm,6.6mm>*{.\hspace{-0.4mm}.\hspace{-0.4mm}.}**@{},
 <0mm,0mm>*{};<6mm,3.6mm>*{.\hspace{-0.1mm}.\hspace{-0.1mm}.}**@{},
<10mm,1.5mm>*{};<8mm,7mm>*{}**@{-},
  <10mm,1.5mm>*{};<9mm,7mm>*{}**@{-},
  <10mm,1.5mm>*{};<10.5mm,6mm>*{}**@{-},
  <10mm,1.5mm>*{};<12mm,7mm>*{}**@{-},
 <10mm,1.5mm>*{};<10.5mm,6.6mm>*{.\hspace{-0.4mm}.\hspace{-0.4mm}.}**@{},
%
<0mm,0mm>*{};<-9.5mm,8.2mm>*{^{I_{ 1}}}**@{},
<0mm,0mm>*{};<-3mm,8.2mm>*{^{I_{ i}}}**@{},
<0mm,0mm>*{};<2mm,8.2mm>*{^{I_{ i+1}}}**@{},
<0mm,0mm>*{};<10mm,8.2mm>*{^{I_{ m}}}**@{},
<-10mm,-1.5mm>*{};<-12mm,-7mm>*{}**@{-},
  <-10mm,-1.5mm>*{};<-11mm,-7mm>*{}**@{-},
  <-10mm,-1.5mm>*{};<-9.5mm,-6mm>*{}**@{-},
  <-10mm,-1.5mm>*{};<-8mm,-7mm>*{}**@{-},
 <-10mm,-1.5mm>*{};<-9.5mm,-6.6mm>*{.\hspace{-0.4mm}.\hspace{-0.4mm}.}**@{},
 <0mm,0mm>*{};<-6.5mm,-3.6mm>*{.\hspace{-0.1mm}.\hspace{-0.1mm}.}**@{},
  <-3mm,-1.5mm>*{};<-5mm,-7mm>*{}**@{-},
  <-3mm,-1.5mm>*{};<-4mm,-7mm>*{}**@{-},
  <-3mm,-1.5mm>*{};<-2.5mm,-6mm>*{}**@{-},
  <-3mm,-1.5mm>*{};<-1mm,-7mm>*{}**@{-},
 <-3mm,-1.5mm>*{};<-2.5mm,-6.6mm>*{.\hspace{-0.4mm}.\hspace{-0.4mm}.}**@{},
  <2mm,-1.5mm>*{};<0mm,-7mm>*{}**@{-},
  <2mm,-1.5mm>*{};<1mm,-7mm>*{}**@{-},
  <2mm,-1.5mm>*{};<2.5mm,-6mm>*{}**@{-},
  <2mm,-1.5mm>*{};<4mm,-7mm>*{}**@{-},
 <2mm,-1.5mm>*{};<2.5mm,-6.6mm>*{.\hspace{-0.4mm}.\hspace{-0.4mm}.}**@{},
 <0mm,0mm>*{};<6mm,-3.6mm>*{.\hspace{-0.1mm}.\hspace{-0.1mm}.}**@{},
<10mm,-1.5mm>*{};<8mm,-7mm>*{}**@{-},
  <10mm,-1.5mm>*{};<9mm,-7mm>*{}**@{-},
  <10mm,-1.5mm>*{};<10.5mm,-6mm>*{}**@{-},
  <10mm,-1.5mm>*{};<12mm,-7mm>*{}**@{-},
 <10mm,-1.5mm>*{};<10.5mm,-6.6mm>*{.\hspace{-0.4mm}.\hspace{-0.4mm}.}**@{},
%
<0mm,0mm>*{};<-9.5mm,-9.2mm>*{^{J_{ 1}}}**@{},
<0mm,0mm>*{};<-3mm,-9.2mm>*{^{J_{j}}}**@{},
<0mm,0mm>*{};<2mm,-9.2mm>*{^{J_{ j+1}}}**@{},
<0mm,0mm>*{};<10mm,-9.2mm>*{^{J_{ n}}}**@{},
\endxy
\right)
=\left\{\Ba{cr}
\begin{xy}
 <0mm,0mm>*{\bullet};<0mm,0mm>*{}**@{},
 <0mm,0mm>*{};<-8mm,5mm>*{}**@{-},
 <0mm,0mm>*{};<-4.5mm,5mm>*{}**@{-},
 <0mm,0mm>*{};<-1mm,5mm>*{\ldots}**@{},
 <0mm,0mm>*{};<4.5mm,5mm>*{}**@{-},
 <0mm,0mm>*{};<8mm,5mm>*{}**@{-},
   <0mm,0mm>*{};<-8.5mm,5.5mm>*{^1}**@{},
   <0mm,0mm>*{};<-5mm,5.5mm>*{^2}**@{},
   <0mm,0mm>*{};<4.5mm,5.5mm>*{^{m\hspace{-0.5mm}-\hspace{-0.5mm}1}}**@{},
   <0mm,0mm>*{};<9.0mm,5.5mm>*{^m}**@{},
 <0mm,0mm>*{};<-8mm,-5mm>*{}**@{-},
 <0mm,0mm>*{};<-4.5mm,-5mm>*{}**@{-},
 <0mm,0mm>*{};<-1mm,-5mm>*{\ldots}**@{},
 <0mm,0mm>*{};<4.5mm,-5mm>*{}**@{-},
 <0mm,0mm>*{};<8mm,-5mm>*{}**@{-},
   <0mm,0mm>*{};<-8.5mm,-6.9mm>*{^1}**@{},
   <0mm,0mm>*{};<-5mm,-6.9mm>*{^2}**@{},
   <0mm,0mm>*{};<4.5mm,-6.9mm>*{^{n\hspace{-0.5mm}-\hspace{-0.5mm}1}}**@{},
   <0mm,0mm>*{};<9.0mm,-6.9mm>*{^n}**@{},
 \end{xy}
& \mbox{for}\ |I_1|=\ldots=|I_m|=|J_1|=\ldots=|J_n|=1\\
0 & \mbox{otherwise}.
\Ea
\right.
\]
Indeed, one can define  $\cF^+$ and, respectively, $\cE\cK^+$ by their values on the generators,
\[
{\cF^+} (\mbox{resp.}\ \cE\cK^+) \left(
\xy
 <0mm,0mm>*{\mbox{$\xy *=<20mm,3mm>\txt{}*\frm{-}\endxy$}};<0mm,0mm>*{}**@{},
  <-10mm,1.5mm>*{};<-12mm,7mm>*{}**@{-},
  <-10mm,1.5mm>*{};<-11mm,7mm>*{}**@{-},
  <-10mm,1.5mm>*{};<-9.5mm,6mm>*{}**@{-},
  <-10mm,1.5mm>*{};<-8mm,7mm>*{}**@{-},
 <-10mm,1.5mm>*{};<-9.5mm,6.6mm>*{.\hspace{-0.4mm}.\hspace{-0.4mm}.}**@{},
 <0mm,0mm>*{};<-6.5mm,3.6mm>*{.\hspace{-0.1mm}.\hspace{-0.1mm}.}**@{},
  <-3mm,1.5mm>*{};<-5mm,7mm>*{}**@{-},
  <-3mm,1.5mm>*{};<-4mm,7mm>*{}**@{-},
  <-3mm,1.5mm>*{};<-2.5mm,6mm>*{}**@{-},
  <-3mm,1.5mm>*{};<-1mm,7mm>*{}**@{-},
 <-3mm,1.5mm>*{};<-2.5mm,6.6mm>*{.\hspace{-0.4mm}.\hspace{-0.4mm}.}**@{},
  <2mm,1.5mm>*{};<0mm,7mm>*{}**@{-},
  <2mm,1.5mm>*{};<1mm,7mm>*{}**@{-},
  <2mm,1.5mm>*{};<2.5mm,6mm>*{}**@{-},
  <2mm,1.5mm>*{};<4mm,7mm>*{}**@{-},
 <2mm,1.5mm>*{};<2.5mm,6.6mm>*{.\hspace{-0.4mm}.\hspace{-0.4mm}.}**@{},
 <0mm,0mm>*{};<6mm,3.6mm>*{.\hspace{-0.1mm}.\hspace{-0.1mm}.}**@{},
<10mm,1.5mm>*{};<8mm,7mm>*{}**@{-},
  <10mm,1.5mm>*{};<9mm,7mm>*{}**@{-},
  <10mm,1.5mm>*{};<10.5mm,6mm>*{}**@{-},
  <10mm,1.5mm>*{};<12mm,7mm>*{}**@{-},
 <10mm,1.5mm>*{};<10.5mm,6.6mm>*{.\hspace{-0.4mm}.\hspace{-0.4mm}.}**@{},
%
<0mm,0mm>*{};<-9.5mm,8.2mm>*{^{I_{ 1}}}**@{},
<0mm,0mm>*{};<-3mm,8.2mm>*{^{I_{ i}}}**@{},
<0mm,0mm>*{};<2mm,8.2mm>*{^{I_{ i+1}}}**@{},
<0mm,0mm>*{};<10mm,8.2mm>*{^{I_{ m}}}**@{},
<-10mm,-1.5mm>*{};<-12mm,-7mm>*{}**@{-},
  <-10mm,-1.5mm>*{};<-11mm,-7mm>*{}**@{-},
  <-10mm,-1.5mm>*{};<-9.5mm,-6mm>*{}**@{-},
  <-10mm,-1.5mm>*{};<-8mm,-7mm>*{}**@{-},
 <-10mm,-1.5mm>*{};<-9.5mm,-6.6mm>*{.\hspace{-0.4mm}.\hspace{-0.4mm}.}**@{},
 <0mm,0mm>*{};<-6.5mm,-3.6mm>*{.\hspace{-0.1mm}.\hspace{-0.1mm}.}**@{},
  <-3mm,-1.5mm>*{};<-5mm,-7mm>*{}**@{-},
  <-3mm,-1.5mm>*{};<-4mm,-7mm>*{}**@{-},
  <-3mm,-1.5mm>*{};<-2.5mm,-6mm>*{}**@{-},
  <-3mm,-1.5mm>*{};<-1mm,-7mm>*{}**@{-},
 <-3mm,-1.5mm>*{};<-2.5mm,-6.6mm>*{.\hspace{-0.4mm}.\hspace{-0.4mm}.}**@{},
  <2mm,-1.5mm>*{};<0mm,-7mm>*{}**@{-},
  <2mm,-1.5mm>*{};<1mm,-7mm>*{}**@{-},
  <2mm,-1.5mm>*{};<2.5mm,-6mm>*{}**@{-},
  <2mm,-1.5mm>*{};<4mm,-7mm>*{}**@{-},
 <2mm,-1.5mm>*{};<2.5mm,-6.6mm>*{.\hspace{-0.4mm}.\hspace{-0.4mm}.}**@{},
 <0mm,0mm>*{};<6mm,-3.6mm>*{.\hspace{-0.1mm}.\hspace{-0.1mm}.}**@{},
<10mm,-1.5mm>*{};<8mm,-7mm>*{}**@{-},
  <10mm,-1.5mm>*{};<9mm,-7mm>*{}**@{-},
  <10mm,-1.5mm>*{};<10.5mm,-6mm>*{}**@{-},
  <10mm,-1.5mm>*{};<12mm,-7mm>*{}**@{-},
 <10mm,-1.5mm>*{};<10.5mm,-6.6mm>*{.\hspace{-0.4mm}.\hspace{-0.4mm}.}**@{},
%
<0mm,0mm>*{};<-9.5mm,-9.2mm>*{^{J_{ 1}}}**@{},
<0mm,0mm>*{};<-3mm,-9.2mm>*{^{J_{j}}}**@{},
<0mm,0mm>*{};<2mm,-9.2mm>*{^{J_{ j+1}}}**@{},
<0mm,0mm>*{};<10mm,-9.2mm>*{^{J_{ n}}}**@{},
\endxy
\right)
=\left\{\Ba{cr}
\begin{xy}
 <0mm,-0.55mm>*{};<0mm,-3mm>*{}**@{-},
 <0mm,0.5mm>*{};<0mm,3mm>*{}**@{-},
 <0mm,0mm>*{\bullet};<0mm,0mm>*{}**@{},
 \end{xy} & \hspace{-9mm}\mbox{for}\ m=n=1, |I_1|=|J_1|=1,\vspace{2mm}\\
0  &\hspace{-9mm}\mbox{for}\ m=n=1,  |I_1|+|J_1|\geq 3,\vspace{2mm}\\
\hspace{-2mm}{\cF} (\mbox{resp.}\ \cE\cK)\left(
\xy
 <0mm,0mm>*{\mbox{$\xy *=<20mm,3mm>\txt{}*\frm{-}\endxy$}};<0mm,0mm>*{}**@{},
  <-10mm,1.5mm>*{};<-12mm,7mm>*{}**@{-},
  <-10mm,1.5mm>*{};<-11mm,7mm>*{}**@{-},
  <-10mm,1.5mm>*{};<-9.5mm,6mm>*{}**@{-},
  <-10mm,1.5mm>*{};<-8mm,7mm>*{}**@{-},
 <-10mm,1.5mm>*{};<-9.5mm,6.6mm>*{.\hspace{-0.4mm}.\hspace{-0.4mm}.}**@{},
 <0mm,0mm>*{};<-6.5mm,3.6mm>*{.\hspace{-0.1mm}.\hspace{-0.1mm}.}**@{},
  <-3mm,1.5mm>*{};<-5mm,7mm>*{}**@{-},
  <-3mm,1.5mm>*{};<-4mm,7mm>*{}**@{-},
  <-3mm,1.5mm>*{};<-2.5mm,6mm>*{}**@{-},
  <-3mm,1.5mm>*{};<-1mm,7mm>*{}**@{-},
 <-3mm,1.5mm>*{};<-2.5mm,6.6mm>*{.\hspace{-0.4mm}.\hspace{-0.4mm}.}**@{},
  <2mm,1.5mm>*{};<0mm,7mm>*{}**@{-},
  <2mm,1.5mm>*{};<1mm,7mm>*{}**@{-},
  <2mm,1.5mm>*{};<2.5mm,6mm>*{}**@{-},
  <2mm,1.5mm>*{};<4mm,7mm>*{}**@{-},
 <2mm,1.5mm>*{};<2.5mm,6.6mm>*{.\hspace{-0.4mm}.\hspace{-0.4mm}.}**@{},
 <0mm,0mm>*{};<6mm,3.6mm>*{.\hspace{-0.1mm}.\hspace{-0.1mm}.}**@{},
<10mm,1.5mm>*{};<8mm,7mm>*{}**@{-},
  <10mm,1.5mm>*{};<9mm,7mm>*{}**@{-},
  <10mm,1.5mm>*{};<10.5mm,6mm>*{}**@{-},
  <10mm,1.5mm>*{};<12mm,7mm>*{}**@{-},
 <10mm,1.5mm>*{};<10.5mm,6.6mm>*{.\hspace{-0.4mm}.\hspace{-0.4mm}.}**@{},
%
<0mm,0mm>*{};<-9.5mm,8.2mm>*{^{I_{ 1}}}**@{},
<0mm,0mm>*{};<-3mm,8.2mm>*{^{I_{ i}}}**@{},
<0mm,0mm>*{};<2mm,8.2mm>*{^{I_{ i+1}}}**@{},
<0mm,0mm>*{};<10mm,8.2mm>*{^{I_{ m}}}**@{},
<-10mm,-1.5mm>*{};<-12mm,-7mm>*{}**@{-},
  <-10mm,-1.5mm>*{};<-11mm,-7mm>*{}**@{-},
  <-10mm,-1.5mm>*{};<-9.5mm,-6mm>*{}**@{-},
  <-10mm,-1.5mm>*{};<-8mm,-7mm>*{}**@{-},
 <-10mm,-1.5mm>*{};<-9.5mm,-6.6mm>*{.\hspace{-0.4mm}.\hspace{-0.4mm}.}**@{},
 <0mm,0mm>*{};<-6.5mm,-3.6mm>*{.\hspace{-0.1mm}.\hspace{-0.1mm}.}**@{},
  <-3mm,-1.5mm>*{};<-5mm,-7mm>*{}**@{-},
  <-3mm,-1.5mm>*{};<-4mm,-7mm>*{}**@{-},
  <-3mm,-1.5mm>*{};<-2.5mm,-6mm>*{}**@{-},
  <-3mm,-1.5mm>*{};<-1mm,-7mm>*{}**@{-},
 <-3mm,-1.5mm>*{};<-2.5mm,-6.6mm>*{.\hspace{-0.4mm}.\hspace{-0.4mm}.}**@{},
  <2mm,-1.5mm>*{};<0mm,-7mm>*{}**@{-},
  <2mm,-1.5mm>*{};<1mm,-7mm>*{}**@{-},
  <2mm,-1.5mm>*{};<2.5mm,-6mm>*{}**@{-},
  <2mm,-1.5mm>*{};<4mm,-7mm>*{}**@{-},
 <2mm,-1.5mm>*{};<2.5mm,-6.6mm>*{.\hspace{-0.4mm}.\hspace{-0.4mm}.}**@{},
 <0mm,0mm>*{};<6mm,-3.6mm>*{.\hspace{-0.1mm}.\hspace{-0.1mm}.}**@{},
<10mm,-1.5mm>*{};<8mm,-7mm>*{}**@{-},
  <10mm,-1.5mm>*{};<9mm,-7mm>*{}**@{-},
  <10mm,-1.5mm>*{};<10.5mm,-6mm>*{}**@{-},
  <10mm,-1.5mm>*{};<12mm,-7mm>*{}**@{-},
 <10mm,-1.5mm>*{};<10.5mm,-6.6mm>*{.\hspace{-0.4mm}.\hspace{-0.4mm}.}**@{},
%
<0mm,0mm>*{};<-9.5mm,-9.2mm>*{^{J_{ 1}}}**@{},
<0mm,0mm>*{};<-3mm,-9.2mm>*{^{J_{j}}}**@{},
<0mm,0mm>*{};<2mm,-9.2mm>*{^{J_{ j+1}}}**@{},
<0mm,0mm>*{};<10mm,-9.2mm>*{^{J_{ n}}}**@{},
\endxy
\right)
& \mbox{for}\ m+n\geq 3\\
\Ea
\right.
\]
It is easy to check that the maps $\cE\cK^+$ and $\cF^+$ commute with the differentials.
In fact the above observation about $\cE\cK^+$ is equivalent to the extension of the Etingof-Kazhdan quantizations from Lie bialgebras to {\em differential}\, Lie bialgebras \cite{GH}. (Of course, one can construct the lifting  $\cF^+$ from the extended Etingof-Kazhdan quantization map $\cE\cK^+$
by an inductive procedure similar to the above construction of the lifting $\cF$ as the prop
$\DefQ^+$ is elementally cofibrant --- in the sense explained in \cite{Ma1} --- over its dg subprop
generated by corollas of non-negative homological degrees.)

\sip

The prop $\widehat{\Lieb}^+_\infty$ is completed and hence has no,
in general, meaningful representations in a graded vector space $V$, but does admit non-trivial continuous representations  in the topological vector space $V[[\hbar]]:=V\ot_\K \K[[\hbar]]$,
where $\hbar$ is a formal variable of homological degree zero. For example, for any representation
$$
\ga: {\Lieb}_\infty \lon \mathsf{End}_V
$$
i.e.\ for any strongly homotopy Lie bialgebra structure in dg vector space $(V,d)$, there is an associated
 well-defined continuous representation $\ga_\hbar$ of the topological dg prop given on the generators by,
$$
\ga_\hbar \left( \resizebox{10mm}{!}{\begin{xy}
 <0mm,0mm>*{\bullet};<0mm,0mm>*{}**@{},
 <0mm,0mm>*{};<-8mm,5mm>*{}**@{-},
 <0mm,0mm>*{};<-4.5mm,5mm>*{}**@{-},
 <0mm,0mm>*{};<-1mm,5mm>*{\ldots}**@{},
 <0mm,0mm>*{};<4.5mm,5mm>*{}**@{-},
 <0mm,0mm>*{};<8mm,5mm>*{}**@{-},
   <0mm,0mm>*{};<-8.5mm,5.5mm>*{^1}**@{},
   <0mm,0mm>*{};<-5mm,5.5mm>*{^2}**@{},
   <0mm,0mm>*{};<4.5mm,5.5mm>*{^{m\hspace{-0.5mm}-\hspace{-0.5mm}1}}**@{},
   <0mm,0mm>*{};<9.0mm,5.5mm>*{^m}**@{},
 <0mm,0mm>*{};<-8mm,-5mm>*{}**@{-},
 <0mm,0mm>*{};<-4.5mm,-5mm>*{}**@{-},
 <0mm,0mm>*{};<-1mm,-5mm>*{\ldots}**@{},
 <0mm,0mm>*{};<4.5mm,-5mm>*{}**@{-},
 <0mm,0mm>*{};<8mm,-5mm>*{}**@{-},
   <0mm,0mm>*{};<-8.5mm,-6.9mm>*{^1}**@{},
   <0mm,0mm>*{};<-5mm,-6.9mm>*{^2}**@{},
   <0mm,0mm>*{};<4.5mm,-6.9mm>*{^{n\hspace{-0.5mm}-\hspace{-0.5mm}1}}**@{},
   <0mm,0mm>*{};<9.0mm,-6.9mm>*{^n}**@{},
 \end{xy}}\right)
 = \hbar^{m+n-2} \ga \left( \resizebox{10mm}{!}{\begin{xy}
 <0mm,0mm>*{\bullet};<0mm,0mm>*{}**@{},
 <0mm,0mm>*{};<-8mm,5mm>*{}**@{-},
 <0mm,0mm>*{};<-4.5mm,5mm>*{}**@{-},
 <0mm,0mm>*{};<-1mm,5mm>*{\ldots}**@{},
 <0mm,0mm>*{};<4.5mm,5mm>*{}**@{-},
 <0mm,0mm>*{};<8mm,5mm>*{}**@{-},
   <0mm,0mm>*{};<-8.5mm,5.5mm>*{^1}**@{},
   <0mm,0mm>*{};<-5mm,5.5mm>*{^2}**@{},
   <0mm,0mm>*{};<4.5mm,5.5mm>*{^{m\hspace{-0.5mm}-\hspace{-0.5mm}1}}**@{},
   <0mm,0mm>*{};<9.0mm,5.5mm>*{^m}**@{},
 <0mm,0mm>*{};<-8mm,-5mm>*{}**@{-},
 <0mm,0mm>*{};<-4.5mm,-5mm>*{}**@{-},
 <0mm,0mm>*{};<-1mm,-5mm>*{\ldots}**@{},
 <0mm,0mm>*{};<4.5mm,-5mm>*{}**@{-},
 <0mm,0mm>*{};<8mm,-5mm>*{}**@{-},
   <0mm,0mm>*{};<-8.5mm,-6.9mm>*{^1}**@{},
   <0mm,0mm>*{};<-5mm,-6.9mm>*{^2}**@{},
   <0mm,0mm>*{};<4.5mm,-6.9mm>*{^{n\hspace{-0.5mm}-\hspace{-0.5mm}1}}**@{},
   <0mm,0mm>*{};<9.0mm,-6.9mm>*{^n}**@{},
 \end{xy}}\right), \ \ \ \ m,n\geq 1.
 $$
The morphism (\ref{1: map of props cF})
 induces a continuous $L_\infty$ quasi-morphism,
 $$
 F^{\hbar}: \Def_{cont}(\widehat{\Lieb}^{_+}_\infty\stackrel{0}{\rar} \mathsf{End}_V[[\hbar]])
 \lon
 \Def_{cont}({\DefQ}^{_+}\stackrel{\rho_0}{\rar} \mathsf{End}_V[[\hbar]])\simeq
 \Def(\sB_\infty^+ \stackrel{\rho_0}{\rar} \mathsf{End^{poly}_{\f_{\mathit V}}}[[\hbar]])
 $$
that is,
$$
 F^{\hbar}: \fl_V[[\hbar]] \lon \mathsf{poly}(\f_V,\f_V)[[\hbar]].
$$
There is a morphism of Lie algebras,
\Beq
\Ba{rccc}\label{1: l_V to l_V[[h]]}
\ii: & \fl_V & \lon & \fl_V[[\hbar]]\\
    & m\in \odot^m(V[-1])\ot \odot^n(V[-1]) &\lon & \hbar^{m+n-2} m
\Ea
\Eeq
such that the composition
$$
F_k^\hbar \circ \wedge^k \ii: \wedge^k ( \fl_V[2]) \lon \mathsf{poly}(\f_V,\f_V)[[\hbar]]
$$
takes values in  the subspace $\mathsf{poly}(\f_V,\f_V)[\hbar]$. Then the family of maps
$$
F=\{F_k:=  F_k^\hbar|_{\hbar=1} \circ \wedge^k\ii\}
$$
gives us finally the required quasi-isomorphism (\ref{1: Formality map F poly}), and hence
(\ref{1:F}). The second proof is completed.

\bip

\bip

{\em Acknowledgement}. {\small This paper is a revised version of the preprint  arXiv:math/0612431
 written by the author during his visit to the MPIM in Bonn in 2006.
 We changed the title, added new references as well as second proofs of some important statements, and clarified some subtleties about completions of the props used in the proofs. We apologize to all authors who cited the previous preprint for such a late revision.

}

\bip

\bip

\bip

\def\cprime{$'$}

  \end{document}